\documentstyle[12pt]{article}

\newlength{\abstractwidth}
\renewcommand{\thefootnote}{\fnsymbol{footnote}}
\renewcommand{\thanks}[1]{\footnote{#1}} 
\newcommand{\starttext}{ \setcounter{footnote}{0}
\renewcommand{\thefootnote}{\arabic{footnote}}}

\newcommand{\be}{\begin{equation}}
\newcommand{\bea}{\begin{eqnarray}}
\newcommand{\eea}{\end{eqnarray}} 
\newcommand{\ee}{\end{equation}}

\def\ba{\begin{eqnarray}}
\def\ea{\end{eqnarray}}

\def\o{\omega}
\def\Re{{\rm Re}}

\def\det{{\rm det}}

\def\log{\,{\rm log}\,}
\def\exp{\,{\rm exp}\,}

\def\o{\omega}

\def\e{\varepsilon}

\def\o{\omega}

\def\na{\nabla}
\def\p{\partial}

\def\ddb{{\partial\bar\partial}}

\def\na{{\nabla}}

\def\[{{\bf [}}
\def\]{{\bf ]}}

\begin{document}
\starttext \baselineskip=15pt \setcounter{footnote}{0}
\newtheorem{theorem}{Theorem}
\newtheorem{lemma}{Lemma}
\newtheorem{definition}{Definition}
\newtheorem{proposition}{Proposition}
\newtheorem{corollary}{Corollary}

\begin{center}
{\Large \bf
THE ANOMALY FLOW AND THE FU-YAU EQUATION
\footnote{Work supported in part by the National Science Foundation Grants DMS-12-66033 and DMS-16-05968. Key words: Hull-Strominger systems, Goldstein-Prokushkin fibration, slope parameter $\alpha'$, conformally balanced, torsion constraints, curvature, Moser iteration, $C^k$ estimates, exponential convergence.}}

\bigskip

{\large Duong H. Phong, Sebastien Picard, and Xiangwen Zhang} \\
\medskip
\begin{abstract} 

{\small The Anomaly flow is shown to converge on toric fibrations with the Fu-Yau ansatz, for both positive and negative values of the slope parameter $\alpha'$. This implies both results of Fu and Yau on the existence of solutions for Hull-Strominger systems, which they proved using different methods depending on the sign of $\alpha'$. It is also the first case 
where the Anomaly flow can even be shown to exist for all time. This is in itself
remarkable from the point of view of the theory of fully nonlinear partial differential equations, as the elliptic terms in the flow are not concave.
}
\end{abstract}
\end{center}

\section{Introduction}
\setcounter{equation}{0}

The Hull-Strominger system \cite{H,S} is a system of equations for supersymmetric compactifications of the heterotic string, which is less restrictive than the Ricci-flat and K\"ahler conditions originally proposed by Candelas, Horowitz, Strominger, and Witten \cite{CHSW}. More specifically, let $Y$ be a compact $3$-fold with a nowhere vanishing holomorphic $(3,0)$-form $\Omega_Y$ and a vector bundle $E\to Y$. Then the Hull-Strominger system is
the following system of equations for a Hermitian metric $\chi$ on $Y$
and a Hermitian metric $H$ on $E$,
\bea
\label{HYM}
&&
F^{2,0}=F^{0,2}=0,\quad \chi^2\wedge F^{1,1}=0\\
\label{GS}
&&i\p\bar\p\chi-{\alpha'\over 4}
{\rm Tr}(Rm\wedge Rm -F\wedge F)=0\\
\label{CB}
&&d(\|\Omega_Y\|_\chi\chi^2)=0.
\eea
Here $\alpha'$ is a constant parameter called the slope, $\|\Omega_Y\|_\chi$ is the norm of $\Omega_Y$ with respect to $\chi$ defined by
$\|\Omega_Y\|_\chi^2=i\Omega_Y\wedge\bar\Omega_Y\, \chi^{-3}$, and $Rm$ and $F$ are respectively the curvatures of the Chern unitary connections of $\chi$ and $H$, viewed as $(1,1)$-forms valued in the bundles of endomorphisms of  $T^{1,0}(Y)$ and $E$
\footnote{The equation (\ref{CB}) was originally written in \cite{H,S} as $d^\dagger\chi=i(\p-\bar\p)\log\|\Omega\|_\chi$. That (\ref{CB}) is an equivalent formulation is an important insight due to J. Li and S.T. Yau \cite{LY}.}. For fixed $\chi$, the first equation (\ref{HYM}) is just the well-known Hermitian-Yang-Mills equation. The essential novelty in the Hull-Strominger system, from the point of view of both non-K\"ahler geometry and non-linear partial differential equations, resides rather in the last two equations (\ref{GS}) and (\ref{CB}). The equation (\ref{CB}) says that the metric $\chi$ is not required to be K\"ahler, but conformally balanced, and the equation (\ref{GS}) defines a new curvature condition, which differs markedly from more familiar conditions such as  Einstein since it is quadratic in the curvature
tensor. By now, many special solutions of the Hull-Strominger system have been found, both in the physics (see e.g. \cite{BBFT,CI,DRS}) and the mathematics literature (see e.g. \cite{AG1, AG2, Fe1, FeY, FIUV1, FIUV2, FTY, FY1,FY2,
Gr, LY, OUV, UV}). The main goal in the present paper is rather to develop PDE techniques towards an eventual general solution.

\smallskip

A first major difficulty in the Hull-Strominger system, symptomatic of non-K\"ahler geometry, is to implement the conformally balanced condition (\ref{CB}) in the absence of an analogue of the $\p\bar\p$-lemma. While balanced metrics can be produced by many Ans\"atze, none seems more natural than the others, and all lead to unwieldy expressions for the equation (\ref{GS}). This is why it was proposed by the authors in \cite{PPZ} to bypass completely the choice of an Ansatz by considering instead the flow
\bea
\label{AFlow}
&&
H^{-1} \, \p_tH= - \Lambda_\chi F(H)\\
&&
\label{GS1}
\p_t(\|\Omega_Y\|_\chi\chi^2)
=
i\p\bar\p\chi- {\alpha' \over 4} {\rm Tr}(Rm(\chi)\wedge Rm(\chi)-F(H)\wedge F(H))
\eea
starting from a metric $H(0)$ on $E$, and a metric $\chi(0)$ on $Y$ which is conformally balanced. Here
$\Lambda_\chi \psi=\chi^2\wedge\psi\,\chi^{-3}$ is the Hodge operator on $(1,1)$-forms $\psi$. 
The point of this flow is that, by Chern-Weil theory, the right hand side of (\ref{GS1}) is a closed $(2,2)$-form, so if the initial metric $\chi(0)$ is conformally balanced, then the metric $\chi(t)$ will remain conformally balanced for all $t$. It suffices then to determine whether the flow exists for all time and converges. Flows of the form (\ref{AFlow}) have been called Anomaly flows in \cite{PPZ}, in recognition of the fact that the equation (\ref{GS}) originates from the famous Green-Schwarz anomaly cancellation mechanism \cite{GS} required for the consistency of superstring theories. 

\smallskip
The flow (\ref{AFlow}) has been shown in \cite{PPZ} to be weakly parabolic when $|\alpha'Rm(\chi)|$ is small, which implies its short-time existence.
However, for any given elliptic system, there are many parabolic flows with it as stationary point, so the true test of whether a particular parabolic flow is the right one is its long-time existence and convergence. In the general theory of 
fully non-linear PDE's, it is customary to select the parabolic flow by some desirable properties, such as, in the case of scalar equations, concavity in the second derivatives \cite{CW, GuanSS, L}. In the present case,
there is no further flexibility, as the particular flow (\ref{AFlow}) provides the only known way of implementing the conformally balanced condition (\ref{CB}) without appealing to any particular Ansatz. In effect, the geometric constraint of the metric being conformally balanced has excluded the customary desirable analytic properties for flows, and it is vital at this stage to develop some new analytic tools. 
For this, we shall consider the Anomaly flow on the model case of Calabi-Eckmann fibrations. This case retains enough features of the general case to provide a valuable guide in the future (see \S \ref{comparison-2d-3d} below). It is also the case where J.X. Fu and S.T. Yau \cite{FY1, FY2} found, by a particularly difficult and delicate analysis, the first non-perturbative, non-K\"ahler solution of the Hull-Strominger system. The precise $C^0$ estimate that they obtained has a great influence on the present work. On the other hand, we shall see that the Anomaly flow requires very different $C^1$, $C^2$, and $C^{2,\alpha}$ estimates, which actually sharpen the elliptic estimates in many ways. In particular, they define a new region, narrower than the space of positive Hermitian forms, which is preserved by the flow and where the unknown metric ultimately belongs. We now describe more precisely our results.

\bigskip
Let $(X,\hat\o)$ be a Calabi-Yau surface, equipped with a nowhere vanishing holomorphic $(2, 0)$-form $\Omega$ normalized to satisfy $\|\Omega\|_{\hat\o}=1$ and two harmonic form $\o_1,\o_2\in H^{1,1}(X,{\bf Z})$. Building on an earlier construction of Calabi and Eckmann \cite{CE}, Goldstein and Prokushkin \cite{GP} have shown how to associate to this data a toric fibration $\pi:Y\to X$, equipped with a $(1,0)$-form $\theta$ with the property that $\Omega_Y=\Omega\wedge\theta$ is holomorphic and non-vanishing on $Y$, and $\o_u=\pi^*(e^u\hat\o)+i\theta\wedge\bar \theta$ is a conformally balanced metric on $Y$ for all $u\in C^\infty(X)$ (see \S \ref{gp-fibration-background} for more details). We shall prove

\begin{theorem}
\label{main1}
Let $\pi:Y\to X$ be a Goldstein-Prokushkin fibration
over a Calabi-Yau surface $(X,\hat\o)$ with Ricci-flat metric $\hat\o$ as above, and $E_X\to (X,\hat\o)$ a stable vector bundle with zero slope and Hermitian-Yang-Mills metric $H_X$. Assume that the cohomological condition
$\int_X \mu = 0$ is satisfied, where $\mu$ is defined in (\ref{mu}) below. 
Set $E=\pi^*(E_X)$ and consider 
the Anomaly flow (\ref{AFlow}, \ref{GS1}) on $(Y,E)$, with initial data 
$H(0)=\pi^*(H_X)$, and $\chi(0)=\pi^*(M\hat\o)+i\theta\wedge
\bar\theta$, where $M$ is a positive constant. Then there exists $M_0 \gg 1$ such that, for all $M\geq M_0$, the flow exists for all time, and converges smoothly to 
a solution $(\chi_\infty,\pi^*(H_X))$ of the Hull-Strominger system for $(Y,E)$.
\end{theorem}

In particular, the theorem recaptures at one stroke the results of Fu and Yau in \cite{FY1,FY2}, where they proved the existence of solutions of the Hull-Strominger system on Goldstein-Prokushkin fibrations satisfying the cohomological condition $\int_X\mu=0$, when $\alpha'>0$ and $\alpha'<0$ respectively.

\medskip
Restricted to a Goldstein-Prokushkin fibration, the Anomaly flow becomes equivalent to the following flow for a metric $\o=ig_{\bar kj}dz^j\wedge d\bar z^k$ on a Calabi-Yau surface $X$, equipped with a nowhere vanishing holomorphic $(2,0)$-form $\Omega$,
\bea \label{FY_parabolic}
\p_t \o
=
-{1\over 2 \|\Omega\|_\o}
\bigg( {R\over 2}- |T|^2 - {\alpha' \over 4} \sigma_2(i{\rm Ric}_\o)  + 2 \alpha' {i\p\bar\p(\|\Omega\|_\o\rho)\over\o^2}
 - 2{\mu\over\o^2} \bigg) \,\o
\eea
where $\sigma_2(\Phi)=\Phi\wedge\Phi\,\o^{-2}$ is the usual determinant of a real $(1,1)$-form $\Phi$, relative to the metric $\o$. The expression $|T|^2$ is the norm of the torsion of $\o$ defined in (\ref{T_omega}) below.
Thus the theorem that we shall actually prove is the following

\begin{theorem}
\label{main2}
Let $(X,\hat\o)$ be a Calabi-Yau surface, equipped with a Ricci-flat metric $\hat\o$ and a nowhere vanishing holomorphic $(2,0)$-form $\Omega$ normalized by $\|\Omega\|_{\hat\o}=1$. Let $\alpha'$ be a non-zero real number, and let $\rho$ and $\mu$ be smooth real $(1,1)$ and $(2,2)$-forms respectively, with $\mu$ satisfying the integrability condition
\bea
\int_X \mu=0.
\eea
Consider the flow (\ref{FY_parabolic}), with an initial metric given by $\o(0)=M \, \hat\o$, where $M$ is a constant. Then there exists $M_0$ large enough so that, for all $M\geq M_0$, the flow
(\ref{FY_parabolic}) exists for all time, and converges exponentially fast to a metric $\o_\infty$ satisfying the Fu-Yau equation
\bea
\label{FY_elliptic}
i\p\bar\p (\o_\infty-\alpha'\|\Omega\|_{\o_\infty} \, \rho)
-
{\alpha' \over 8} {\rm Ric}_{\o_\infty} \wedge {\rm Ric}_{\o_\infty} +\mu=0, 
\eea
and the normalization $\int_X \| \Omega \|_{\o_\infty} {\o_\infty^2 \over 2!} = M.$
\end{theorem}

\medskip

Rewriting the flow as a flow of the conformal factor $u$, we can discuss now more precisely its features in the context of the general theory of non-linear parabolic PDE's. The flow (\ref{FY_parabolic}) can be expressed as
\bea
\label{FY_parabolic_u}
\p_t u = {1 \over 2} \bigg({\Delta_{\hat\o}} u + \alpha'  e^{-u} \hat{\sigma}_2(i \ddb u) - 2 \alpha' e^{-u} {i \ddb (e^{-u} \rho) \over \hat{\o}^2}+ |Du|^2_{\hat{\o}}  + e^{-u} \tilde{\mu} \bigg)
\eea
where $\tilde\mu=2\mu\,\hat\o^{-2}$ is a time-independent scalar function, and both the Laplacian $\Delta_{\hat\o}$ and the determinant $\hat\sigma_2$ are written with respect to the fixed metric $\hat\o$.
Setting the right hand side to $0$ gives the equation solved by Fu and Yau \cite{FY1,FY2}, so the Anomaly flow is indeed a parabolic version of the Fu-Yau equation. Moreover, the equation (\ref{FY_parabolic_u}) can be rewritten in the form
\bea
\label{FY_parabolic_u1}
2\alpha' e^u \p_t u = {\left(e^u \hat \o + e^{-u}\rho + \alpha' i\p\bar\p u\right)^2\over \hat\o^2} + w(\mu, \rho, u, Du ),
\eea
and it is a parabolic complex Monge-Amp\`ere type equation. However, unlike the K\"ahler-Ricci flow where the elliptic term is $\log\det(g_{i\bar j} + u_{i\bar j})$, the equation (\ref{FY_parabolic_u1}) has none of the desirable concavity properties of elliptic and parabolic equations. In particular, none of the techniques used in \cite{FY1,FY2} for the elliptic case (as well as the ones for the more general equations in \cite{PPZ2,PPZ3}) can be adapted here besides the Moser iteration technique for the $C^0$ estimate. We shall see below that the proof of Theorem \ref{main2} relies instead, in an essential manner, on the geometric formulation of (\ref{FY_parabolic}), and that the estimates are obtained using the metric which evolves with the flow. 

\medskip

 Various geometric flows have been studied in non-K\"ahler complex geometry (see e.g. \cite{Chu,Gi,SW,ST,Su,TW1,TW2} and references therein). The main difficulty in studying (\ref{FY_parabolic}) is that it is quadratic in the Ricci curvature. This creates substantial problems in applying known techniques to try and obtain estimates on the torsion, curvature, and derivatives of curvature. To overcome these issues, we first start the flow with a metric with vanishing Ricci curvature and torsion, and the objective is to show that for suitably large normalization of the initial metric, we can prevent the terms which are nonlinear in curvature from growing too large and dominating the behavior of the flow. Proposition \ref{all_estimates} in \S 7 shows exactly how the estimates on Ricci curvature and torsion depend on the normalization of the initial metric. More precisely, the metric $\omega(t)$ will be proved to belong at all times in the following set
 \bea
  C^{-1}M \, \hat\omega \leq \omega(t) \leq CM\, \hat\omega, \quad |T(\omega(t))| \leq C M^{-1/2}, \quad
|\alpha'Ric_{\o(t)}|\leq C M^{-1/2}
 \eea
 for a constant $C$ depending only on the geometric data defining the Goldstein-Prokushkin fibration, the Hermitian-Yang-Mills metric $H_X$, and the slope $\alpha'$, but independent of the normalization constant $M$. In particular, it is a set much narrower and more specific than the cone of metrics on $X$. 
 
\medskip
The normalization of the initial metric actually sets a scale for the problem, and Proposition \ref{all_estimates} is an example of estimates with scales. The key underlying strategy of the present paper can be viewed as the use of estimates with scales to tame higher powers of the curvature tensor. This strategy appears both flexible and powerful: since the initial posting online in 2016 of the original version of the present paper, we have applied it successfully in our paper ``Fu-Yau Hessian equations", arXiv:1801.09842, in fact with the same test functions used in the proof of Proposition \ref{all_estimates}, to solve the Fu-Yau equation and its Hessian generalizations in all dimensions, for both signs of the parameter $\alpha'$. See also the paper ``The Fu-Yau equation in higher
dimensions" by J. Chu, L. Huang, and X.H. Zhu, arXiv:1801.09351.

\medskip
The paper is organized as follows. In Section \S 2, we provide the background on Goldstein-Prokushkin fibrations, and show how to reduce the Anomaly flow in this case to the flow (\ref{FY_parabolic}) for metrics on a Calabi-Yau surface. Sections \S 3-\S 6 are devoted to successive estimates: the uniform boundedness of the metrics $\o$ in \S 3, the estimates for the torsion in Section \S 4, the estimates for the curvature, and higher order derivatives of both the torsion and the curvature in Sections \S 5-\S 6. Finally, long time existence is shown in Section \S 7 and the convergence of the flow is proved in Section \S 8.

\

\section{Anomaly flows on Goldstein-Prokushkin fibrations}
\setcounter{equation}{0}

\subsection{The Goldstein-Prokushkin fibration} \label{gp-fibration-background}

We would like to restrict the Anomaly flow from a general $3$-fold $Y$ to the special case of a Goldstein-Prokushkin fibration $\pi:Y\to X$. 
We begin by recalling the basic properties of Goldstein-Prokushkin fibrations that we need. 

\smallskip

Let $(X,\hat\o)$ be a compact Calabi-Yau manifold of dimension $2$, with $\hat\o$ a Ricci-flat K\"ahler metric, and $\Omega$ a nowhere vanishing holomorphic $(2,0)$-form, normalized so that 
\bea
1=\|\Omega\|_{\hat\o}^2= \Omega\wedge \overline{\Omega}\,\hat\o^{-2}.
\eea
Let $\o_1, \o_2 \in 2\pi\,H^2(X,{\bf Z})$ be two $(1,1)$-forms such that $\o_1 \wedge \hat\o = \o_2 \wedge \hat\o=0$. From this data, Goldstein and Prokushkin \cite{GP} construct a compact $3$-fold $Y$ which is a toric fibration $\pi: Y \rightarrow X$ over $X$ equipped with a $(1,0)$ form $\theta$ on $Y$ 
satisfying
\bea
\label{theta}
\bar\p\theta=\pi^*(\o_1+i\o_2)\qquad \p\theta=0.
\eea
Furthermore, the $(3,0)$-form
\bea
\Omega_Y=\sqrt 3\,\Omega\wedge \theta
\eea
is holomorphic and nowhere vanishing, and 
the $(1,1)$-form 
\bea
\chi_0=\pi^*(\hat\o)+i\theta\wedge\bar\theta
\eea
is positive-definite on $Y$. Observe that
\bea
i\Omega_Y\wedge\overline{\Omega_Y}
=
3 \Omega\wedge\overline{\Omega}\wedge\,i\theta\wedge\bar\theta
=
\|\Omega\|_{\hat\o}^2 (3\hat\o^2\wedge\,i\theta\wedge\bar\theta)
=
\|\Omega\|_{\hat\o}^2 \,\chi_0^3.
\eea
Thus, defining the norm $\|\Omega_Y\|_{\chi}$ of the holomorphic form $\Omega_Y$ on $Y$ with respect to a metric $\chi$ as $\|\Omega_Y \|_\chi^2=i\Omega_Y \wedge\bar\Omega_Y \wedge\chi^{-3}$, we have $\|\Omega_Y\|_{\chi_0}=\|\Omega\|_{\hat\o}=1$. Consequently,
\bea
\|\Omega_Y\|_{\chi_0}\chi_0^2
=
\hat \o^2+2i\hat\o \wedge \theta\wedge\bar\theta.
\eea
This implies that $d(\|\Omega_Y\|_{\chi_0}\chi_0^2)=0$ by (\ref{theta}) and the fact that $\o_1$ and $\o_2$ wedged with $\hat{\o}$ gives zero. Thus $\chi_0$ is a conformally balanced metric on $Y$.

\smallskip
More generally, for any smooth function $u$ on $Y$, introduce the following metrics $\o_u$ and $\chi_u$ on the manifolds $X$ and $Y$ respectively,
\bea
\o_u=e^u\hat\o,
\qquad
\chi_u=\pi^*(e^u\hat\o)+i\theta\wedge\bar\theta.
\eea
Then the same arguments that we just used show that
\bea
\|\Omega_Y\|_{\chi_u}=\|\Omega\|_{\o_u}=e^{-u},
\eea
and furthermore,
\bea \label{Omega-chi_u}
\|\Omega_Y\|_{\chi_u}\chi_u^2
=
\|\Omega\|_{\o_u}\o_u^2
+
2i\hat\o\wedge \theta\wedge\bar\theta.
\eea
This shows that $d(\|\Omega_Y\|_{\chi_u}\chi_u^2)=0$  since $d (\|\Omega\|_{\o_u}\o_u^2)$ is the pull-back of a differential form of rank $5$ defined on the $4$-dimensional manifold $X$, and $2i d (\hat\o\wedge \theta\wedge\bar\theta)$ is zero as before in view of (\ref{theta}). It follows that the metric $\chi_u$ is a conformally balanced metric on $Y$ for any choice of $u$.

\subsection{The Fu-Yau Ansatz}

In \cite{FY1, FY2}, Fu and Yau obtain a solution of the Hull-Strominger system in the following manner. Let $\pi:Y\to X$ be a Goldstein-Prokushkin fibration, constructed as described above from a Calabi-Yau surface $(X,\hat\o)$, equipped with two integer-valued harmonic $(1,1)$-forms $\o_1/(2\pi)$ and $\o_2/(2\pi)$.

\medskip

Let $E_X\to X$ be a stable holomorphic vector bundle over $X$, with slope
$\int_X c_1(E_X)\wedge \hat\o=0$. Then by the Donaldson-Uhlenbeck-Yau \cite{D, UY} theorem, $E_X$ admits a metric $H_X$ with respect to $\hat\o$ satisfying the Hermitian-Yang-Mills equation $\hat\o\wedge F(H_X)=0$.  Let $E=\pi^*(E_X)\to Y$ be the pull-back bundle over $Y$, and let $H=\pi^*(H_X)$. Since 
\bea
\chi_u^2\wedge F(H)=\pi^*(e^u\hat\o\wedge F(H_X))\wedge (e^u\hat\o+2i\theta\wedge\bar\theta)=0,
\eea 
it follows that $H$ is Hermitian-Yang-Mills with respect to $\chi_u$, for any $u$. Now recall that $\chi_u$ is conformally balanced for any $u$. This means that, if we look for a solution of the Strominger system under the form $(Y,E)$, equipped with the metrics $\chi_u$ and $H$, then the only equation which is left to solve is the Green-Schwarz anomaly cancellation equation (\ref{GS}), with the scalar function $u$ defining the metric $\chi_u$ as the unknown.

\medskip
The key property that allows this approach to work is that, for metrics of the form $\chi_u$ in a Goldstein-Prokushkin fibration, the equation on $Y$
\bea
i\p\bar\p\chi_u- {\alpha' \over 4} {\rm Tr}(Rm(\chi_u)\wedge Rm(\chi_u)
-F\wedge F)=0
\eea
descends to an equation on the base $X$. This was established by Fu and Yau \cite{FY1} and a summary of their results is as follows.

\medskip
First, the term $i\p\bar\p\chi_u$ is readily worked out, using the properties (\ref{theta}) of the form $\theta$,
\bea
i\p\bar\p\chi_u
=
i\p\bar\p \o_u
-\bar\p\theta\wedge\p\bar\theta.
\eea
Next, the quadratic term in the curvature tensor can be worked out to be (Proposition 8 in \cite{FY1})
\bea
{\rm Tr}(Rm(\chi_u)\wedge Rm(\chi_u))
=
{\rm Tr}(Rm(\o_u)\wedge Rm(\o_u))
+\p\bar\p(\|\Omega\|_{\o_u}{\rm Tr}(\bar\p B\wedge\p B^*\cdot \hat\o^{-1})).
\eea
Here $B$ is a $(1,0)$-form depending on the data $(\hat\o,\o_1,\o_2)$, 
which is only locally defined. However the full expression
${\rm Tr}(\bar\p B\wedge\p B^*\cdot \hat\o^{-1})$ is not only globally well-defined on $Y$, but it is the pull-back of a globally defined real $(1,1)$-form $\rho$ on $X$,
\bea
{1 \over 4} {\rm Tr}(\bar\p B\wedge\p B^*\cdot \hat\o^{-1})= \pi^*(\rho).
\eea
On the other hand, from $\o_u=e^u\hat\o$, it follows that
\bea
Rm(\o_u)
=
-\p\bar\p u\otimes I
+
Rm(\hat\o)
\eea
and hence, in view of the fact that the metric $\hat{\o}$ is Ricci-flat,
\bea
{\rm Tr}(Rm(\o_u)\wedge Rm(\o_u))
=
{\rm Tr}(Rm(\hat\o)\wedge Rm(\hat\o))
+2\p\bar\p u\wedge \p\bar\p u.
\eea
Altogether, the Green-Schwarz anomaly cancellation equation can be written as the following equation on the base manifold $X$,
\bea
\label{FY_elliptic}
0&=& i\p\bar\p(\o_u-\alpha'\|\Omega\|_{\o_u}\rho)
-
{\alpha' \over 2} (\p\bar\p u)\wedge (\p\bar\p u)+\mu
\eea
where we have set
\bea
\label{mu}
\mu =-\bar\p\theta\wedge \p\bar\theta
-
{\alpha' \over 4} {\rm Tr}(Rm(\hat\o) \wedge Rm(\hat\o))+{\alpha' \over 4}{\rm Tr}(F(H_X)\wedge F(H_X)).
\eea
which is a well-defined $(2,2)$-form on $X$. The equation (\ref{FY_elliptic}) is the Fu-Yau equation. Clearly, a necessary condition for the existence of solutions is
\bea
\int_X \mu=0.
\eea 
This condition was shown to be sufficient by Fu and Yau in \cite{FY1} for $\alpha'>0$ and in \cite{FY2} for $\alpha'<0$. Examples of fibrations $\pi: Y \rightarrow X$ and vector bundles $E \rightarrow X$ satisfying $\int_X \mu=0$ are exhibited in \cite{FY1,FY2}.

\subsection{Reduction of the Anomaly flow by the Fu-Yau Ansatz}

We consider now the Anomaly flow (\ref{GS1}) on a Goldstein-Prokushkin fibration $\pi:Y\to X$, equipped with the holomorphic $(3,0)$-form $\Omega_Y$ and restricted to metrics of the form $\chi=\chi_u$. Recall that $F=F(\pi^*(H_X))$ is fixed, and that $\Lambda_\chi F=0$. Thus the metric $H_X$ remains fixed along the flow, and we need only concentrate on the flow for the metric $\chi_u$.

\medskip
We work out both sides of the flow (\ref{GS1}) in this setting. Recall that we have set $\o_u=e^u\hat\o$ which is a metric on $X$. From the earlier equation (\ref{Omega-chi_u}) and the fact that the term $\hat\o\wedge\theta\wedge \bar\theta$ is time-independent, it follows at once that
\bea
\p_t(\|\Omega_Y\|_{\chi_u}\chi_u^2)
=
\p_t(\|\Omega\|_{\o_u}\o_u^2).
\eea
On the other hand, the same formulas derived by Fu and Yau \cite{FY1} for their reduction of the anomaly equation on $Y$ to an equation on $X$ and which we described in the previous section give
\bea
&&
i\p\bar\p\chi_u
-
{\alpha' \over 4} {\rm Tr}(Rm(\chi_u)\wedge Rm(\chi_u)-F\wedge F)
\nonumber\\
&&
\qquad
=
i\p\bar\p( \o_u- \alpha' \|\Omega\|_{\o_u}\rho)
-
{\alpha' \over 2} \p\bar\p u\wedge\p\bar\p u
+
\mu
\eea
with $\mu$ the $(2,2)$-form defined by (\ref{mu}). Thus the Anomaly flow for Goldstein-Prokushkin fibrations is equivalent to the flow for metrics on $X$ given by
\bea \label{flow_forms}
\p_t(\|\Omega\|_{\o_u}\o_u^2)
=
i\p\bar\p( \o_u-\alpha' \|\Omega\|_{\o_u}\rho)
-
{\alpha' \over 2} \p\bar\p u\wedge\p\bar\p u
+
\mu.
\eea
Since we wish to apply techniques of geometric flows, it is useful to re-express the flow entirely in terms of curvature. If we denote by ${\rm Ric}_{\o_u}$ the Chern-Ricci tensor of $\o_u$, we have
\bea
{\rm Ric}_{\o_u}=-2\p\bar\p u
\eea
since the metric $\hat\o$ is Ricci-flat. Thus the Anomaly flow can be rewritten as
the following flow of metrics on $X$,
\bea
\label{Flow1}
\p_t(\|\Omega\|_{\o}\o^2)
=
i\p\bar\p( \o- \alpha' \|\Omega\|_{\o}\rho)
-
{\alpha' \over 8} {\rm Ric}_\o \wedge {\rm Ric}_\o
+
\mu
\eea
which we can take now as our starting point. Here we have suppressed the subindex $u$ in $\o_u$.

\medskip
A technical issue in Anomaly flows is that they are formulated in terms of flows for $\|\Omega\|_{\o}\,\o^2$, and not of $\o$ itself. This issue was addressed in all generality in \cite{PPZ1} in dimension $3$. For the above Anomaly flow on the surface $X$ arising from the Goldstein-Prokushkin fibration, the metric $\o$ is already characterized by its volume form, and we can proceed more directly as follows.

\medskip
First, using the two-dimensional identity $2\p_t\o\wedge\o=(\p_t\log\o^2)\o^2$ and the fact that $\o^2 = \| \Omega \|_{\o}^{-2} \hat{\o}^2$, we can rewrite the left hand side as
\bea \label{p_tOmega-log}
\p_t(\|\Omega\|_\o\o^2)
=
\|\Omega\|_\o(\p_t\log\|\Omega\|_\o \o^2+2\p_t\o\wedge\o)
=
- \|\Omega\|_\o(\p_t\log \| \Omega \|_\o )\o^2.
\eea

\smallskip
Next, we work out the right hand side more explicitly. Following \cite{PPZ1}, we define the torsion $T(\o)={1\over 2}T_{\bar kpq}\,dz^q\wedge dz^p\wedge d\bar z^k$ of a Hermitian metric $\o$ by
\bea
T=i\p\o, \qquad \bar T=-i\bar\p\o,
\eea
and we also introduce the $(1,0)$-form $T_m$ and the $(0,1)$-form $\bar T_{\bar m}$ by
\bea
\label{TbarT}
T_m=g^{j\bar k}T_{\bar kjm},
\qquad
\bar T_{\bar m}=g^{\bar j k}\bar T_{k\bar j\bar m}.
\eea
Then 
\bea
(i\p\bar\p\o)_{\bar kj\bar\ell m}
= R_{\bar kj\bar\ell m}-R_{\bar km\bar\ell j}+R_{\bar\ell m\bar kj}-R_{\bar\ell j\bar km}
+g^{s\bar r}\,T_{\bar r mj}\bar T_{s\bar k\bar\ell}.
\eea
In general, there are several notions of Ricci curvature for Hermitian metrics, given by
\bea
\label{Riccis}
&&
R_{\bar kj}=R_{\bar kj}{}^p{}_p,\qquad \tilde R_{\bar kj}=R^p{}_p{}_{\bar kj},
\qquad
R_{\bar kj}'=R_{\bar k}{}^p{}_p{}_j,
\qquad
R_{\bar kj}''=R^p{}_{j\bar k}{}_p.
\eea
For metrics of the form $\o=e^u\hat\o$, where $\hat\o$ is K\"ahler and Ricci-flat, the following important relations between torsion and curvature hold,
\bea
\label{torsion1}
T_q(\o)=\p_q\log \|\Omega\|_\o, \qquad\bar T_{\bar q}(\o)=\p_{\bar q}\log\|\Omega\|_\o
\eea
and
\label{curvaturetorsion}
\bea
&&
R_{\bar kj}(\o)=2\nabla_{\bar k}T_j(\o)=2\na_j\bar T_{\bar k}(\o),
\nonumber\\
&&
R'_{\bar kj}(\o)=R''_{\bar kj}(\o)={1\over 2}R_{\bar kj}(\o).
\eea
Also, because $\hat\o$ is K\"ahler, we have
\bea
T(\o)=i\p u\wedge \o
\eea
so that the $(1,0)$-form $T_m$ actually determines in our case the full $(2,1)$ torsion tensor $T(\o)$. Henceforth, unless explicitly indicated otherwise, we shall designate by $T$ the $(1,0)$-form $T_mdz^m$ rather than the $(2,1)$-form $i\p\o$. For example, the norm $|T|^2$ will designate the expression
\bea
\label{T_omega}
|T|^2=g^{m\bar\ell}T_m\bar T_{\bar\ell}
\eea
rather than $|i\p\o|^2$ (which can be verified to be equal to $2|T|^2$).

\medskip

Using these relations, and the fact that we are in dimension $2$, we find
\bea
i\p\bar\p\o= {1\over 2}(-R+ 2|T|^2){\o^2\over 2}.
\eea
Substituting this equation and (\ref{p_tOmega-log}) in the flow (\ref{Flow1}), we obtain
\bea
\p_t\log \|\Omega\|_\o
=
{1\over \|\Omega\|_\o}
\bigg({R\over 2}- |T|^2 +2\alpha' {i\p\bar\p(\|\Omega\|_\o\rho)\over\o^2}
-{\alpha' \over 4} \sigma_2(i {\rm Ric}_\o )- \|\Omega\|_\o^2 \, \tilde\mu \bigg),
\eea
where we have introduced the time-independent, scalar function $\tilde\mu$ by $\mu=\tilde\mu{\hat\o^2\over 2}$, and the $\sigma_2$ operator with respect to the evolving metric
\be \label{sigma2_Ric_defn}
2 \sigma_2 (i {\rm Ric}_\o ) \, {\o^2 \over 2} = i {\rm Ric}_\o \wedge i {\rm Ric}_\o.
\ee
Since the metric $\o=e^u\hat\o$ is entirely determined by the conformal factor $e^u$, this flow for the volume form is equivalent to the flow of metrics (\ref{FY_parabolic}) quoted in the Introduction. The flow in terms of the conformal factor $u$ is easily worked out to be given by the equation (\ref{FY_parabolic_u}). 

\medskip

\subsection{Comparisons between the $3$-dimensional Anomaly flow and its $2$-dimensional reduction} \label{comparison-2d-3d}

It may be noteworthy that the flow (\ref{Flow1}) retains many of the features of the original Anomaly flow in $3$-dimensions. Indeed, as shown in \cite{PPZ1}, the conformally balanced condition (\ref{CB}) in dimension $3$ implies that the Hermitian metric $\chi$ on the $3$-fold $Y$ satisfies exactly the same relations (\ref{torsion1}) and (\ref{curvaturetorsion}) between torsion and curvature as the metric $\o=e^u\hat\o$ on the surface $X$,
\bea \label{T_j-p_jOmega}
T_q(\chi)=\p_q\log \|\Omega_Y\|_{\chi}, \qquad\bar T_{\bar q}(\chi)=\p_{\bar q}\log\|\Omega_Y\|_{\chi}
\eea
and
\bea \label{na-T_j}
&&
R_{\bar kj}(\chi)=2\nabla_{\bar k}T_j(\chi)=2\na_j\bar T_{\bar k}(\chi),
\nonumber\\
&&
R'_{\bar kj}(\chi)=R''_{\bar kj}(\chi)={1\over 2}R_{\bar kj}(\chi).
\eea
This suggests that the flow (\ref{FY_parabolic}) is interesting not just as a special case of the general Anomaly flow, but also as a good model for developing general methods for studying the flow.

\subsection{Starting the Flow}

In \cite{PPZ} general conditions were given for the short-time existence of the Anomaly flow, using the Nash-Moser implicit function theorem. However, the short-time existence of the flow can be seen more directly from the parabolicity of the flow, which holds when the form
\be \label{defn_o'}
\o' = e^u \hat{\o} + \alpha' e^{-u} \rho + \alpha' i \ddb u >0,
\ee
is positive definite. This can be seen from the scalar equation (\ref{FY_parabolic_u}). We will always assume that we start the flow from a large constant multiple of the background metric
\be
u(x,0) = \log M \gg 1, \ \ \ \o(0) = e^{u(0)} \hat{\o}= M \hat{\o}.
\ee
Recall that $\mu$ is defined in (\ref{mu}). In all that follows, we will assume that the cohomological condition
\be
\int_X \mu = 0,
\ee
is satisfied. 
Integrating (\ref{flow_forms}) and using the fact that $\| \Omega \|_{\o} \o^2 = e^u \hat{\o}^2$ gives the following conservation law
\be
{\p \over \p t} \int_X e^u \, {\hat{\o}^2 \over 2!} =0.
\ee
Hence
\be \label{conservation_law}
\int_X e^u \, {\hat{\o}^2 \over 2!} = M,
\ee
along the flow.

\

\section{The $C^0$ estimate of the conformal factor}
\setcounter{equation}{0}
In this section, we will work with equation (\ref{flow_forms}), since it will be easier to work with differential forms to obtain integral estimates. We let $\hat{\o}$ denote the fixed background K\"ahler form of $X$. We can rescale $\hat{\o}$ such that $\int_X {\hat{\o}^2 \over 2!} = 1$. We will omit the background volume form $\hat{\o}^2 \over 2!$ when integrating scalar functions. All norms in the current section will be taken with respect to the background metric $\hat{\omega}$. The starting point for the Moser iteration argument is to compute the quantity
\be
\int_X i \ddb (e^{-ku}) \wedge \o',
\ee
in two different ways. Recall that $\o'$ is defined in (\ref{defn_o'}). On one hand, by the definition of $\o'$ and Stokes' theorem, we have
\be
\int_X i \ddb (e^{-ku}) \wedge \o'  = \int_X \{ e^u \hat{\o} + \alpha' e^{-u} \rho \} \wedge i \ddb (e^{-ku}) .
\ee
Expanding
\bea \label{int_identity1}
\int_X i \ddb (e^{-ku}) \wedge \o'  &=& k^2 \int_X e^{-ku} \{ e^u \hat{\o}  + \alpha' e^{-u} \rho \} \wedge i \p u \wedge \bar{\p} u   \nonumber\\
&& -k \int_X e^{-ku} \{ e^u \hat{\o}  + \alpha' e^{-u} \rho \} \wedge i \ddb u  . 
\eea
On the other hand, without using Stokes' theorem, we obtain
\bea \label{int_identity2}
\int_X i \ddb (e^{-ku}) \wedge \o'  &=& k^2 \int_X e^{-ku} i \p u \wedge \bar{\p} u \wedge \o'   -k \int_X  e^{-ku} \{ e^u\hat{\o} + \alpha' e^{-u} \rho \} \wedge i \ddb u  \nonumber\\
&&- \alpha'k \int_X  e^{-ku} i \ddb u \wedge i \ddb u .
\eea
We equate (\ref{int_identity1}) and (\ref{int_identity2})
\bea
0 &=& - k^2 \int_X e^{-ku} i \p u \wedge \bar{\p} u \wedge \o'  +k^2 \int_X e^{-ku} \{ e^u \hat{\o}  + \alpha' e^{-u} \rho \} \wedge i \p u \wedge \bar{\p} u \nonumber\\
&&+  \alpha' k \int_X  e^{-ku} i \ddb u \wedge i \ddb u  .
\eea
Using equation (\ref{flow_forms}) and that $\| \Omega \|_{\o} \o^2 = e^u \hat{\o}^2$,
\bea
0 &=&   - k^2 \int_X e^{-ku} i \p u \wedge \bar{\p} u \wedge \o' +k^2 \int_X e^{-ku} \{ e^u \hat{\o}  + \alpha' e^{-u} \rho \} \wedge i \p u \wedge \bar{\p} u \nonumber\\
&& - 2k \int_X e^{-ku} \mu  -2k \int_X  e^{-ku} i \ddb (e^u \hat{\o} - \alpha' e^{-u} \rho)  + 4 k \int_X e^{-(k-1)u} \p_t u {\hat{\o}^2 \over 2!}.
\eea
Expanding out terms and dividing by $2k$ yields
\bea
0 &=&   - {k \over 2} \int_X e^{-ku} i \p u \wedge \bar{\p} u \wedge \o'+ {k \over 2} \int_X e^{-ku} \{ e^u \hat{\o}  + \alpha' e^{-u} \rho \} \wedge i \p u \wedge \bar{\p} u  -   \int_X e^{-ku} \mu  \nonumber\\
&&- \int_X e^{-(k-1)u} i \ddb u \wedge \hat{\o} 
 - \int_X e^{-(k-1)u} i \p u \wedge \bar{\p} u \wedge \hat{\o} - \alpha' \int_X  e^{-(k+1)u} i \ddb u \wedge \rho \nonumber\\
&&+  \alpha' \int_X  e^{-(k+1)u} i \p u \wedge \bar{\p} u \wedge \rho  +  \alpha' \int_X e^{-(k+1)u} i \ddb \rho  \nonumber\\
&&- 2 \alpha' \Re \int_X e^{-(k+1)u} i \p u \wedge \bar{\p} \rho  + 2 \int_X e^{-(k-1)u} \p_t u {\hat{\o}^2 \over 2!} .
\eea
Integration by parts gives
\bea \label{integrals_byparts}
0 &=&    - {k \over 2} \int_X e^{-ku} i \p u \wedge \bar{\p} u \wedge \o'  - {k \over 2} \int_X e^{-ku} \{ e^u \hat{\o}  + \alpha' e^{-u} \rho \} \wedge i \p u \wedge \bar{\p} u  \nonumber\\
&& - \int_X e^{-ku} \mu   +  \alpha' \int_X e^{-(k+1)u} i \ddb \rho-  \alpha' \Re \int_X e^{-(k+1)u} i \p u \wedge \bar{\p} \rho \nonumber\\
&&   +  2\int_X e^{-(k-1)u} \p_t u {\hat{\o}^2 \over 2!}.
\eea
One more integration by parts yields the following identity:
\bea \label{integrals_simplified}
&& {k \over 2} \int_X e^{-ku} \{ e^u \hat{\o} + \alpha' e^{-u} \rho \} \wedge i \p u \wedge \bar{\p} u  +  {\p \over \p t}{2 \over k-1} \int_X e^{-(k-1)u} {\hat{\o}^2 \over 2!} \\&=& - {k \over 2}  \int_X e^{-ku} i \p u \wedge \bar{\p} u \wedge \o'  -  \int_X e^{-ku} \mu   + (\alpha' - {\alpha' \over k+1}) \int_X  e^{-(k+1)u}  i \ddb \rho  . \nonumber 
\eea
The identity (\ref{integrals_simplified}) will be useful later to control the infimum of $u$, but to control the supremum of $u$, we replace $k$ with $-k$ in (\ref{integrals_simplified}). Then, for $k \neq 1$,
\bea \label{sup_key_identity}
&& {k \over 2} \int_X e^{(k+1)u} \{ \hat{\o} + \alpha' e^{-2u} \rho \} \wedge i \p u \wedge \bar{\p} u +  {\p \over \p t}{2 \over k+1} \int_X e^{(k+1)u} {\hat{\o}^2 \over 2!} \\&=& - {k \over 2}  \int_X e^{ku} i \p u \wedge \bar{\p} u \wedge \o' +  \int_X e^{ku} \mu   - (\alpha' - {\alpha' \over 1-k}) \int_X  e^{(k-1)u}  i \ddb \rho  . \nonumber 
\eea
\subsection{Estimating the supremum}
\begin{proposition} \label{sup_estimate}
Start the flow with initial data $e^{u(x,0)} = M$. Suppose the flow exists for $t \in [0,T)$ with $T>0$, and that $\inf_X e^u \geq 1$ and $\alpha' e^{-2u} \rho \geq  - {1 \over 2} \hat{\o}$ for all time $t \in [0,T)$. Then
\be 
\sup_{X \times [0,T)} e^u \leq C_1 M,
\ee
where $C_1$ only depends on $(X, \hat{\o})$, $\rho$, $\mu$, $\alpha'$.
\end{proposition}
{\it Proof:}  As long as the flow exists, we have
\be \label{pos_estimate}
i \p u \wedge \bar{\p} u \wedge \o'  \geq 0.
\ee
Let $\beta = {n \over n-1} = 2$. We can use (\ref{pos_estimate}), (\ref{sup_key_identity}), and $\alpha' e^{-2u} \rho \geq  - {1 \over 2} \hat{\o}$ to derive the following estimate for any $k \geq \beta$
\be \label{basic_e^u_int_estimate}
{k \over 4} \int_X e^{(k+1)u} |Du|^2 + {\p \over \p t}{2 \over k+1} \int_X e^{(k+1)u} \leq (\| \mu\|_{L^\infty} + 2 |\alpha'| \| \rho \|_{C^2}) \left( \int_X e^{ku} + \int_X e^{(k-1)u} \right).
\ee
Here we omit the background volume form ${\hat{\o}^2 \over 2!}$ when integrating scalars. 

\medskip

We now consider two cases: the case of small time and the case of large time. We must consider both these cases carefully because the objective is more than just to control $e^u$ uniformly in time; rather, we need to establish that $e^u$ stays comparable to the scale $M$ for all times. 
\smallskip
\par We begin with the estimate for large time. Suppose $T \in [n,n+1]$ for an integer $n \geq 1$. Let $n-1< \tau < \tau' < T$. Let $\zeta(t) \geq 0$ be a monotone function which is zero for $t \leq \tau$, identically $1$ for $t \geq \tau'$, and $|\zeta'| \leq 2 (\tau'-\tau)^{-1}$. Multiplying inequality (\ref{basic_e^u_int_estimate}) by $\zeta$ gives, for any $k \geq \beta$,
\bea
& \ & {k   \zeta \over 4} \int_X  e^{(k+1)u} |Du|^2 + {\partial \over \partial t} {2\zeta \over k+1} \int_X  e^{(k+1)u} \nonumber\\
&\leq&  (\| \mu\|_{L^\infty} + 2|\alpha'| \| \rho \|_{C^2}) \bigg\{ \zeta \int_X e^{(k-1)u} +  \zeta \int_X  e^{ku} \bigg\} + {2 \zeta' \over k+1} \int_X  e^{{(k+1)}u} .
\eea
Let $\tau' < s \leq T$. Integrating from $\tau$ to $s$ yields
\bea
&\ & {k  \over 4} \int_{\tau'}^{s} \int_X  e^{(k+1)u} |Du|^2 + {2 \over k+1} \int_X  e^{(k+1)u}(s)\\
&\leq& C \bigg\{ \int_{\tau}^T\int_X e^{(k-1)u} +  \int_{\tau}^T\int_X  e^{ku} + {1 \over \tau' - \tau} \int_{\tau}^T \int_X  e^{{(k+1)}u} \bigg\},
\eea
for any $k \geq \beta$, where $C$ only depends on $\alpha'$, $\rho$, $\mu$. We rearrange this inequality to obtain, for $k \geq \beta+1$,
\bea 
&\ & {(k-1) \over k} \int_{\tau'}^{s} \int_X |D e^{{k \over 2} u}|^2   + \int_X  e^{ku}(s)\nonumber\\
&\leq&  C k  \bigg\{ \int_{\tau}^T \int_X e^{(k-2)u} +  \int_{\tau}^T \int_X e^{(k-1)u} +  {1 \over \tau' - \tau} \int_{\tau}^T \int_X  e^{k u} \bigg\}.
\eea
Using $e^{-u} \leq 1$, 
\be \label{iteration_sup_est_cutoff}
 \int_{\tau'}^{s} \int_X  |D e^{{k \over 2} u}|^2 +  \int_X  e^{k u}(s) \leq C k \bigg\{ 1 +  {1 \over \tau' - \tau} \bigg\} \bigg\{ \int_{\tau}^T \int_X  e^{k u}\bigg\}.
\ee
The Sobolev inequality gives us 
\be
 \left( \int_X  e^{k \beta u} \right)^{1 \over \beta} \leq C'_X  \left( \int_X | e^{{k \over 2}u}|^2  + \int_X |D e^{{k \over 2}u}|^2 \right) ,
\ee
where $C_X'$ is the Sobolev constant on manifold $(X, \hat{\o})$.
Let $\beta^*$ be such that ${1 \over \beta} + {1 \over \beta^*} = 1$. By H\"older's inequality and the Sobolev inequality,
\bea \label{Holder-Sobolev}
\int_{\tau'}^{T} \int_X e^{k u} e^{{k \over \beta^*} u} &\leq& \int_{\tau'}^{T} \bigg(  \int_X e^{k \beta u} \bigg)^{1/\beta} \bigg( \int_X e^{ku} \bigg)^{1/\beta^*} \nonumber\\
&\leq& C_X' \sup_{t \in [\tau',T]} \bigg( \int_X e^{ku} \bigg)^{1/\beta^*} \int_{\tau'}^{T} \bigg\{ \int_X e^{ku} + \int_X |D e^{{k \over 2}u}|^2 \bigg\}.
\eea
Using estimate (\ref{iteration_sup_est_cutoff}), and defining $\gamma = 1+{1 \over \beta^*} = 1+{1 \over 2}$, we have for $k \geq 1+ \beta$,
\be
\bigg( \int_{\tau'}^{T} \int_X  e^{ \gamma k u} \bigg)^{1/ \gamma} \leq  Ck \bigg\{ 1+  {1 \over \tau' - \tau} \bigg\}  \int_{\tau}^T \int_X  e^{k u} .
\ee
We will iterate with $\tau_k = (n-1)+ \theta_1 - \gamma^{-k} (\theta_1 - \theta_2)$, for fixed $0< \theta_2 < \theta_1 \leq 1$.
\be
\bigg( \int_{\tau_{k+1}}^{T} \int_X  e^{ \gamma^{k+1} u} \bigg)^{1/ \gamma^{k+1}} \leq   \bigg\{ C \gamma^k +  (\theta_1-\theta_2)^{-1} {C \gamma^{2k} \over 1 - \gamma^{-1}} \bigg\}^{1 / \gamma^k} \bigg\{ \int_{\tau_k}^T \int_X  e^{\gamma^k u} \bigg\}^{1 / \gamma^k}.
\ee
Iterating, and using $\sum_i \gamma^{-i} = 3$, we see that for $p = \gamma^{\kappa_0} \geq 1 + \beta$, there holds 
\be
\sup_{X \times [n-1+\theta_1,T]} e^{u} \leq {C \over (\theta_1-\theta_2)^3} \, \| e^{u} \|_{L^{p}(X \times [n-1+\theta_2,T])},
\ee
where $C$ only depends on $(X,\hat{\o})$, $\rho$, $\mu$, and $\alpha'$. A standard argument can be used to relate the $L^{p}$ norm of $e^u$ to $\int_X e^u = M$. Indeed, by Young's inequality,
\bea
\sup_{X \times [n-1 + \theta_1,T]} e^{u} &\leq&  C (\theta_1-\theta_2)^{-3} \bigg( \sup_{X \times [n-1+\theta_2,T]} e^{(1 -1/p)u} \bigg) \bigg( \int_{X \times [n-1+\theta_2,T] } e^{u} \bigg)^{1/ p} \nonumber\\
&\leq& {1 \over 2}  \sup_{X \times [n-1+ \theta_2,T]} e^u + C (\theta_1-\theta_2)^{-3p} \int_{X \times [n-1,T] } e^{u},
\eea
for all $0<\theta_2<\theta_1 \leq 1$. We iterate this inequality with $\theta_0=1$ and $\theta_{i+1} = \theta_i - {1 \over 2}(1-\eta)\eta^{i+1} $, where $1/2 < \eta^{3p} < 1$. Then for each $k>1$,
\be
\sup_{X \times [n,T]} e^{u} \leq {1 \over 2^k} \left( \sup_{X \times[n-1+\theta_k,T]} e^u \right) + {2^{3p} \, CM \over (1-\eta)^{3p} \eta^{3p}} \sum_{i=0}^{k-1} \left( {1 \over 2 \eta^{3p}} \right)^{i}.
\ee
Taking the limit as $k \rightarrow \infty$, we obtain a constant $C$ depending only on $(X, \hat{\o})$, $\rho$, $\mu$, $\alpha'$ such that
\be \label{e^u-largeT-est}
\sup_{X \times [n,T]} e^{u} \leq C M,
\ee
for any $T \in [n,n+1]$ and integer $n \geq 1$.
\smallskip
\par Next, we adapt the previous estimate to the small time region $[0,T] \subseteq [0,1]$. The argument is similar in essence, and we provide all details for completeness. Integrating (\ref{basic_e^u_int_estimate}) from $0$ to $0 < s \leq T$ yields
\bea
 {k  \over 4} \int_{0}^{s} \int_X  e^{(k+1)u} |Du|^2 + {2 \over k+1} \int_X  e^{(k+1)u}(s) \leq C \bigg\{ \int_{0}^T\int_X e^{(k-1)u} +  \int_{0}^T\int_X  e^{ku} + M^{k+1} \bigg\}, \nonumber
\eea
for any $k \geq \beta$, where $C$ only depends on $\alpha'$, $\rho$, $\mu$. We rearrange this inequality to obtain, for $k \geq \beta+1$,
\bea 
{(k-1) \over k} \int_{0}^{s} \int_X |D e^{{k \over 2} u}|^2   + \int_X  e^{ku}(s) \leq  C k  \bigg\{ \int_{0}^T \int_X e^{(k-2)u} +  \int_0^T \int_X e^{(k-1)u} + M^k \bigg\}.
\eea
Using $e^{-u} \leq 1$, we obtain the following estimate, which holds uniformly for all $0 < s \leq T$.
\be \label{iteration_sup_est}
 \int_{0}^{s} \int_X  |D e^{{k \over 2} u}|^2 +  \int_X  e^{k u}(s) \leq C k  \bigg\{ \int_{0}^T \int_X  e^{k u} + M^k\bigg\}.
\ee
As estimate in (\ref{Holder-Sobolev}), by the H\"older and Sobolev inequalities there holds
\be
\int_{0}^{T} \int_X e^{k u} e^{{k \over \beta^*} u} \leq C_X' \sup_{s \in [0,T]} \bigg( \int_X e^{ku} \bigg)^{1/\beta^*} \int_{0}^{T} \bigg\{ \int_X e^{ku} + \int_X |D e^{{k \over 2}u}|^2 \bigg\}.
\ee
Recall that $\gamma = 1+{1 \over \beta^*}$. Thus for $k \geq 1+ \beta$,
\be \label{k-improve-by-gamma}
\int_{0}^{T} \int_X e^{k \gamma u} \leq  (Ck)^\gamma \bigg( \int_{0}^T \int_X  e^{k u} + M^k\bigg)^\gamma.
\ee
Therefore
\be
\bigg( \int_{0}^{T} \int_X e^{k \gamma u} + M^{k \gamma} \bigg)^{1/\gamma} \leq \bigg\{ (Ck)^\gamma \bigg( \int_{0}^T \int_X  e^{k u} + M^k\bigg)^\gamma + M^{k \gamma} \bigg\}^{1/\gamma},
\ee
and hence
\be
\bigg( \int_{0}^{T} \int_X  e^{k \gamma u} + M^{k \gamma} \bigg)^{1/ \gamma} \leq  Ck \bigg\{ \int_{0}^T \int_X  e^{k u} + M^k \bigg\}.
\ee
It follows that for all $\gamma^k \geq 1 + \beta$,
\be
\bigg( \int_{0}^{T} \int_X  e^{ \gamma^{k+1} u} + M^{\gamma^{k+1}} \bigg)^{1/ \gamma^{k+1}} \leq   \bigg\{ C \gamma^k \bigg\}^{1 / \gamma^k} \bigg\{ \int_{0}^T \int_X  e^{\gamma^k u} + M^{\gamma^k} \bigg\}^{1 / \gamma^k}.
\ee
Iterating, we see that for all $k$ such that $\gamma^k \geq \gamma^{\kappa_0} \geq 1 + \beta$,
\be
\bigg( \int_{0}^{T} \int_X  e^{ \gamma^{k+1} u} \bigg)^{1/ \gamma^{k+1}} \leq \bigg\{ \prod_{i=\kappa_0}^k \bigg( C \gamma^i \bigg)^{1 / \gamma^i} \bigg\}  \bigg\{ \int_{0}^T \int_X  e^{\gamma^{\kappa_0} u} + M^{\gamma^{\kappa_0}} \bigg\}^{1 / \gamma^{\kappa_0}}.
\ee
Sending $k \rightarrow \infty$, we obtain for $p = \gamma^{\kappa_0}$,
\be
\sup_{X \times [0,T]} e^{u} \leq C ( \| e^{u} \|_{L^{p}(X \times [0,T])} + M),
\ee
where $C$ only depends on $(X,\hat{\o})$, $\rho$, $\mu$, and $\alpha'$. Lastly, we relate the $L^{p}$ norm of $e^u$ to $\int_X e^u = M$. By the previous estimate
\bea
\sup_{X \times [0,T]} e^{u} \leq C \bigg( \sup_{X \times [0,T]} e^{(p -1)u} \bigg)^{1 \over p} \bigg( \int_{X \times [0,T] } e^{u} \bigg)^{1/ p} + CM.
\eea
We absorb the supremum term on the right-hand side using Young's inequality. Therefore,
\be \label{e^u-smallT-est}
\sup_{X \times [0,T]} e^{u} \leq CTM + CM \leq C M,
\ee
for any $0<T \leq 1$, and $C$ only depends on $(X,\hat{\o})$, $\rho$, $\mu$, and $\alpha'$.
\smallskip
\par By combining (\ref{e^u-largeT-est}) and (\ref{e^u-smallT-est}), we conclude the proof of Proposition \ref{sup_estimate}. Q.E.D.

\subsection{Estimating the infimum}
We introduce the constant
\be
\theta = {1 \over 2C_1 -1}.
\ee
Note that since $C_1 \geq 1$, we must have $0 < \theta \leq 1$. Fix a small constant $0<\delta<1$ such that
\be \label{delta_defn}
\delta < {\theta \over 4 C_X (|\alpha'|\|\rho\|_{C^2}+\|\mu\|_{C^0})}, \ \ {\textit and } \  \ \alpha' \delta^2 \rho \geq -{1\over 2} \hat{\o},
\ee
where $C_X$ is the Poincar\'e constant for the reference K\"ahler manifold $(X, \hat{\o})$. Define
\be
S_\delta := \{ t \in [0,T) : \sup_X e^{-u} \leq \delta \}.
\ee
Recall that we start the flow at $u_0=\log M$. It follows that if $M > \delta^{-1}$, then the flow starts in the region $S_\delta$. At any time $\hat{t} \in S_\delta$, we consider $U = \{ z \in X \ : \ e^{-u} \leq {2 \over M} \}$. Then by Proposition \ref{sup_estimate},
\be
M = \int_U e^u + \int_{X \backslash U} e^u \leq |U| \sup_X e^u + (1 - |U|) {M \over 2} \leq C_1 M |U| + (1- |U|) {M \over 2}.
\ee
It follows that at any $\hat{t}$,
\be
|U| > \theta >0 .
\ee
 We will also need the constant $C_0>1$ defined by
\be
C_0 = {1 \over 1 - {\theta \over 4}} (1+{2 \over \theta}) \bigg( {2 \over \theta^2} \bigg).
\ee
\subsubsection{Integral estimate}
\begin{proposition} \label{prop_int_est}
Start the flow at $u_0 = \log M$, where $M$ is large enough such that the flow starts in the region $S_\delta$. Suppose $[0,T] \subseteq S_\delta$. Then on $[0,T]$, there holds
\be
\int_X e^{-u} \leq {2 C_0 \over M}.
\ee
\end{proposition}
{\it Proof:} At $t=0$, we have $\int_X e^{-u} = {1 \over M} < {2C_0 \over M}$. Suppose $\hat{t} \in S_\delta$ is the first time when we reach $\int_X e^{-u} = {2C_0 \over M}$. Then we must have
\be
{\partial \over \partial t} \bigg|_{t=\hat{t}} \int_X e^{-u} \geq 0.
\ee
Setting $k=2$ in (\ref{integrals_simplified}) and dropping the negative term involving $\o' \geq 0$, we have
\bea
\int_X  e^{-u} \{\hat \o + \alpha' e^{-2u} \rho \}\wedge i \p u \wedge i \bar{\p} u + 2{\partial \over \partial t} \int_X  e^{-u} \leq \left(  |\alpha'| \|\rho\|_{C^2} \int_X e^{-3u} +  \|\mu\|_{C^0} \int_X  e^{-2u} \right). \nonumber
\eea
Since ${\partial \over \partial t}\bigg|_{t=\hat{t}}  \int_X  e^{-u} \geq 0$, and $e^{-u}\leq \delta<1$, there holds at $\hat{t}$,
\be
\int_X |D e^{-{u \over 2}}|^2 \leq (|\alpha'| \|\rho\|_{C^2}+\|\mu\|_{C^0}) \delta \int_X e^{-u}.
\ee
By the Poincar\'e inequality
\be
\int_X e^{-u} - \bigg( \int_X e^{- {u \over 2}} \bigg)^2 = \int_X \bigg| e^{- {u \over 2}} - \int_X e^{- {u \over 2}} \bigg|^2 \leq C_X \int_X |D e^{- {u \over 2}} |^2.
\ee
By (\ref{delta_defn}), we have
\be
\int_X e^{-u} - \bigg( \int_X e^{- {u \over 2}} \bigg)^2 \leq {\theta \over 4} \int_X e^{-u},
\ee
and it implies
\be \label{reverse_poincare}
\int_X e^{-u} \leq {1 \over 1 - {\theta \over 4}} \bigg( \int_X e^{- {u \over 2}} \bigg)^2.
\ee
We may use the measure estimate and (\ref{reverse_poincare}) to obtain
\bea
\bigg( \int_X e^{- {u \over 2}} \bigg)^2 &\leq& (1 + {2 \over \theta}) \bigg( \int_U e^{- {u \over 2}} \bigg)^2 + (1 + {\theta \over 2}) \bigg( \int_{X \backslash U} e^{- {u \over 2}} \bigg)^2 \nonumber\\
&\leq& (1+{2 \over \theta}) |U| \int_U e^{-u} + (1+{\theta \over 2})(1-|U|) \int_{X \backslash U} e^{-u} \nonumber\\
&\leq& (1+{2 \over \theta}) {2 \over M} + (1+{\theta \over 2}) (1-\theta) {1 \over 1- {\theta \over 4}} \bigg( \int_X e^{- {u \over 2}} \bigg)^2.
\eea
Thus
\be
\bigg( \int_X e^{- {u \over 2}} \bigg)^2 \leq (1+{2 \over \theta}) {2 \over M} \bigg( {1 \over 1 - (1+ {\theta \over 2}) (1-\theta) (1-{\theta \over 4})^{-1}} \bigg).
\ee
For any $\theta \geq 0$, we have the elementary estimate
\be
(1+ {\theta \over 2}) (1-\theta) (1-{\theta \over 4})^{-1} \leq 1 - \theta^2.
\ee
Using this and (\ref{reverse_poincare}),
\be
\int_X e^{-u} \leq {1 \over 1 - {\theta \over 4}} (1+{2 \over \theta}) \bigg( {2 \over \theta^2} \bigg) {1 \over M}  = {C_0 \over M}.
\ee
This contradicts that $\int_X e^{-u} = {2C_0 \over M}$ at $\hat{t}$. It follows that $\int_X e^{-u}$ stays less than ${2 C_0 \over M}$ for all time $t \in S_\delta$.

\subsubsection{Iteration}

\begin{proposition} \label{inf_estimate}
Start the flow with initial data $e^{u(x,0)} = M$. Suppose the flow exists for $t \in [0,T)$ with $T>0$, and $[0,T)\subseteq S_\delta$. Then
\be
\sup_{X \times [0,T)} e^{-u} \leq {C_2 \over M},
\ee
where $C_2$ only depends on $(X,\hat{\o})$, $\rho$, $\mu$, $\alpha'$.
\end{proposition}
{\it Proof:} We can drop the negative terms involving $\o' \geq 0$ and use $\alpha' e^{-2u} \rho \geq - {1 \over 2} \hat{\o}$ in (\ref{integrals_simplified}) to obtain the estimate, for $k\geq 2$,
\be\label{basic_e^-u_int_estimate}
{k \over 4} \int_X  e^{-(k-1)u} |Du|^2 + {\partial \over \partial t} {2 \over k-1} \int_X  e^{-(k-1)u} \leq C \bigg( \int_X e^{-(k+1)u}  + \int_X  e^{-ku} \bigg).
\ee
As in the upper bound on $e^u$, we split the argument into the cases of large time and small time, and first consider the case of large time. 
\smallskip
\par Suppose $T \in [n,n+1]$ for an integer $n \geq 1$. Let $n-1 < \tau < \tau' < T$. Let $\zeta(t) \geq 0$ be a monotone function which is zero for $t \leq \tau$ and identically $1$ for $t \geq \tau'$. Multiplying (\ref{basic_e^-u_int_estimate}) by $\zeta$ gives
\be
{k \zeta \over 4} \int_X  e^{-(k-1)u} |Du|^2 + {\partial \over \partial t} {2\zeta \over k-1} \int_X  e^{-(k-1)u} \leq C \bigg\{ \zeta \int_X e^{-(k+1)u} +  \zeta \int_X  e^{-ku} + \zeta' \int_X  e^{-{(k-1)}u} \bigg\}.
\ee
Let $\tau' < s \leq T$. Integrating from $\tau$ to $s$
\bea
&\ & {k \over 4} \int_{\tau'}^{s} \int_X  e^{-(k-1)u} |Du|^2 + {2 \over k-1} \int_X  e^{-(k-1)u}(s)\nonumber\\
&\leq& C \bigg\{ \int_{\tau}^T \int_X e^{-(k+1)u} +  \int_{\tau}^T \int_X e^{-ku} +  {1 \over \tau' - \tau} \int_{\tau}^T \int_X  e^{-{(k-1)}u} \bigg\}.
\eea
We rearrange this inequality to obtain, for $k \geq 1$,
\bea
\int_{\tau'}^{s} \int_X |D e^{-{k \over 2} u}|^2   +  2\int_X  e^{-ku}(s) \leq C k \bigg\{ \int_{\tau}^T \int_X e^{-(k+2)u} +  \int_{\tau}^T \int_X e^{-(k+1)u} +  {1 \over \tau' - \tau} \int_{\tau}^T \int_X  e^{-ku} \bigg\}. \nonumber
\eea
Since $e^{-u}\leq \delta <1$, we have
\be \label{inf_iteration_est1}
\int_{\tau'}^{s} \int_X |D e^{-{k \over 2} u}|^2   +  2\int_X  e^{-ku}(s) \leq C  k \bigg\{ 1+ {1 \over \tau' - \tau} \bigg\} \bigg\{ \int_{\tau}^T \int_X e^{-ku}  \bigg\}. 
\ee
Recall that we denote $\beta= {n \over n-1} = 2$, $\beta^*$ such that ${1 \over \beta} + {1 \over \beta^*} =1$, and $\gamma = 1 + {1 \over \beta^*}$. By the Sobolev inequality
\bea
\int_{\tau'}^{T} \int_X e^{-ku} e^{-{k \over \beta^*} u} &\leq& \int_{\tau'}^{T} \bigg(  \int_X e^{-k \beta u} \bigg)^{1/\beta} \bigg( \int_X e^{-ku} \bigg)^{1/\beta^*} \nonumber\\
&\leq& C \sup_{t \in [\tau',T]} \bigg( \int_X e^{-ku} \bigg)^{1/\beta^*} \int_{\tau'}^{T} \bigg\{ \int_X e^{-ku} + \int_X |D e^{-{k \over 2}u}|^2 \bigg\}.
\eea
Using estimate (\ref{inf_iteration_est1}), we arrive at
\be 
\bigg( \int_{\tau'}^{T} \int_X  e^{-\gamma k  u} \bigg)^{1/\gamma} \leq C  k \bigg\{1+  {1 \over \tau' - \tau} \bigg\}  \bigg\{ \int_{\tau}^T \int_X e^{-ku}   \bigg\}.
\ee
Iterating with $\tau_k = (1 - \gamma^{-(k+1)}) + (n-1)$, 
\be
\bigg( \int_{\tau_{k+1}}^T \int_X e^{- \gamma^{k+1}u} \bigg)^{1/\gamma^{k+1}} \leq \bigg\{ C \gamma^k + {C \gamma^{2k} \over 1 - \gamma^{-1}} \bigg\}^{1/\gamma^k} \bigg\{ \int_{\tau_k}^T \int_X e^{\gamma^k u} \bigg\}^{1/\gamma^k}.
\ee
Note $\tau_k \geq n- {2 \over 3}$. Sending $k \rightarrow \infty$, we have the $C^0$ estimate
\be
\sup_{X \times [n,T]} e^{-u} \leq C \| e^{-u} \|_{L^{1}(X \times [n-{2 \over 3},T])}.
\ee
By Proposition \ref{prop_int_est}, for $n \leq T \leq n+1$ and $n \geq 1$, we obtain
\be \label{e^-u-largeT-est}
\sup_{X \times [n,T]} e^{-u} \leq {C \over M}.
\ee
Next, we consider the small time region $[0,T] \subseteq [0,1]$. Integrating (\ref{basic_e^-u_int_estimate}) from $0$ to $0<s<T$, we obtain
\bea
{k \over 4} \int_{0}^{s} \int_X  e^{-(k-1)u} |Du|^2 + {2 \over k-1} \int_X  e^{-(k-1)u}(s) \leq C \bigg\{ \int_{0}^T \int_X e^{-(k+1)u} +  \int_{0}^T \int_X e^{-ku} + {M^{-(k-1)} \over k-1} \bigg\}. \nonumber
\eea
We rearrange this inequality to obtain, for $k \geq 1$,
\bea
\int_{0}^{s} \int_X |D e^{-{k \over 2} u}|^2   +  2\int_X  e^{-ku}(s)
\leq C k \bigg\{ \int_{0}^T \int_X e^{-(k+2)u} +  \int_{0}^T \int_X e^{-(k+1)u} + M^{-k} \bigg\}.
\eea
Since $e^{-u}\leq \delta <1$, we have
\be \label{inf_iteration_est2}
\int_{0}^{s} \int_X |D e^{-{k \over 2} u}|^2   +  2\int_X  e^{-ku}(s) \leq C  k \bigg\{ \int_{0}^T \int_X e^{-ku} + M^{-k} \bigg\}. 
\ee
As before, by the Sobolev inequality
\be
\int_{0}^{T} \int_X e^{-ku} e^{-{k \over \beta^*} u} \leq C \sup_{s \in [0,T]} \bigg( \int_X e^{-ku} \bigg)^{1/\beta^*} \int_{0}^{T} \bigg\{ \int_X e^{-ku} + \int_X |D e^{-{k \over 2}u}|^2 \bigg\}.
\ee
Combining this with (\ref{inf_iteration_est2}) yields
\be 
\bigg( \int_{0}^{T} \int_X  e^{-\gamma k  u} + M^{-\gamma k} \bigg)^{1/\gamma} \leq C  k  \bigg\{ \int_{0}^T \int_X e^{-ku} + M^{-k}  \bigg\}.
\ee
Iterating, we obtain the $C^0$ estimate
\be
\sup_{X \times [0,T]} e^{-u} \leq C \| e^{-u} \|_{L^{1}(X \times [0,T])} + CM^{-1}.
\ee
By Proposition \ref{prop_int_est}, for $0<T \leq 1$ we obtain
\be \label{e^-u-smallT-est}
\sup_{X \times [0,T)} e^{-u} \leq C T M^{-1} + C M^{-1} \leq {C \over M}.
\ee
By combining (\ref{e^-u-largeT-est}) and (\ref{e^-u-smallT-est}), we conclude the proof of Proposition \ref{inf_estimate}. Q.E.D.

\begin{theorem} \label{C0-thm}
Suppose the flow exists for $t \in [0,T)$, and initially starts with $u_0 = \log M$. There exists $M_0 \gg 1$ such that for all $M \geq M_0$, there holds
\be
 \sup_{X\times[0,T)} e^{u} \leq C_1 M, \ \ \sup_{X\times[0,T)} e^{-u} \leq {C_2 \over M},
\ee
where $C_1, C_2$ only depend on $(X, \hat{\o})$, $\rho$, $\mu$, $\alpha'$.
\end{theorem}
{\it Proof:} By Proposition \ref{sup_estimate} and Proposition \ref{inf_estimate}, the estimates hold as long as we stay in $S_\delta$. Choose $M_0$ such that
\be
{C_2 \over M_0} < {\delta\over 2},
\ee
where recall $\delta$ is defined in (\ref{delta_defn}). Then at $t=0$, we have $e^{-u_0} < \delta$, and the estimate is preserved on $[0,T)$. The theorem follows. Q.E.D.

\

\section{Evolution of the torsion}
\setcounter{equation}{0}
Before proceeding, we clearly state the conventions and notation that will be used for the maximum principle estimates of Sections \S 4-6. All norms from this point on will be with respect to the evolving metric $\o = e^u \hat{\o}$, unless denoted otherwise. We will write $\o = i g_{\bar{k} j} dz^j \wedge d \bar{z}^k$. We will use the Chern connection of $\o$ to differentiate
\be
\na_{\bar{k}} V^\alpha = \p_{\bar{k}} V^\alpha, \ \ \na_k V^\alpha = g^{\alpha \bar{\beta}} \p_k (g_{\bar{\beta} \gamma} V^\gamma).
\ee
The curvature of the metric $\o$ is 
\be \label{curv_defn}
R_{\bar{k} j}{}^\alpha{}_\beta = - \p_{\bar{k}} (g^{\alpha \bar{\gamma}} \p_j g_{\bar{\gamma} \beta}) = \hat{R}_{\bar{k} j}{}^\alpha{}_\beta - u_{\bar{k} j} \delta^\alpha{}_\beta.
\ee
The torsion tensor of the metric $\o$ is $T_{\bar{k} mj} = \p_m g_{\bar{k} j} - \p_j g_{\bar{k} m}$, and since $\hat{\o}$ has zero torsion, we may compute
\be \label{T^lambda}
T^\lambda{}_{mj} = g^{\lambda \bar{k}} T_{\bar{k} mj} = u_m \delta^\lambda{}_j - u_j \delta^\lambda{}_m.
\ee
We note the following formulas for the torsion and Chern-Ricci curvature of the evolving metric
\be \label{T_j-defn}
R_{\bar{k} j} = R_{\bar{k} j}{}^\alpha{}_\alpha= -2 u_{\bar{k} j}, \ \ T_j = T^\lambda{}_{\lambda j} = - \p_j u.
\ee
Recall that $|T|^2$ refers to the norm of $T_j$, as noted in (\ref{T_omega}). We will often use the following commutation formulas to exchange covariant derivatives
\be
[\na_j, \na_{\bar{k}}] \, V_i = - R_{\bar{k} j}{}^p{}_i V_p, \ \ [\na_j, \na_k] \, V_i = - T^\lambda{}_{jk} \na_\lambda V_i.
\ee
To handle the differentiation of the equation, we will rewrite the terms involving $\rho$ in the flow (\ref{FY_parabolic}). Compute
\bea
- \alpha' i \ddb (e^{-u} \rho) &=& -\alpha' e^{-u} i \ddb \rho + 2 \alpha' \Re \{ e^{-u} i \p u \wedge \bar{\p} \rho \} \nonumber\\
&& + \alpha' e^{-u} i \ddb u \wedge \rho - \alpha' i e^{-u} \p u \wedge \bar{\p} u \wedge \rho.
\eea
We introduce the notation
\be
- \alpha' i \ddb (e^{-u} \rho) = \bigg( -\alpha' e^{-u} \psi_\rho + \alpha' e^{-u} \Re \{ b^i_\rho u_i \} + \alpha'  e^{-u} \tilde{\rho}^{j \bar{k}} u_{\bar{k} j} - \alpha' e^{-u} \tilde{\rho}^{p \bar{q}} u_p \bar{u}_{\bar{q}}  \bigg) {\hat{\o}^2 \over 2},
\ee
where $\psi_\rho(z)$, $b_\rho^i(z)$, $\tilde{\rho}^{j \bar{k}}(z)$ are defined one by one corresponding to the previous expression. We note that $\psi_\rho$, $b^i_\rho$, $\tilde{\rho}^{j \bar{k}}$ are bounded in $C^\infty$ by constants depending only on the form $\rho$ and the background metric $\hat{\o}$. We also note that $\tilde{\rho}^{j \bar{k}}$ is Hermitian since $\rho$ is real. We may rewrite this expression as
\be
- \alpha' i \ddb (e^{-u} \rho) = \bigg( -\alpha' e^{-3u} \psi_\rho - \alpha' e^{-3u} \Re \{ b^i_\rho T_i \} - {\alpha' \over 2} e^{-3u} \tilde{\rho}^{j \bar{k}} R_{\bar{k} j} - \alpha' e^{-3u} \tilde{\rho}^{p \bar{q}} T_p \bar{T}_{\bar{q}}  \bigg) {\o^2 \over 2}.
\ee
 With all the introduced notation, we can write the flow (\ref{FY_parabolic}) in the following way.
\be \label{Flow3}
\p_t g_{\bar{k} j} = {1 \over 2 \| \Omega \|_{\o}} \bigg(- {R \over 2} -  {\alpha' \over 2} \| \Omega \|_{\o}^3 \tilde{\rho}^{p \bar{q}} R_{\bar{q} p}    + {\alpha' \over 4} \sigma_2(i {\rm Ric}_\o) + |T|^2+ \| \Omega \|_{\o}^2 \, \nu  \bigg) \, g_{\bar{k} j},
\ee
where 
\be \label{defn_nu}
\nu = -\alpha' \| \Omega \|_{\o} \psi_\rho - \alpha' \| \Omega \|_{\o} \Re \{ b^i_\rho T_i \}  - \alpha'\| \Omega \|_{\o}  \tilde{\rho}^{p \bar{q}} T_p \bar{T}_{\bar{q}}  +  \tilde{\mu}.
\ee
In the following, we will use $\|\Omega\|$ to replace $\|\Omega\|_{\o}$ for simplicity, if there is no confusing of the notation.

\subsection{Torsion tensor}
Using $\| \Omega \|=e^{-u}$ and $g_{\bar{k} j}=e^u \hat{g}_{\bar{k} j}$, (\ref{Flow3}) implies the following evolution of $\| \Omega \|$,
\be \label{p_t-log|Omega|-from-eq}
\p_t \log \| \Omega \| =  {1 \over 2 \| \Omega \|} \bigg({R \over 2} + {\alpha' \over 2} \| \Omega \|^3 \tilde{\rho}^{p \bar{q}} R_{\bar{q} p} -  |T|^2  - {\alpha' \over 4} \sigma_2(i {\rm Ric}_\o) -  \| \Omega \|^2 \, \nu  \bigg).
\ee
Using (\ref{torsion1}) and (\ref{p_t-log|Omega|-from-eq}), we evolve
\bea
\p_t T_j &=& \p_j \p_t \log \| \Omega \| \nonumber\\
&=& \na_j \bigg\{ {1 \over 2 \| \Omega \|} \bigg({R \over 2} + {\alpha' \over 2} \| \Omega \|^3 \tilde{\rho}^{p \bar{q}} R_{\bar{q} p} -  |T|^2 - {\alpha' \over 4} \sigma_2(i {\rm Ric}_\o) -  \| \Omega \|^2 \, \nu  \bigg) \bigg\}.
\eea
Using $\p_j \| \Omega \| = \| \Omega \| T_j$ and the definition of $\nu$ (\ref{defn_nu}), a straightforward computation gives
\bea \label{p_t-T1}
\p_t T_j &=&  {1 \over 2 \| \Omega \|} \bigg\{-{1 \over 2} T_j R + T_j |T|^2  +  {\alpha' \over 4} T_j \sigma_2(i {\rm Ric}_\o)  \nonumber\\
&& + {1 \over 2} \na_j R + {\alpha' \over 2}  \| \Omega \|^3 \tilde{\rho}^{p \bar{q}} \na_j R_{\bar{q} p} -  \na_j |T|^2 - {\alpha' \over 4} \na_j \sigma_2(i {\rm Ric}_\o)  + E_j \bigg\},
\eea
where 
\bea \label{defn_E}
E_j &=& 2 \alpha' \| \Omega \|^3 \psi_\rho T_j + 2 \alpha' \| \Omega \|^3  \Re \{ b^i_{\rho} T_i \} T_j + \alpha' \| \Omega \|^3  \tilde{\rho}^{p \bar{q}} R_{\bar{q} p} T_j \nonumber\\
&&+ 2 \alpha' \| \Omega \|^3 (\tilde{\rho}^{p \bar{q}} T_p \bar{T}_{\bar{q}} )  T_j -  \| \Omega \|^2 \tilde{\mu} T_j + \alpha' \| \Omega \|^3 \na_j \psi_\rho \nonumber\\
&& +  \alpha' \| \Omega \|^3  \Re \{ \na_j b^i_{\rho} T_i \} + \alpha' \| \Omega \|^3  \Re \{ b^i_{\rho} \na_j T_i \} + {\alpha' \over 2}  \| \Omega \|^3 (\na_j \tilde{\rho}^{p \bar{q}}) R_{\bar{q} p} \nonumber\\
&& + \alpha'  \| \Omega \|^3 (\na_j \tilde{\rho}^{p \bar{q}}) T_p \bar{T}_{\bar{q}} + \alpha'  \| \Omega \|^3 \tilde{\rho}^{p \bar{q}} \na_j T_p \bar{T}_{\bar{q}}+ {\alpha' \over 2}  \| \Omega \|^3 \tilde{\rho}^{p \bar{q}} T_p R_{\bar{q} j} \nonumber\\
&& - \| \Omega \|^2 \na_j \tilde{\mu}.
\eea
Our reason for treating $E_j$ as an error term is that the $C^0$ estimate tells us that $\| \Omega \| = e^{-u} \ll 1$ if we start the flow from a large enough constant $\log M$. As we will see, the terms appearing in $E_j$ will only slightly perturb the coefficients of the leading terms in the proof of Theorem \ref{torsion_estimate}.
\medskip
\par
We need to express the highest order terms in (\ref{p_t-T1}) as the linearized operator acting on torsion. First, we write the Ricci curvature in terms of the conformal factor 
\be
\na_j R_{\bar{q} p} = -2 \na_j \na_p \na_{\bar{q}} u .
\ee
Exchanging covariant derivatives
\be
-2 \na_j \na_p \na_{\bar{q}} u = -2 \na_p \na_{\bar{q}} \na_j u - 2 T^\lambda{}_{pj} \na_\lambda \na_{\bar{q}} u.
\ee
It follows from (\ref{T_j-defn}) that
\be
\na_j R_{\bar{q} p} = 2\na_p \na_{\bar{q}} T_j + T^\lambda{}_{pj} R_{\bar{q} \lambda} .
\ee
Hence
\bea \label{na_R-sigma2Ric}
& \ &  \na_j R - {\alpha' \over 2} \na_j \sigma_2(i {\rm Ric}_\o) + \alpha'  \| \Omega \|^3 \tilde{\rho}^{p \bar{q}} \na_j R_{\bar{q} p} \nonumber\\
&=& g^{p \bar{q}} \na_j R_{\bar{q} p} + \alpha'  \| \Omega \|^3 \tilde{\rho}^{p \bar{q}} \na_j R_{\bar{q} p} - {\alpha' \over 2} \sigma_2^{p \bar{q}} \na_j R_{\bar{q} p} \nonumber\\
&=& 2 F^{p \bar{q}} \na_p \na_{\bar{q}} T_j +  F^{p \bar{q}} T^\lambda{}_{pj} R_{\bar{q} \lambda},
\eea
where we introduced the notation
\be \label{sigma2^pq-defn}
\sigma_2^{p \bar{q}} = R  \,g^{p \bar{q}} - R^{p \bar{q}},
\ee
and
\be \label{F^pq-defn}
F^{p \bar{q}} = g^{p \bar{q}} + \alpha' \| \Omega \|^3 \tilde{\rho}^{p \bar{q}} - {\alpha' \over 2} (R  \,g^{p \bar{q}} - R^{p \bar{q}}).
\ee
The tensor $F^{p \bar{q}}$ is Hermitian, and in Section \S \ref{section_curvature} we will show that $F^{p \bar{q}}$ stays close to $g^{p \bar{q}}$ along the flow. Substituting (\ref{na_R-sigma2Ric}) into (\ref{p_t-T1})
\bea \label{p_t-T}
\p_t T_j &=&  {1 \over 2 \| \Omega \|} \bigg\{ F^{p \bar{q}} \na_p \na_{\bar{q}} T_j -  \na_j |T|^2 - {1 \over 2} T_j R + {\alpha'  \over 4} T_j \sigma_2(i {\rm Ric}_\o) \nonumber\\
&&+ {1 \over 2} F^{p \bar{q}} T^\lambda{}_{pj} R_{\bar{q} \lambda} +   T_j |T|^2  +E_j \bigg\}.
\eea
\par Before proceeding, let us discuss $\sigma_2^{p \bar{q}}$ and $F^{p \bar{q}}$ using convenient coordinates. Suppose we work at a point where the evolving metric $g_{i \bar j} = \delta_{ij}$ and $R_{\bar{k} j}$ is diagonal. Let $A^i{}_j = g^{i \bar{k}} R_{\bar{k} j}$. The function $\sigma_2(A^i{}_j)$ maps a Hermitian endomorphism to the second elementary symmetric polynomial of its eigenvalues. We are working in dimension $n=2$, so $\sigma_2(A^i{}_j)$ is the product of the two eigenvalues of $A$. Our operator $\sigma_2(i {\rm Ric}_\o)$ defined in (\ref{sigma2_Ric_defn}) is with respect to the evolving metric $\o$, so denoting $A^i{}_j = g^{i \bar{k}} R_{\bar{k} j}$, we have $\sigma_2(i {\rm Ric}_\o)=\sigma_2(A)$. We define $\sigma_2^{p \bar{q}} = {\p \sigma_2 \over \p A^k{}_p} g^{k \bar{q}}$. It is well-known that ${\p \sigma_2 \over \p A^1{}_1} = A^2{}_2$, ${\p \sigma_2 \over \p A^2{}_2} = A^1{}_1$, and ${\p \sigma_2 \over \p A^1{}_2}=0$ if $A$ is diagonal. Then in our case,
\be
\sigma_2^{1 \bar{1}} = R_{\bar{2} 2}, \ \ \sigma_2^{2 \bar{2}} =  R_{\bar{1} 1}, \ \ \sigma_2^{1 \bar{2}}=\sigma_2^{2 \bar{1}}=0.
\ee
We obtain
\bea \label{F^ii-in-coords}
&\ & F^{1 \bar{1}} = 1 + \alpha' \| \Omega \|^3 \tilde{\rho}^{1 \bar{1}} - {\alpha' \over 2} R_{\bar{2} 2}, \ \ F^{2 \bar{2}} = 1 + \alpha' \| \Omega \|^3 \tilde{\rho}^{2 \bar{2}} - {\alpha' \over 2} R_{\bar{1} 1}, \nonumber\\
& \ & \ \ F^{1 \bar{2}} = \alpha' \| \Omega \|^3 \tilde{\rho}^{1 \bar{2}}, \ \   F^{2 \bar{1}} =\alpha' \| \Omega \|^3 \tilde{\rho}^{2 \bar{1}}.
\eea

\subsection{Norm of the torsion}
We will compute
\be
\p_t |T|^2 = \p_t \{ g^{i \bar{j}} T_i \bar{T}_{\bar{j}} \}.
\ee
We have
\be
\p_t g^{i \bar{j}} = - g^{i \bar{\lambda}} g^{\gamma \bar{j}} \p_t g_{\bar{\lambda} \gamma} =  {1 \over 2 \| \Omega \| } \bigg( {R \over 2} +{\alpha' \over 2} \| \Omega \|^3 \tilde{\rho}^{p \bar{q}} R_{\bar{q} p}  - {\alpha' \over 4} \sigma_2(i {\rm Ric}_\o) -  |T|^2 -  \| \Omega \|^2 \, \nu  \bigg) \, g^{i\bar{j} }.
\ee
Hence
\bea \label{p_t-normT1}
\p_t |T|^2 &=&  2 \Re \langle \p_t T, T \rangle \nonumber\\
&&+  {| T |^2 \over 2 \| \Omega \|} \bigg( {R \over 2} +{\alpha' \over 2} \| \Omega \|^3 \tilde{\rho}^{p \bar{q}} R_{\bar{q} p}  - {\alpha' \over 4} \sigma_2(i {\rm Ric}_\o)  -  |T|^2 - \| \Omega \|^2 \, \nu  \bigg)
\eea
Next, using the notation $|W|^2_{Fg}=F^{p \bar{q}} g^{i \bar{j}} W_{pi} \bar{W}_{\bar{q} \bar{j}}$,
\bea
F^{p \bar{q}} \na_p \na_{\bar{q}} |T|^2 &=& F^{p \bar{q}} g^{i \bar{j}} \na_p \na_{\bar{q}} T_i \bar{T}_{\bar{j}} + F^{p \bar{q}} g^{i \bar{j}}  T_i \na_{p} \na_{\bar{q}} \bar{T}_{\bar{j}} + |\na T|^2_{Fg} + |\overline{\na} T|^2_{Fg} \nonumber\\
&=& F^{p \bar{q}} g^{i \bar{j}} \na_p \na_{\bar{q}} T_i \bar{T}_{\bar{j}} +  g^{i \bar{j}}  T_i \overline{F^{q \bar{p}} \na_{q} \na_{\bar{p}}T_{j}} + F^{p \bar{q}} g^{i \bar{j}}  T_i R_{\bar{q} p \bar{j}}{}^{\bar{\lambda}} \bar{T}_{\bar{\lambda}}\nonumber\\
&& + |\na T|^2_{Fg} + |\overline{\na} T|^2_{Fg}.
\eea
We introduce the notation $\Delta_F = F^{p \bar{q}} \na_p \na_{\bar{q}}$. We have shown
\be \label{Delta-normT}
\Delta_F | T |^2 =  2 \Re \langle \Delta_F T,T \rangle  + |\na T |^2_{Fg} + |\overline{\na} T |^2_{Fg} +F^{p \bar{q}} g^{i \bar{j}}  T_i R_{\bar{q} p \bar{j}}{}^{\bar{\lambda}} \bar{T}_{\bar{\lambda}} .
\ee
Combining (\ref{p_t-T}), (\ref{p_t-normT1}), and (\ref{Delta-normT}), we obtain
\bea \label{p_t-|T|^2}
\p_t |T|^2 &=&   {1 \over 2 \| \Omega \|} \bigg\{ \Delta_F | T |^2 - |\na T |^2_{Fg} - |\overline{\na} T |^2_{Fg} -  2 \Re \{ g^{i \bar{j}} \na_i |T|^2 \bar{T}_{\bar{j}}\}   \nonumber\\
&& -{1 \over 2} R |T|^2 + {\alpha' \over 4} \sigma_2(i {\rm Ric}_\o) |T|^2 + \Re \{ F^{p \bar{q}} g^{i \bar{j}} T^\lambda{}_{pi} R_{\bar{q} \lambda} \bar{T}_{\bar{j}} \} \nonumber\\
&& -F^{p \bar{q}} g^{i \bar{j}}  T_i R_{\bar{q} p \bar{j}}{}^{\bar{\lambda}} \bar{T}_{\bar{\lambda}}  + |T|^4 + {\alpha' \over 2} \| \Omega \|^3 \tilde{\rho}^{p \bar{q}} R_{\bar{q} p} |T|^2 \nonumber\\
&&- \| \Omega \|^2  |T|^2\nu  + 2 \Re \langle E,T \rangle \bigg\}.
\eea

\subsection{Estimating the torsion}
\begin{theorem} \label{torsion_estimate}
There exists $M_0 \gg 1$ such that all $M \geq M_0$ have the following property. Start the flow with a constant function $u_0 = \log M$. If 
\be
|\alpha' {\rm Ric}_\o| \leq 10^{-6} 
\ee
along the flow, then there exists $C_3>0$ depending only on $(X,\hat{\o})$, $\rho$, $\tilde{\mu}$ and $\alpha'$, such that
\be
|T|^2 \leq {C_3 \over M} \ll 1.
\ee
\end{theorem}
Denote $\Lambda=1+{1 \over 8}$. We will study the test function 
\be
G = \log |T|^2 - \Lambda \log \| \Omega \|.
\ee
Taking the time derivative gives us
\be \label{log_torsion_1}
\p_t G = {\p_t |T|^2 \over |T|^2} - \Lambda \p_t \log \| \Omega \|.
\ee
Computing using (\ref{torsion1}) and (\ref{F^pq-defn}),
\bea
\Delta_F \log \| \Omega \| &=& F^{p \bar{q}} \na_p \bar{T}_{\bar{q}} = {1 \over 2} F^{p \bar{q}} R_{\bar{q} p} \nonumber\\
&=& {1 \over 2} R - {\alpha' \over 4} \sigma_2^{p \bar{q}} R_{\bar{q} p} +{\alpha' \over 2} \| \Omega \|^3 \tilde{\rho}^{p \bar{q}} R_{\bar{q} p} \nonumber\\
&=& {1 \over 2} R - {\alpha' \over 2} \sigma_2 (i {\rm Ric}_\o)+{\alpha' \over 2} \| \Omega \|^3 \tilde{\rho}^{p \bar{q}} R_{\bar{q} p}.
\eea
Therefore by (\ref{p_t-log|Omega|-from-eq})
\be \label{p_t-log|Omega|}
\p_t \log \| \Omega \| = {1 \over 2 \| \Omega \|} \bigg\{ \Delta_F \log \| \Omega \| - |T|^2  + {\alpha' \over 4} \sigma_2(i {\rm Ric}_\o) -  \| \Omega \|^2 \, \nu \bigg\}.
\ee
Substituting (\ref{p_t-|T|^2}) and (\ref{p_t-log|Omega|}) into (\ref{log_torsion_1}), we have
\bea \label{log_torsion_2}
\p_t G &=&   {1 \over 2 \| \Omega \|} \bigg\{ \Delta_F G +{|\na |T|^2 |^2_F \over |T|^4}- {|\na T |^2_{Fg}\over |T|^2} - {|\overline{\na} T |^2_{Fg} \over |T|^2}   -  {2 \over |T|^2} \Re \{ g^{i \bar{j}} \na_i |T|^2 \bar{T}_{\bar{j}}\}  \nonumber\\
&& - {1 \over 2} R  + {\alpha' \over 4} \sigma_2(i {\rm Ric}_\o)   +  {1 \over |T|^2} \Re \{F^{p \bar{q}} g^{i \bar{j}} T^\lambda{}_{pi} R_{\bar{q} \lambda} \bar{T}_{\bar{j}} \}  \nonumber\\
&& - {1 \over |T|^2} F^{p \bar{q}} g^{i \bar{j}}  T_i R_{\bar{q} p \bar{j}}{}^{\bar{\lambda}} \bar{T}_{\bar{\lambda}} + |T|^2  + {\alpha' \over 2} \| \Omega \|^3 \tilde{\rho}^{p \bar{q}} R_{\bar{q} p} - \| \Omega \|^2 \nu  \nonumber\\
&&+ {2 \over |T|^2} \Re \langle E, T \rangle  + \Lambda |T|^2 - {\alpha' \over 4} \Lambda \sigma_2(i {\rm Ric}_\o) +   \Lambda \| \Omega \|^2 \, \nu   \bigg\}.
\eea
Let $(p,t_0)$ be the point in $X \times [0,T]$ where $G$ attains its maximum. Since we start the flow at $t=0$ with a constant function $u_0=\log M$, the torsion is zero at the initial time. It follows that $t_0>0$. The following computation will be done at this point $(p,t_0)$, and we note that $|T|^2>0$ at $(p,t_0)$. The critical equation $\na G =0$ gives
\be \label{critical_eq}
0 = {\na_i |T|^2 \over |T|^2} - \Lambda T_i.
\ee
Using (\ref{na-T_j}), this can be rewritten in the following way 
\be
{\langle \na_i T, T \rangle \over |T|^2} = \Lambda T_i - {\langle T, \na_{\bar{i}} T \rangle \over |T|^2} = \Lambda T_i - {1 \over 2 |T|^2} g^{j \bar{k}} T_j R_{\bar{k} i}.
\ee
Therefore, by Cauchy-Schwarz and the critical equation,
\bea \label{cauchyschwarz}
- {|\na T|^2_{Fg} \over |T|^2} &\leq& - \bigg| {\langle \na T, T \rangle \over |T|^2} \bigg|^2_F= - \bigg| \Lambda T_i - {1 \over 2 |T|^2} g^{j \bar{k}} T_j R_{\bar{k} i} \bigg|^2_F \nonumber\\
&=& - \Lambda^2 |T|^2_F - {1 \over 4 |T|^4} \bigg|  g^{j \bar{k}} T_j R_{\bar{k} i} \bigg|^2_F +  {\Lambda \over |T|^2} \Re \{ F^{p \bar{q}} g^{j \bar{k}} T_j R_{\bar{k} p} \bar{T}_{\bar{q}} \}.
\eea
Here we used the notation $|V|^2_F = F^{p \bar{q}} V_p \bar{V}_{\bar{q}}$. We may also expand the following term using the definition of $F^{p \bar{q}}$,
\bea \label{expand_Ric^2}
4 |\overline{\na} T |^2_{Fg} &=&  F^{p \bar{q}} g^{i \bar{j}} R_{\bar{q} i} R_{\bar{j} p} = |{\rm Ric}_\o|^2 - {\alpha' \over 2} g^{i \bar{j}} \sigma_2^{p \bar{q}}R_{\bar{q} i} R_{\bar{j} p} + \alpha' \| \Omega \|^3 g^{i \bar{j}} \tilde{\rho} ^{p \bar{q}}R_{\bar{q} i} R_{\bar{j} p}.
\eea
Set $\e = 1/100$. Using (\ref{cauchyschwarz}) and (\ref{expand_Ric^2}), and the critical equation (\ref{critical_eq}) once more on the first and last term, we obtain
\bea
& \ & {|\na |T|^2 |^2_F \over |T|^4} - (1-\e){|\na T |^2_{Fg}\over |T|^2} - {|\overline{\na} T |^2_{Fg} \over |T|^2} -  {2 \over |T|^2} \Re \{ g^{i \bar{j}} \na_i |T|^2 \bar{T}_{\bar{j}}\} \nonumber\\
&\leq& \Lambda^2 |T|^2_F - (1-\e)\Lambda^2 |T|^2_F - (1-\e){1 \over 4 |T|^4} \bigg|  g^{j \bar{k}} T_j R_{\bar{k} i} \bigg|^2_F  \nonumber\\
&&+  (1-\e){\Lambda \over |T|^2} \Re \{ F^{p \bar{q}} g^{j \bar{k}} T_j R_{\bar{k} p} \bar{T}_{\bar{q}} \} - {1 \over 4} {|{\rm Ric}_\o|^2 \over |T|^2} + {\alpha' \over 8 |T|^2}  g^{i \bar{j}} \sigma_2^{p \bar{q}}R_{\bar{q} i} R_{\bar{j} p} \nonumber\\
&&- {\alpha' \over 4 |T|^2} \| \Omega \|^3 g^{i \bar{j}} \tilde{\rho} ^{p \bar{q}}R_{\bar{q} i} R_{\bar{j} p} - 2 \Lambda |T|^2.
\eea
Substituting this inequality into (\ref{log_torsion_2}), our main inequality becomes
\bea \label{log_torsion_3}
\p_t G &\leq&   {1 \over 2 \| \Omega \|} \bigg\{ \Delta_F G  -\e {|\na T |^2_{Fg}\over |T|^2}- {1 \over 4} {|{\rm Ric}_\o|^2 \over |T|^2} - (\Lambda -1) |T|^2 + \e \Lambda^2 |T|^2_F  - {1 \over 2} R  \nonumber\\
&& - {\alpha' \over 4} (\Lambda-1) \sigma_2(i {\rm Ric}_\o) + {\alpha' \over 8 |T|^2}   g^{i \bar{j}} \sigma_2^{p \bar{q}}R_{\bar{q} i} R_{\bar{j} p} +  (1-\e){\Lambda \over |T|^2} \Re \{ F^{p \bar{q}} g^{j \bar{k}} T_j R_{\bar{k} p} \bar{T}_{\bar{q}} \} \nonumber\\
&&    - {1 \over |T|^2} F^{p \bar{q}} g^{i \bar{j}}  T_i R_{\bar{q} p \bar{j}}{}^{\bar{\lambda}} \bar{T}_{\bar{\lambda}}+  {1 \over |T|^2} \Re \{ F^{p \bar{q}} g^{i \bar{j}} T^\lambda{}_{pi} R_{\bar{q} \lambda} \bar{T}_{\bar{j}} \}\nonumber\\
&&  - {(1 -\e) \over 4 |T|^4} \bigg|  g^{j \bar{k}} T_j R_{\bar{k} i} \bigg|^2_F - {\alpha' \over 4 |T|^2} \| \Omega \|^3 g^{i \bar{j}} \tilde{\rho} ^{p \bar{q}}R_{\bar{q} i} R_{\bar{j} p}+ {\alpha' \over 2} \| \Omega \|^3 \tilde{\rho}^{p \bar{q}} R_{\bar{q} p} \nonumber\\
&& + (\Lambda-1) \| \Omega \|^2 \nu  + {2 \over |T|^2} \Re \langle E, T \rangle   \bigg\},
\eea
which holds at $(p,t_0)$. Next, we use (\ref{curv_defn}) to write the evolving curvature as
\be
R_{\bar{q} p \bar{j}}{}^{\bar{\lambda}}  = \hat{R}_{\bar{q} p \bar{j}}{}^{\bar{\lambda}}+ {1 \over 2} R_{\bar{q} p} \delta_j{}^\lambda.
\ee
This identity allows us to write
\be
- {1 \over |T|^2} F^{p \bar{q}} g^{i \bar{j}}  T_i R_{\bar{q} p \bar{j}}{}^{\bar{\lambda}} T_{\bar{\lambda}}  = - {1 \over |T|^2} F^{p \bar{q}} g^{i \bar{j}}  T_i \hat{R}_{\bar{q} p \bar{j}}{}^{\bar{\lambda}} \bar{T}_{\bar{\lambda}} - {1 \over 2} F^{p \bar{q}} R_{\bar{q} p}.
\ee
Next, by (\ref{T^lambda}), the torsion can be written as
\be
 T^\lambda{}_{pi} =  T_i \delta^\lambda{}_p - T_p \delta^\lambda{}_i,
\ee
so we may rewrite
\bea
 {1 \over |T|^2} \Re \{ F^{p \bar{q}} g^{i \bar{j}} T^\lambda{}_{pi} R_{\bar{q} \lambda} \bar{T}_{\bar{j}} \} &=&  F^{p \bar{q}} R_{\bar{q} p} -   {1 \over |T|^2} \Re \{ F^{p \bar{q}} g^{i \bar{j}} R_{\bar{q} i} \bar{T}_{\bar{j}} T_p \}. 
\eea
Together, we have
\bea \label{log_torsion_4}
& \ & - {1 \over |T|^2} F^{p \bar{q}} g^{i \bar{j}}  T_i R_{\bar{q} p \bar{j}}{}^{\bar{\lambda}} \bar{T}_{\bar{\lambda}} + {1 \over |T|^2} \Re \{ F^{p \bar{q}} g^{i \bar{j}} T^\lambda{}_{pi} R_{\bar{q} \lambda} \bar{T}_{\bar{j}} \} \nonumber\\
&=& -{1 \over |T|^2} F^{p \bar{q}} g^{i \bar{j}}  T_i \hat{R}_{\bar{q} p \bar{j}}{}^{\bar{\lambda}} \bar{T}_{\bar{\lambda}} + {1 \over 2} F^{p \bar{q}} R_{\bar{q} p} - {1 \over |T|^2} \Re \{ F^{p \bar{q}} g^{i \bar{j}} R_{\bar{q} i} \bar{T}_{\bar{j}} T_p \} \nonumber\\
&=&  -{1 \over |T|^2} F^{p \bar{q}} g^{i \bar{j}}  T_i \hat{R}_{\bar{q} p \bar{j}}{}^{\bar{\lambda}} \bar{T}_{\bar{\lambda}} + {1 \over 2} R - {1 \over 2} \alpha' \sigma_2(i {\rm Ric}_\o) \nonumber\\
&& + {\alpha' \over 2} \| \Omega \|^3 \tilde{\rho}^{p \bar{q}} R_{\bar{q} p} -  {1 \over |T|^2} \Re \{ F^{p \bar{q}} g^{i \bar{j}} R_{\bar{q} i} \bar{T}_{\bar{j}} T_p \}. 
\eea
We also compute
\be \label{log_torsion_5}
|T|_F^2 = |T|^2 + \alpha' \| \Omega \|^3 \tilde{\rho}^{p \bar{q}} T_p T_{\bar{q}} - {\alpha' \over 2} \sigma_2^{p \bar{q}} T_p \bar{T}_{\bar{q}}.
\ee
Substituting (\ref{log_torsion_4}) and (\ref{log_torsion_5}) in the main inequality (\ref{log_torsion_3}), we see that the terms of order $R$ have cancelled. 
\bea 
\p_t G &\leq&   {1 \over 2 \| \Omega \|} \bigg\{ \Delta_F G  -\e {|\na T |^2_{Fg}\over |T|^2}- {1 \over 4} {|{\rm Ric}_\o|^2 \over |T|^2} - (\Lambda -1 - \e \Lambda^2) |T|^2   \nonumber\\
&&-  \e \Lambda^2 {\alpha' \over 2} \sigma_2^{p \bar{q}} T_p \bar{T}_{\bar{q}} - (1 + \Lambda) {\alpha' \over 4} \sigma_2(i {\rm Ric}_\o)  +{\alpha' \over 8 |T|^2}   g^{i \bar{j}} \sigma_2^{p \bar{q}}R_{\bar{q} i} R_{\bar{j} p}     \nonumber\\
&& +  (\Lambda-\e\Lambda -1){1 \over |T|^2} \Re \{ F^{p \bar{q}} g^{j \bar{k}} T_j R_{\bar{k} p} \bar{T}_{\bar{q}} \} - {(1 -\e) \over 4 |T|^4} \bigg|  g^{j \bar{k}} T_j R_{\bar{k} i} \bigg|^2_F \nonumber\\
&& - {1 \over |T|^2} F^{p \bar{q}} g^{i \bar{j}}  T_i \hat{R}_{\bar{q} p \bar{j}}{}^{\bar{\lambda}} \bar{T}_{\bar{\lambda}} - {\alpha' \over 4 |T|^2} \| \Omega \|^3 g^{i \bar{j}} \tilde{\rho} ^{p \bar{q}}R_{\bar{q} i} R_{\bar{j} p} + \e \Lambda^2 \alpha' \| \Omega \|^3 \tilde{\rho}^{p \bar{q}} T_p \bar{T}_{\bar{q}}  \nonumber\\
&&+ \alpha'  \| \Omega \|^3 \tilde{\rho}^{p \bar{q}} R_{\bar{q} p} + (\Lambda-1) \| \Omega \|^2 \nu + {2 \over |T|^2} \Re \langle E, T \rangle   \bigg\}.
\eea
We now substitute $\Lambda = 1 + {1 \over 8}$ and $\e = {1 \over 100}$. Then
\bea \label{log_torsion_6}
\p_t G &\leq&   {1 \over 2 \| \Omega \|} \bigg\{ \Delta_F G  -{1 \over 100} {|\na T |^2_{Fg}\over |T|^2}  - {1 \over 4} {|{\rm Ric}_\o|^2 \over |T|^2} - {1 \over 9} |T|^2   - \bigg( {9 \over 8}\bigg)^2 {1 \over 100} {\alpha' \over 2} \sigma_2^{p \bar{q}} T_p \bar{T}_{\bar{q}} \nonumber\\
&&- {17 \over 16} {\alpha' \over 2} \sigma_2(i {\rm Ric}_\o) + {\alpha' \over 8 |T|^2}   g^{i \bar{j}} \sigma_2^{p \bar{q}}R_{\bar{q} i} R_{\bar{j} p}  + \bigg({1 \over 8} - {9 \over 800} \bigg) {1 \over |T|^2} \Re \{ F^{p \bar{q}} g^{j \bar{k}} T_j R_{\bar{k} p} \bar{T}_{\bar{q}} \} \nonumber\\
&&- {99 \over 400} {1 \over |T|^4} \bigg|  g^{j \bar{k}} T_j R_{\bar{k} i} \bigg|^2_F  - {1 \over |T|^2} F^{p \bar{q}} g^{i \bar{j}}  T_i \hat{R}_{\bar{q} p \bar{j}}{}^{\bar{\lambda}} \bar{T}_{\bar{\lambda}}  - {\alpha' \over 4 |T|^2} \| \Omega \|^3 g^{i \bar{j}} \tilde{\rho} ^{p \bar{q}}R_{\bar{q} i} R_{\bar{j} p}\nonumber\\
&& + \alpha'  \| \Omega \|^3 \tilde{\rho}^{p \bar{q}} R_{\bar{q} p} + {1 \over 100} \bigg({9 \over 8} \bigg)^2 \alpha' \| \Omega \|^3 \tilde{\rho}^{p \bar{q}} T_p \bar{T}_{\bar{q}}   + {1 \over 8} \| \Omega \|^2 \nu + {2 \over |T|^2} \Re \langle E, T \rangle   \bigg\}
\eea
We are assuming in the hypothesis of Theorem \ref{torsion_estimate} that $|\alpha' {\rm Ric}_\o| < 10^{-6}$. By Theorem \ref{C0-thm}, we know that $\| \Omega \| \leq {C_2 \over M} \ll 1$, so for $M$ large enough we can assume
\be \label{g-onemillion-F}
(1-10^{-6}) g^{i \bar{j}} \leq F^{i \bar{j}} \leq (1+10^{-6}) g^{i \bar{j}}.
\ee
One way to see this inequality is by writing $F^{i \bar{j}}$ in coordinates (\ref{F^ii-in-coords}). Using (\ref{g-onemillion-F}), we can estimate
\bea
& \ & - {17 \over 16} {\alpha' \over 2} \sigma_2(i {\rm Ric}_\o) - {\alpha' \over 8 |T|^2}   g^{i \bar{j}} \sigma_2^{p \bar{q}}R_{\bar{q} i} R_{\bar{j} p} + \bigg({1 \over 8} - {9 \over 800} \bigg) {1 \over |T|^2} \Re\{ F^{p \bar{q}} g^{j \bar{k}} T_j R_{\bar{k} p} \bar{T}_{\bar{q}} \} \nonumber\\
&\leq& {17 \over 16} {1 \over 2} |\alpha' {\rm Ric}_\o| \, |{\rm Ric}_\o|   + {1 \over 8} | \alpha' {\rm Ric}_\o| \, {|{\rm Ric}_\o|^2 \over |T|^2}  + {1 \over 7}|{\rm Ric}_\o|  \nonumber\\
&\leq& {1 \over (2)(3)} |{\rm Ric}_\o| + {1 \over 100} {|{\rm Ric}_\o|^2 \over |T|^2} \nonumber\\
&\leq&  {1 \over (2)(3)^2} |T|^2 + \bigg({1 \over 100} + {1 \over (2)(2)^2}\bigg) {|{\rm Ric}_\o|^2 \over |T|^2}.
\eea
We also notice
\be
 - \bigg( {9 \over 8}\bigg)^2 {1 \over 100} {\alpha' \over 2} \sigma_2^{p \bar{q}} T_p \bar{T}_{\bar{q}} \leq |\alpha' {\rm Ric}| |T|^2 \leq {1 \over 10^6} |T|^2,
\ee
and
\be
 -{1 \over |T|^2} F^{p \bar{q}} g^{i \bar{j}}  T_i \hat{R}_{\bar{q} p \bar{j}}{}^{\bar{\lambda}} \bar{T}_{\bar{\lambda}} \leq C e^{-u} = C \| \Omega \|.
\ee
Substituting these estimates into (\ref{log_torsion_6}) gives
\bea
\p_t G &\leq&   {1 \over 2 \| \Omega \|} \bigg\{ \Delta_F G  -{1 \over 200} {|\na T |^2 \over |T|^2}  - {1 \over 100} {|{\rm Ric}_\o|^2 \over |T|^2} - {1 \over 100} |T|^2 \nonumber\\
&&  +C \| \Omega \|    + {\alpha' \over 4 |T|^2} \| \Omega \|^3 g^{i \bar{j}} \tilde{\rho} ^{p \bar{q}}R_{\bar{q} i} R_{\bar{j} p} + \alpha'  \| \Omega \|^3 \tilde{\rho}^{p \bar{q}} R_{\bar{q} p} \nonumber\\
&& + {1 \over 100} \bigg({9 \over 8} \bigg)^2 \alpha' \| \Omega \|^3 \tilde{\rho}^{p \bar{q}} T_p \bar{T}_{\bar{q}} + {1 \over 8} \| \Omega \|^2 \nu + {2 \over |T|^2} \Re \langle E, T \rangle   \bigg\}.
\eea
By the definition of $E$ (\ref{defn_E}) and $\nu$ (\ref{defn_nu}), the terms on the last two lines can only slightly perturb the coefficients of the first line since $\| \Omega \| = e^{-u} \leq {C_2 \over M} \ll 1$ for $M \gg 1$ large enough. We recall that $\tilde{\rho}^{p \bar{q}}$ and $b^i_{\rho}$ are bounded in $C^\infty$ in terms of the background metric $\hat{g}$, so for example,
\be \label{tilde-rho-bi}
\| \Omega \| \tilde{\rho}^{p \bar{q}} \leq C e^{-u} \hat{g}^{p \bar{q}} = C g^{p \bar{q}}, \ \ \| \Omega \|^{1/2} |b^i_{\rho} T_i| \leq C |T|.
\ee
This allows us to bound certain terms such as
\be
\alpha'  \| \Omega \|^3 \tilde{\rho}^{p \bar{q}} R_{\bar{q} p} \leq C \| \Omega \|^2 |{\rm Ric}_\o| \leq {C \over 2} \| \Omega \|^2 {|{\rm Ric}_\o|^2 \over |T|^2} + {C \over 2} \| \Omega \|^2 |T|^2,
\ee
and
\be
\alpha' \| \Omega \|^3 \Re \{ b^i_{\rho} T_i \} \leq C \| \Omega \|^2 |T| \leq C \| \Omega \|^2 {|T|^2 \over 2} + {C \over 2} \| \Omega \|^2.
\ee
Covariant derivatives with respect to the evolving metric act like $\na_i = \p_i - T_i$, so we can bound terms such as
\be
{2 \over |T|^2} {\alpha' \over 2} \| \Omega \|^3 g^{j \bar{k}} (\na_j \tilde{\rho}^{p \bar{q}}) R_{\bar{q} p} \bar{T}_{\bar{k}} \leq C\| \Omega \|^2 {|{\rm Ric}_\o| \over |T|} + C\| \Omega \|^2 {|{\rm Ric}_\o| \over |T|} |T|.
\ee
The inequality $2ab \leq a^2+b^2$ can be used to absorb terms into the first line. We also bound terms
\be
-{2 \over |T|^2} \| \Omega \|^2 g^{j \bar{k}} \na_j \tilde{\mu} \bar{T}_{\bar{k}} \leq C \| \Omega \|^2 {\| \Omega \|^{1/2} \over |T|}.
\ee
Using these estimates, it is possible to show that at the maximum point $(p,t_0)$ of $G$, for $\| \Omega \| \leq {C_2 \over M} \ll 1$, there holds
\be
0 \leq {1 \over 2 \| \Omega \|} \bigg\{ \Delta_F G - {1 \over 200} |T|^2 + C  \| \Omega \|  \bigg( 1 + {\| \Omega \|^{1/2} \over |T|} \bigg) \bigg\}.
\ee
By (\ref{g-onemillion-F}), $\Delta_F G \leq 0$ at the maximum $(p,t_0)$ of $G$, hence
\be
|T|^2 \leq C \| \Omega \| \leq  {C \over M}.
\ee
Therefore
\be
G \leq G(p,t_0) \leq \log {C \over M} + \Lambda u(p).
\ee
By Theorem \ref{C0-thm},
\bea
|T|^2 &\leq& {C \over M} \exp \{ \Lambda (u(p)-u) \} \nonumber\\
&\leq& {C \over M} \bigg(\sup_{X \times [0,T)} e^{u} \bigg)^\Lambda \bigg(\sup_{X \times [0,T)} e^{-u} \bigg)^\Lambda \nonumber\\
&\leq&  {C \over M} (C_2 C_1)^\Lambda \ll 1.
\eea
This proves Theorem \ref{torsion_estimate}.

\

\section{Evolution of the curvature} \label{section_curvature}
\setcounter{equation}{0}
\subsection{Ricci curvature}
In this subsection, we flow the Ricci curvature of the evolving Hermitian metric $e^u \hat{g}$. We will use the well-known general formula for the evolution of the curvature tensor
\be
\p_t R_{\bar{k} j}{}^\alpha{}_\beta = - \na_{\bar{k}} \na_j (g^{\alpha \bar{\gamma}} \p_t g_{\bar{\gamma} \beta}).
\ee
Recall that we defined $R_{\bar{k} j} = R_{\bar{k} j}{}^\alpha{}_\alpha$, hence substituting (\ref{Flow3}) yields
\be
\p_t R_{\bar{k} j} =   - \na_{\bar{k}} \na_j \bigg\{ {1 \over 2 \| \Omega \|} \bigg( - R  - \alpha' \| \Omega \|^3 \tilde{\rho}^{p \bar{q}} R_{\bar{q} p}  + {\alpha' \over 2} \sigma_2(i {\rm Ric}_\o) +  2 |T|^2 + 2 \| \Omega \|^2 \nu   \bigg) \bigg\}.
\ee
Expanding out terms gives
\bea
\p_t R_{\bar{k} j} &=&  {1 \over 2 \| \Omega \|} \bigg\{  \na_{\bar{k}} \na_j R   + \alpha' \| \Omega \|^3 \tilde{\rho}^{p \bar{q}} \na_{\bar{k}} \na_j R_{\bar{q} p}   - \na_{\bar{k}} \na_j{\alpha' \over 2} \sigma_2(i {\rm Ric}_\o) - 2 \na_{\bar{k}} \na_j|T|^2 \nonumber\\
&&+ \alpha' \na_{\bar{k}}( \| \Omega \|^3 \tilde{\rho}^{p \bar{q}})  \na_j R_{\bar{q} p} + \alpha' \na_j ( \| \Omega \|^3 \tilde{\rho}^{p \bar{q}}) \na_{\bar{k}} R_{\bar{q} p}
 \nonumber\\
&& + \alpha' \na_{\bar{k}} \na_j ( \| \Omega \|^3 \tilde{\rho}^{p \bar{q}})  R_{\bar{q} p} - \na_{\bar{k}} \na_j 2 \| \Omega \|^2 \nu \bigg\}  \nonumber\\
&& -{\na_j \| \Omega \| \over 2 \| \Omega \|^2} \na_{\bar{k}} \bigg\{    R  + \alpha'( \| \Omega \|^3 \tilde{\rho}^{p \bar{q}}  R_{\bar{q} p})   - {\alpha' \over 2} \sigma_2(i {\rm Ric}_\o) - 2 |T|^2 - 2 \| \Omega \|^2 \nu \bigg\} \nonumber\\
&& -{\na_{\bar{k}} \| \Omega \| \over 2 \| \Omega \|^2} \na_{j}  \bigg\{  R  + \alpha' ( \| \Omega \|^3 \tilde{\rho}^{p \bar{q}}  R_{\bar{q} p})   -{\alpha' \over 2} \sigma_2(i {\rm Ric}_\o) -2  |T|^2 -2 \| \Omega \|^2 \nu \bigg\}\nonumber\\
&& +\bigg\{ {-\na_{\bar{k}}\na_j \| \Omega \| \over 2 \| \Omega \|^2} + {2 \na_{\bar{k}}\| \Omega \| \na_j \| \Omega \| \over 2 \| \Omega \|^3} \bigg\} \bigg\{  R   + \alpha' \| \Omega \|^3 \tilde{\rho}^{p \bar{q}} R_{\bar{q} p} \nonumber\\
&&  - {\alpha' \over 2} \sigma_2(i {\rm Ric}_\o) - 2|T|^2 - 2 \| \Omega \|^2 \nu \bigg\}.
\eea
Using $\na_j \| \Omega \| = \| \Omega \| T_j$, 
\bea \label{p_t-Ric1}
\p_t R_{\bar{k} j} &=&  {1 \over 2 \| \Omega \|} \bigg\{  \na_{\bar{k}} \na_j R   + \alpha' \| \Omega \|^3 \tilde{\rho}^{p \bar{q}} \na_{\bar{k}} \na_j R_{\bar{q} p}  - \na_{\bar{k}} \na_j{\alpha' \over 2} \sigma_2(i {\rm Ric}_\o) \nonumber\\
&&- 2 \na_{\bar{k}} \na_j|T|^2 + \alpha' \na_{\bar{k}}( \| \Omega \|^3 \tilde{\rho}^{p \bar{q}})  \na_j R_{\bar{q} p} + \alpha' \na_j ( \| \Omega \|^3 \tilde{\rho}^{p \bar{q}}) \na_{\bar{k}} R_{\bar{q} p} \nonumber\\
&& + \alpha' \na_{\bar{k}} \na_j ( \| \Omega \|^3 \tilde{\rho}^{p \bar{q}})  R_{\bar{q} p}   - 2\na_{\bar{k}} \na_j \Big\{  \| \Omega \|^2  \nu \Big\}  - T_j  \na_{\bar{k}}  R  \nonumber\\
&&- \alpha' T_j \na_{\bar{k}}( \| \Omega \|^3 \tilde{\rho}^{p \bar{q}}  R_{\bar{q} p}) + 2 T_j \na_{\bar{k}} |T|^2  + {\alpha' \over 2} T_j \na_{\bar{k}} \left(\sigma_2(i {\rm Ric}_\o)\right) + 2 T_j \na_{\bar{k}} \Big\{ \| \Omega \|^2 \nu \Big\}   \nonumber\\
&& - T_{\bar{k}}  \na_{j}  R  - \alpha' T_{\bar{k}} \na_{j}( \| \Omega \|^3 \tilde{\rho}^{p \bar{q}}  R_{\bar{q} p}) +2 T_{\bar{k}} \na_j |T|^2   + {\alpha' \over 2}T_{\bar{k}}\na_j\left( \sigma_2(i {\rm Ric}_\o) \right)\nonumber\\
&& + 2 T_{\bar{k}}\na_j  \Big\{ \| \Omega \|^2 \nu \Big\} +  R T_j T_{\bar{k}}  + \alpha' T_j T_{\bar{k}}( \| \Omega \|^3 \tilde{\rho}^{p \bar{q}}  R_{\bar{q} p}) - 2 |T|^2 T_j T_{\bar{k}}  \nonumber\\
&& - {\alpha' \over 2} \sigma_2(i {\rm Ric}_\o)T_j T_{\bar{k}} - 2 T_j T_{\bar{k}} \Big\{ \| \Omega \|^2 \nu \Big\} -  R \na_{\bar{k}} T_j - \alpha' \na_{\bar{k}} T_j ( \| \Omega \|^3 \tilde{\rho}^{p \bar{q}}  R_{\bar{q} p})  \nonumber\\
&&+ 2 |T|^2 \na_{\bar{k}}T_j + {\alpha' \over 2} \sigma_2(i {\rm Ric}_\o)\na_{\bar{k}}T_j  + 2 \na_{\bar{k}}T_j \Big\{ \| \Omega \|^2 \nu \Big\}  \bigg\}.
\eea
We now study the highest order terms, namely
\be
\na_{\bar{k}} \na_j R_{\bar{q} p}  =  - 2 \na_{\bar{k}} \na_j \na_p \na_{\bar{q}} u.
\ee
We will use the following commutation formula for covariant derivatives in Hermitian geometry
\bea
\na_{\bar{k}} \na_j \na_p \na_{\bar{q}} u &=& \na_p \na_{\bar{q}} \na_j \na_{\bar{k}} u + T^\lambda{}_{pj} \na_{\bar{q}} \na_\lambda \na_{\bar{k}} u  + \bar{T}^{\bar{\lambda}}{}_{\bar{q} \bar{k}}   \na_p  \na_j \na_{\bar{\lambda}} u \nonumber\\
&& + R_{\bar{k} j}{}^{\lambda}{}_p u_{\bar{q} \lambda}  - R_{\bar{q} p \bar{k}}{}^{\bar{\lambda}} u_{\bar{\lambda} j} + \bar{T}^{\bar{\lambda}}{}_{\bar{q} \bar{k}}  T^\gamma{}_{pj}  u_{\bar{\lambda} \gamma}.
\eea
Using $R_{\bar{q} p} = - 2 u_{\bar{q} p}$, we obtain
\bea
\na_{\bar{k}} \na_j R_{\bar{q} p} &=& \na_p \na_{\bar{q}} R_{\bar{k} j} + T^\lambda{}_{pj} \na_{\bar{q}} R_{\bar{k} \lambda}  + \bar{T}^{\bar{\lambda}}{}_{\bar{q} \bar{k}}   \na_p  R_{\bar{\lambda} j} \nonumber\\
&& + R_{\bar{k} j}{}^{\lambda}{}_p R_{\bar{q} \lambda}  - R_{\bar{q} p \bar{k}}{}^{\bar{\lambda}} R_{\bar{\lambda} j} + \bar{T}^{\bar{\lambda}}{}_{\bar{q} \bar{k}}  T^\gamma{}_{pj}  R_{\bar{\lambda} \gamma}.
\eea
Hence
\bea \label{nana-R}
\na_{\bar{k}} \na_j R &=& g^{p \bar{q}}  \na_p \na_{\bar{q}} R_{\bar{k} j} + g^{p \bar{q}} T^\lambda{}_{pj} \na_{\bar{q}} R_{\bar{k} \lambda} +  g^{p \bar{q}} \bar{T}^{\bar{\lambda}}{}_{\bar{q} \bar{k}} \na_{p} R_{\bar{\lambda} j}\nonumber\\
&&+ R_{\bar{k} j}{}^{p \bar{q}} R_{\bar{q} p} - R^p{}_{p \bar{k}}{}^{\bar{\lambda}} R_{\bar{\lambda} j} + g^{p \bar{q}} \bar{T}^{\bar{\lambda}}{}_{\bar{q} \bar{k}}  T^\gamma{}_{pj}  R_{\bar{\lambda} \gamma}.
\eea
Differentiating $\sigma_2^{p \bar{q}}$ (\ref{sigma2^pq-defn}) leads to the following definition
\be \label{sigma_2-pqrs}
\sigma_2^{p \bar{q}, r \bar{s}} = g^{p \bar{q}} g^{r \bar{s}} - g^{p \bar{s}} g^{r \bar{q}}. 
\ee
With this notation, we now differentiate $\sigma_2(i {\rm Ric}_\o)$ twice.
\bea \label{nana-sigma2}
\na_{\bar{k}} \na_j \sigma_2(i {\rm Ric}_\o) &=& \na_{\bar{k}} (\sigma_2^{p \bar{q}} \na_j R_{\bar{q} p}) \nonumber\\
&=& \sigma_2^{p \bar{q}}  \na_{\bar{k}} \na_j R_{\bar{q} p} + \sigma_2^{p \bar{q}, r \bar{s}} \na_{\bar{k}} R_{\bar{s} r} \na_j R_{\bar{q} p} \nonumber\\
&=&  \sigma_2^{p \bar{q}} \na_p \na_{\bar{q}} R_{\bar{k} j} + \sigma_2^{p \bar{q}, r \bar{s}} \na_{\bar{k}} R_{\bar{s} r} \na_j R_{\bar{q} p} + \sigma_2^{p \bar{q}} T^\lambda{}_{pj} \na_{\bar{q}} R_{\bar{k} \lambda}  \nonumber\\
&& + \sigma_2^{p \bar{q}} \bar{T}^{\bar{\lambda}}{}_{\bar{q} \bar{k}}  \na_p  R_{\bar{\lambda} j}  + \sigma_2^{p \bar{q}} R_{\bar{k} j}{}^{\lambda}{}_p R_{\bar{q} \lambda}  -\sigma_2^{p \bar{q}}  R_{\bar{q} p \bar{k}}{}^{\bar{\lambda}} R_{\bar{\lambda} j} \nonumber\\
&&+ \sigma_2^{p \bar{q}} \bar{T}^{\bar{\lambda}}{}_{\bar{q} \bar{k}}  T^\gamma{}_{pj}  R_{\bar{\lambda} \gamma}.
\eea
By (\ref{nana-R}) and (\ref{nana-sigma2}), and proceeding similarly for the $\rho$ term, we obtain
\bea \label{nana_R-sigma2Ric}
& \ & \na_{\bar{k}} \na_j R  + \alpha' \| \Omega \|^3 \tilde{\rho}^{p \bar{q}} \na_{\bar{k}} \na_j R_{\bar{q} p} - {\alpha' \over 2} \na_{\bar{k}} \na_j \sigma_2(i {\rm Ric}_\o) \nonumber\\
&=&  F^{p \bar{q}} \na_p \na_{\bar{q}} R_{\bar{k} j} -{\alpha' \over 2} \sigma_2^{p \bar{q}, r \bar{s}} \na_{\bar{k}} R_{\bar{s} r} \na_j R_{\bar{q} p} + F^{p \bar{q}} T^\lambda{}_{pj} \na_{\bar{q}} R_{\bar{k} \lambda} + F^{p \bar{q}} \bar{T}^{\bar{\lambda}}{}_{\bar{q} \bar{k}}  \na_p  R_{\bar{\lambda} j}   \nonumber\\
&&+ F^{p \bar{q}} R_{\bar{k} j}{}^{\lambda}{}_p R_{\bar{q} \lambda}  -F^{p \bar{q}}  R_{\bar{q} p \bar{k}}{}^{\bar{\lambda}} R_{\bar{\lambda} j} + F^{p \bar{q}} \bar{T}^{\bar{\lambda}}{}_{\bar{q} \bar{k}}  T^\gamma{}_{pj}  R_{\bar{\lambda} \gamma},
\eea
where the definition of $F^{p \bar{q}}$ was given in (\ref{F^pq-defn}).
\par Using (\ref{na-T_j}), we may convert derivatives of torsion $\overline{\na} T$ into curvature terms, but terms $\na T$ are of different type and must be treated separately. For example
\bea \label{na-barna-|T|^2}
- 2 \na_{\bar{k}} \na_j|T|^2 &=& - 2 g^{p \bar{q}} \na_{\bar{k}} \na_j T_p \bar{T}_{\bar{q}} - 2 g^{p \bar{q}} \na_j T_p \na_{\bar{k}} \bar{T}_{\bar{q}} - {1 \over 2} g^{p \bar{q}} R_{\bar{k} p} R_{\bar{q} j}- g^{p \bar{q}}   T_p \na_{\bar{k}} R_{\bar{q} j} \nonumber\\
&=&  - g^{p \bar{q}}  \na_j R_{\bar{k} p} \bar{T}_{\bar{q}} - 2 g^{p \bar{q}} R_{\bar{k} j}{}^\lambda{}_p T_\lambda \bar{T}_{\bar{q}}  - 2 g^{p \bar{q}} \na_j T_p \na_{\bar{k}} \bar{T}_{\bar{q}} \nonumber\\
&& - {1 \over 2} g^{p \bar{q}} R_{\bar{k} p} R_{\bar{q} j} -  g^{p \bar{q}}   T_p \na_{\bar{k}} R_{\bar{q} j} .
\eea
Substituting (\ref{nana_R-sigma2Ric}) and (\ref{na-barna-|T|^2}) into (\ref{p_t-Ric1}),
\bea \label{p_t-Ric2}
\p_t R_{\bar{k} j} &=&  {1 \over 2 \| \Omega \|} \bigg\{  F^{p \bar{q}}  \na_p \na_{\bar{q}} R_{\bar{k} j} -{\alpha' \over 2} \sigma_2^{p \bar{q}, r \bar{s}} \na_{\bar{k}} R_{\bar{s} r} \na_j R_{\bar{q} p} \nonumber\\
&&+ 2 \alpha'  \| \Omega \|^3  \tilde{\rho}^{p \bar{q}}  \na_j T_p \na_{\bar{k}} \bar{T}_{\bar{q}} - 2 g^{p \bar{q}} \na_j T_p \na_{\bar{k}} \bar{T}_{\bar{q}} + Y_{\bar{k} j}  \bigg\}.
\eea
where $Y_{\bar{k} j}$ contains various combinations of torsion and curvature terms, but is linear in first derivatives of curvature and torsion and does not contain higher order derivatives of curvature and torsion. Explicitly,
\bea \label{Y-full}
Y_{\bar{k} j} &=& F^{p \bar{q}} T^\lambda{}_{pj} \na_{\bar{q}} R_{\bar{k} \lambda} + F^{p \bar{q}} \bar{T}^{\bar{\lambda}}{}_{\bar{q} \bar{k}}  \na_p  R_{\bar{\lambda} j}   + F^{p \bar{q}} R_{\bar{k} j}{}^{\lambda}{}_p R_{\bar{q} \lambda}  -F^{p \bar{q}}  R_{\bar{q} p \bar{k}}{}^{\bar{\lambda}} R_{\bar{\lambda} j}  \nonumber\\
&& + F^{p \bar{q}} \bar{T}^{\bar{\lambda}}{}_{\bar{q} \bar{k}}  T^\gamma{}_{pj}  R_{\bar{\lambda} \gamma}  - g^{p \bar{q}}  \na_j R_{\bar{k} p} \bar{T}_{\bar{q}} - 2 g^{p \bar{q}} R_{\bar{k} j}{}^\lambda{}_p T_\lambda \bar{T}_{\bar{q}}   - {1 \over 2} g^{p \bar{q}} R_{\bar{k} p} R_{\bar{q} j}  \nonumber\\
&&-  g^{p \bar{q}}   T_p \na_{\bar{k}} R_{\bar{q} j}  + \alpha' \na_{\bar{k}}( \| \Omega \|^3 \tilde{\rho}^{p \bar{q}})  \na_j R_{\bar{q} p} + \alpha' \na_j ( \| \Omega \|^3 \tilde{\rho}^{p \bar{q}}) \na_{\bar{k}} R_{\bar{q} p}  \nonumber\\
 &&+ \alpha'  (\na_{\bar{k}} \na_j \| \Omega \|^3 \tilde{\rho}^{p \bar{q}})  R_{\bar{q} p} - T_j F^{p \bar{q}} \na_{\bar{k}} R_{\bar{q} p} - \alpha' T_j \na_{\bar{k}}( \| \Omega \|^3 \tilde{\rho}^{p \bar{q}})  R_{\bar{q} p}+  g^{p \bar{q}} T_j R_{\bar{k} p} \bar{T}_{\bar{q}}  \nonumber\\
&& + 2 g^{p \bar{q}} T_j T_p  \na_{\bar{k}} \bar{T}_{\bar{q}} - 2  \bigg\{ -\alpha' \na_{\bar{k}} \na_j (\| \Omega \|^3 \psi_\rho) - {\alpha' \over 2}  \Re \{ \| \Omega \|^3 b^i_\rho  \na_j  R_{\bar{k} i} \}  \nonumber\\
&& - \alpha'  \Re \{ \| \Omega \|^3 b^i_\rho  R_{\bar{k} j}{}^\lambda{}_i T_\lambda \} - \alpha'  \Re \{\na_{\bar{k}}( \| \Omega \|^3 b^i_\rho)  \na_j T_i \}  \nonumber\\
&& - \alpha'  \Re \{\na_{j}( \| \Omega \|^3 b^i_\rho)  \na_{\bar{k}} T_i \}  - \alpha'  \Re \{\na_{\bar{k}} \na_j (\| \Omega \|^3 b^i_\rho) T_i \}   +\na_{\bar{k}} \na_j ( \| \Omega \|^2 \tilde{\mu}) \bigg\}  \nonumber\\
&& + \Bigg\{ 2 \alpha' (\na_{\bar{k}} \na_j \| \Omega \|^3  \tilde{\rho}^{p \bar{q}}) T_p \bar{T}_{\bar{q}} + 2 \alpha' \na_{\bar{k}}  ( \| \Omega \|^3  \tilde{\rho}^{p \bar{q}} )\na_j (T_p \bar{T}_{\bar{q}})  \nonumber\\
&&+ 2 \alpha'  \na_j  (\| \Omega \|^3  \tilde{\rho}^{p \bar{q}}) \na_{\bar{k}} (T_p \bar{T}_{\bar{q}}) +  \alpha'  \| \Omega \|^3  \tilde{\rho}^{p \bar{q}}\na_j R_{\bar{k} p} \bar{T}_{\bar{q}}  + 2 \alpha'  \| \Omega \|^3  \tilde{\rho}^{p \bar{q}} R_{\bar{k} j}{}^\lambda{}_p T_\lambda \bar{T}_{\bar{q}} \nonumber\\
&& +  \alpha'  \| \Omega \|^3  \tilde{\rho}^{p \bar{q}}  T_p \na_{\bar{k}} R_{\bar{q} j} + {\alpha' \over 2}  \| \Omega \|^3  \tilde{\rho}^{p \bar{q}} R_{\bar{k} p} R_{\bar{q} j} \bigg\} \nonumber\\
&&+ 2 T_j \na_{\bar{k}} \bigg\{ -\alpha' \| \Omega \|^3 \psi_\rho - \alpha' \| \Omega \|^3 \Re \{ b^i_\rho T_i \}  - \alpha'\| \Omega \|^3  \tilde{\rho}^{p \bar{q}} T_p \bar{T}_{\bar{q}}  +  \| \Omega \|^2 \tilde{\mu} \bigg\} \nonumber\\
&& - T_{\bar{k}}  F^{p \bar{q}} \na_{j}  R_{\bar{q} p}  - \alpha' T_{\bar{k}} \na_{j}( \| \Omega \|^3 \tilde{\rho}^{p \bar{q}})  R_{\bar{q} p}  +2 g^{p \bar{q}} T_{\bar{k}} \na_j T_p \bar{T}_{\bar{q}} + g^{p \bar{q}} T_{\bar{k}}  T_p R_{\bar{q} j}  \nonumber\\
&& + 2 T_{\bar{k}}\na_j \bigg\{ -\alpha' \| \Omega \|^3 \psi_\rho - \alpha' \| \Omega \|^3 \Re \{ b^i_\rho T_i \}  - \alpha'\| \Omega \|^3  \tilde{\rho}^{p \bar{q}} T_p \bar{T}_{\bar{q}}  +  \| \Omega \|^2 \tilde{\mu} \bigg\} \nonumber\\
&& +  R T_j T_{\bar{k}}  + \alpha' T_j T_{\bar{k}}( \| \Omega \|^3 \tilde{\rho}^{p \bar{q}}  R_{\bar{q} p}) - 2 |T|^2 T_j T_{\bar{k}} - {\alpha' \over 2} \sigma_2(i {\rm Ric}_\o)T_j T_{\bar{k}} \nonumber\\
&& - 2  T_j T_{\bar{k}} \bigg\{ -\alpha' \| \Omega \|^3 \psi_\rho - \alpha' \| \Omega \|^3 \Re \{ b^i_\rho T_i \}  - \alpha'\| \Omega \|^3  \tilde{\rho}^{p \bar{q}} T_p \bar{T}_{\bar{q}}  +  \| \Omega \|^2 \tilde{\mu} \bigg\} \nonumber\\
&& -  {1 \over 2} R R_{\bar{k} j} - {\alpha' \over 2}  R_{\bar{k} j} ( \| \Omega \|^3 \tilde{\rho}^{p \bar{q}}  R_{\bar{q} p}) +  |T|^2 R_{\bar{k} j} + {\alpha' \over 4} \sigma_2(i {\rm Ric}_\o)  R_{\bar{k} j}  \nonumber\\
&&+ R_{\bar{k} j} \bigg\{ -\alpha' \| \Omega \|^3 \psi_\rho - \alpha' \| \Omega \|^3 \Re \{ b^i_\rho T_i \}  - \alpha'\| \Omega \|^3  \tilde{\rho}^{p \bar{q}} T_p \bar{T}_{\bar{q}}  +  \| \Omega \|^2 \tilde{\mu} \bigg\}.
\eea
The terms in brackets indicate terms which come from substituting the definition of $\nu$ (\ref{defn_nu}).

\subsection{Evolving the norm of the curvature}
We will compute
\be
\p_t | {\rm Ric}_\o |^2 = \p_t \{  g^{k \bar{\ell}}g^{i \bar{j}} R_{\bar{\ell} i} \overline{R_{\bar{k} j}} \}.
\ee
We have
\bea
\p_t g^{i \bar{j}} &=& - g^{i \bar{\lambda}} g^{\gamma \bar{j}} \p_t g_{\bar{\lambda} \gamma} \nonumber\\
 &=&  {1 \over 2 \| \Omega \|} \bigg( {R \over 2} +{\alpha' \over 2} \| \Omega \|^3 \tilde{\rho}^{p \bar{q}} R_{\bar{q} p} - {\alpha' \over 4} \sigma_2(i {\rm Ric}_\o)-  |T|^2  - \| \Omega \|^2 \, \nu  \bigg) \, g^{i\bar{j} }.
\eea
Hence
\bea \label{p_t-|Ric|1}
\p_t | {\rm Ric}_\o |^2 &=& 2 \Re \langle \p_t {\rm Ric}_\o, {\rm Ric}_\o \rangle \nonumber\\
&&+  {| {\rm Ric}_\o |^2 \over 2 \| \Omega \|} \bigg( R  + \alpha' \| \Omega \|^3 \tilde{\rho}^{p \bar{q}} R_{\bar{q} p}  - {\alpha' \over 2} \sigma_2(i {\rm Ric}_\o)-  2|T|^2 - 2 \| \Omega \|^2 \, \nu  \bigg).
\eea
Next,
\bea
F^{p \bar{q}} \na_p \na_{\bar{q}} | {\rm Ric}_\o |^2 &=&  g^{k \bar{\ell}}g^{i \bar{j}} F^{p \bar{q}} \na_p \na_{\bar{q}} R_{\bar{\ell} i} R_{\bar{j} k} +  g^{k \bar{\ell}}g^{i \bar{j}} R_{\bar{\ell} i}  F^{p \bar{q}}  \na_p \na_{\bar{q}} R_{\bar{j} k} \nonumber\\
&& + |\na {\rm Ric}_\o |^2_{Fgg} + |\overline{\na} {\rm Ric}_\o |^2_{Fgg} \nonumber\\
&=& g^{k \bar{\ell}}g^{i \bar{j}} F^{p \bar{q}} \na_p \na_{\bar{q}} R_{\bar{\ell} i} \overline{R_{\bar{k} j}} +  g^{k \bar{\ell}}g^{i \bar{j}} R_{\bar{\ell} i} \overline{ F^{q \bar{p}} \na_q \na_{\bar{p}} R_{\bar{k} j}}\nonumber\\
&& -  g^{k \bar{\ell}}g^{i \bar{j}} R_{\bar{\ell} i}  F^{p \bar{q}} R_{\bar{q} p}{}^{\lambda}{}_k R_{\bar{j} \lambda} +  g^{k \bar{\ell}}g^{i \bar{j}} R_{\bar{\ell} i}  F^{p \bar{q}} R_{\bar{q} p \bar{j}}{}^{\bar{\lambda}} R_{\bar{\lambda} k}  \nonumber\\
&& + |\na {\rm Ric}_\o |^2_{Fgg} + |\overline{\na} {\rm Ric}_\o |^2_{Fgg}.
\eea
We have shown
\bea
\Delta_F | {\rm Ric}_\o |^2 &=&  2 \Re \langle \Delta_F {\rm Ric}_\o,{\rm Ric}_\o \rangle  + |\na {\rm Ric}_\o |^2_{Fgg} + |\overline{\na} {\rm Ric}_\o |^2_{Fgg} \nonumber\\
&& -  g^{k \bar{\ell}}g^{i \bar{j}} R_{\bar{\ell} i}  F^{p \bar{q}} R_{\bar{q} p}{}^{\lambda}{}_k R_{\bar{j}\lambda }+  g^{k \bar{\ell}}g^{i \bar{j}} R_{\bar{\ell} i}  F^{p \bar{q}} R_{\bar{q} p \bar{j}}{}^{\bar{\lambda}} R_{\bar{\lambda} k} .
\eea
Substituting (\ref{p_t-Ric2}) into (\ref{p_t-|Ric|1}) gives
\bea \label{p_t-Ric-full}
\p_t | {\rm Ric}_\o |^2 &=&  {1 \over 2 \| \Omega \|} \bigg\{ \Delta_F | {\rm Ric}_\o |^2  - |\na {\rm Ric}_\o |^2_{Fgg} - |\overline{\na} {\rm Ric}_\o |^2_{Fgg} \\
&&  -\alpha' \Re \{ g^{j \bar{\ell}} g^{m \bar{k}} \sigma_2^{p \bar{q}, r \bar{s}}R_{\bar{\ell} m} \na_{\bar{k}} R_{\bar{s} r} \na_j R_{\bar{q} p} \} \nonumber\\
&& + 4 \alpha' \Re \{  g^{j \bar{\ell}} g^{m \bar{k}}R_{\bar{\ell} m} \| \Omega \|^3  \tilde{\rho}^{p \bar{q}}  \na_j T_p \na_{\bar{k}} \bar{T}_{\bar{q}} \}\nonumber\\
&& - 4 \Re \{ g^{j \bar{\ell}} g^{m \bar{k}}R_{\bar{\ell} m}g^{p \bar{q}} \na_j T_p \na_{\bar{k}} \bar{T}_{\bar{q}} \} + 2 \Re \{ g^{j \bar{\ell}} g^{m \bar{k}}R_{\bar{\ell} m}Y_{\bar{k} j} \} \nonumber\\
&&+  g^{k \bar{\ell}}g^{i \bar{j}} R_{\bar{\ell} i}  F^{p \bar{q}} R_{\bar{q} p}{}^{\lambda}{}_k R_{\bar{j}\lambda }-  g^{k \bar{\ell}}g^{i \bar{j}} R_{\bar{\ell} i}  F^{p \bar{q}} R_{\bar{q} p \bar{j}}{}^{\bar{\lambda}} R_{\bar{\lambda} k} + | {\rm Ric}_\o |^2R \nonumber\\
&&+ \alpha' \| \Omega \|^3 \tilde{\rho}^{p \bar{q}} R_{\bar{q} p}| {\rm Ric}_\o |^2 - 2 |T|^2 | {\rm Ric}_\o |^2  - {\alpha' \over 2} \sigma_2(i {\rm Ric}_\o)| {\rm Ric}_\o |^2 \nonumber\\
&& - 2 \| \Omega \|^2 | {\rm Ric}_\o |^2\, \nu \bigg\}. \nonumber
\eea

\

\subsection{Estimating Ricci curvature}
\begin{lemma} \label{Ric-delta-epsilon-lemma}
Let $0<\delta, \epsilon < {1 \over 2}$ be such that $- {1 \over 4} g^{p \bar{q}} < \alpha' \delta^2 \| \Omega \| \tilde{\rho}^{p \bar{q}} < {1 \over 4} g^{p \bar{q}}$, and
\be
\| \Omega \|^2 \leq \delta, \ \ |T|^2 \leq \delta, \ \  |\alpha' {\rm Ric}_\o| \leq \epsilon,
\ee
at a point $(p,t_0)$. Let $\Lambda >1$ be any constant. Then at $(p,t_0)$ there holds
\bea \label{lemma-Ric-estimate}
&&\p_t ( |\alpha' {\rm Ric}_\o |^2+ \Lambda |T|^2) \nonumber
\\
&\leq&  {1 \over 2 \| \Omega \|} \bigg\{ \Delta_F ( |\alpha' {\rm Ric}_\o |^2 + \Lambda |T|^2) - \bigg( {1 \over 2} - 2 \epsilon \bigg) | \alpha' \na {\rm Ric}_\o |^2 \\
&&-\bigg( {\Lambda \over 4}  - (5+C\delta^2) \epsilon |\alpha'|^{-1}  \bigg) |\na T|^2 - {\Lambda \over 8} | {\rm Ric}_\o |^2 + C(1+\Lambda)\epsilon \delta  + C \epsilon^2 + C \Lambda \delta \bigg\}, \nonumber
\eea
for some constant $C$ only depending on $\tilde{\mu}$, $\rho$, $\alpha'$, and the background manifold $(X,\hat{\o})$.
\end{lemma}
{\it Proof:} Since $\epsilon$ and $\delta$ are assumed to be small, we have
\be
F^{p \bar{q}} = g^{p \bar{q}} + \alpha' \| \Omega \|^3 \tilde{\rho}^{p \bar{q}} - {\alpha' \over 2} \sigma_2^{p \bar{q}}, \ \   {1 \over 2} g^{p \bar{q}} < F^{p \bar{q}}< {3 \over 2} g^{p \bar{q}}.
\ee
We note the following estimate
\be \label{sigma2-pqrs-est1}
 -\alpha' \Re \{g^{j \bar{\ell}} g^{m \bar{k}} \sigma_2^{p \bar{q}, r \bar{s}}R_{\bar{\ell} m} \na_{\bar{k}} R_{\bar{s} r} \na_j R_{\bar{q} p} \} \leq  |\alpha' {\rm Ric}_\o| \,  |\na {\rm Ric}_\o|^2 .
\ee
We will estimate and group terms in (\ref{p_t-Ric-full}) and (\ref{Y-full}). We will convert $F^{p\bar q}$ into the metric $g^{p\bar q}$, and handle $\tilde{\rho}^{p \bar{q}}$ and $b^i$ as in (\ref{tilde-rho-bi}). We will also use that the norm of the full torsion $T(\o)=i \p \o$ is $2|T|^2$, $\na_i \| \Omega \| = \| \Omega \| T_i$, $\na_{\bar{k}} \na_i \| \Omega \| = \| \Omega \| T_i \bar{T}_{\bar{k}} + 2^{-1} \|\Omega\| R_{\bar{k} j}$, and $\| \Omega \| \leq 1$. 
\bea \label{pt-Ric^2-est1}
 \p_t | {\rm Ric}_\o |^2 
&\leq&  {1 \over 2 \| \Omega \|} \bigg\{ \Delta_F | {\rm Ric}_\o |^2  - {1 \over 2} |\na {\rm Ric}_\o |^2 - {1 \over 2} |\overline{\na} {\rm Ric}_\o |^2   \\
&&+ |\alpha' {\rm Ric}_\o | |\na {\rm Ric}_\o |^2 +  (4 + C \| \Omega \|^2)|{\rm Ric}_\o | |\na T|^2 \bigg\} \nonumber\\
&& +{C \over 2 \| \Omega \|} \bigg\{  |T| |{\rm Ric}_\o ||\na {\rm Ric}_\o |  +  \| \Omega\|^2 (1+ |T|) |{\rm Ric}_\o ||\na {\rm Ric}_\o | \nonumber\\
&& + (|{\rm Ric}_\o |+|{\rm Ric}_\o |^2) |T|^2 |\na T|+  |Rm|  |{\rm Ric}_\o |^2  + |Rm| |{\rm Ric}_\o | |T|^2 \nonumber\\
&&  + |{\rm Ric}_\o |^2 |T|^2  + |{\rm Ric}_\o | |T|^4 + |{\rm Ric}_\o |^3 (|T|+1)^2  + |{\rm Ric}_\o |^4 \nonumber\\
&&+ \| \Omega \|^2 |{\rm Ric}_\o|  (|T|+1)^4( |{\rm Ric}_\o| +|Rm| +|\na T| +1) \bigg\}.\nonumber
\eea
First, we estimate
\be
C (|{\rm Ric}_\o |+|{\rm Ric}_\o |^2) |T|^2 |\na T| \leq |{\rm Ric}_\o | |\na T|^2 + {C^2 \over 2} |{\rm Ric}_\o |(1+|{\rm Ric}_\o |)^2 |T|^4.
\ee
\be
C |T| |{\rm Ric}_\o ||\na {\rm Ric}_\o | \leq {1 \over 2} |\alpha' {\rm Ric}_\o | |\na {\rm Ric}_\o |^2 + {C^2 \over 2 |\alpha'|} |{\rm Ric}_\o | |T|^2,
\ee
We may estimate, using $|T| \leq 1$,
\be
C \| \Omega \|^2 |{\rm Ric}_\o| \, (|T|+1)^4 |\na T| \leq  \| \Omega \|^2 \, |{\rm Ric}_\o | |\na T|^2 + {C^2 \over 4 } (2)^8 \| \Omega \|^2 |{\rm Ric}_\o|,
\ee
\be
C  \| \Omega\|^2 (1+ |T|) |{\rm Ric}_\o ||\na {\rm Ric}_\o |\leq  {1 \over 2} |\alpha' {\rm Ric}_\o | |\na {\rm Ric}_\o |^2 + {1\over 2 |\alpha'|} ( 2 C \| \Omega\|^2 )^2  | {\rm Ric}_\o | .
\ee
Recall that
\be
R_{\bar{k} j}{}^\alpha{}_\beta = \hat{R}_{\bar{k} j}{}^\alpha{}_\beta + {1 \over 2} R_{\bar{k} j} \, \delta^\alpha{}_\beta.
\ee
Hence, using $\| \Omega \| \leq 1$, $|T| \leq 1$ and $|\alpha' {\rm Ric}_\o|\leq 1$ on lower order terms, from (\ref{pt-Ric^2-est1}) and the above estimates, we get
\bea
\p_t | {\rm Ric}_\o |^2 &\leq&  {1 \over 2 \| \Omega \|} \bigg\{ \Delta_F | {\rm Ric}_\o |^2  - \bigg( {1 \over 2} -  2 |\alpha' {\rm Ric}_\o| \bigg) |\na {\rm Ric}_\o |^2 + (5+ C \|\Omega\|^2) |{\rm Ric}_\o | |\na T|^2 \bigg\} \nonumber\\
&& +{C \over 2 \| \Omega \|} \bigg\{  |{\rm Ric}_\o| |T|^2 + |{\rm Ric}_\o|^2 + \| \Omega \|^2 |{\rm Ric}_\o| \bigg\}.
\eea
In terms of $0< \epsilon, \delta < 1$, we have
\bea \label{p_t-|Ric|^2-deltaepsilon}
\p_t \,  |\alpha' {\rm Ric}_\o |^2 &\leq&  {1 \over 2 \| \Omega \|} \bigg\{ \Delta_F \,  |\alpha' {\rm Ric}_\o |^2  - \bigg( {1 \over 2} -  2 \epsilon \bigg) |\alpha' \na {\rm Ric}_\o |^2 \nonumber\\
&&  + (5+C\delta^2)  \epsilon |\alpha'|^{-1} |\na T|^2 + C\delta \, \epsilon +C \epsilon^2  \bigg\}.
\eea
Using the evolution of the torsion (\ref{p_t-|T|^2})
\bea
\p_t |T|^2 &=&   {1 \over 2 \| \Omega \|} \bigg\{ \Delta_F | T |^2 - |\na T |^2_{Fg} - {1 \over 4} | {\rm Ric}_\o |^2_{Fg} -  2 \Re \{g^{i \bar{j}} g^{p \bar{q}} \na_i T_p \bar{T}_{\bar{q}} \bar{T}_{\bar{j}}\}   \nonumber\\
&& - \Re \{ g^{i \bar{j}} g^{p \bar{q}}  T_p R_{\bar{q} i} \bar{T}_{\bar{j}} \} - {1 \over 2} R |T|^2  + {\alpha' \over 4} \sigma_2(i {\rm Ric}_\o) |T|^2 \nonumber\\
&& +  \Re \{F^{p \bar{q}} g^{i \bar{j}} T^\lambda{}_{pi} R_{\bar{q} \lambda} \bar{T}_{\bar{j}} \}  -F^{p \bar{q}} g^{i \bar{j}}  T_i (\hat{R}_{\bar{q} p \bar{j}}{}^{\bar{\lambda}}+R_{\bar{q} p} \delta_j^{\bar{\lambda}}) T_{\bar{\lambda}} +  |T|^4 \nonumber\\
&&
+ {\alpha' \over 2} \| \Omega \|^3 \tilde{\rho}^{p \bar{q}} R_{\bar{q} p} |T|^2  -   |T|^2 \| \Omega \|^2 \, \nu + 2 \Re \langle E,T \rangle \bigg\}.
\eea
Estimating by replacing $F^{p \bar{q}}$ by the evolving metric $g^{p\bar q}$,
\bea
\p_t |T|^2 &\leq&   {1 \over 2 \| \Omega \|} \bigg\{ \Delta_F | T |^2 - {1 \over 2} |\na T |^2 - {1 \over 8} | {\rm Ric}_\o |^2  + 2 |\na T| |T|^2 + C|{\rm Ric}_\o| |T|^2  \nonumber\\
&& +|R| |T|^2  + {|\alpha'| \over 4} |{\rm Ric}_\o|^2  |T|^2 + \|\Omega\||\hat{Rm}|_{\hat{g}} |T|^2+  |T|^4 \nonumber\\
&& +C   \| \Omega \|^2 (|T|^4+|T|^3+ |T|^2+|T|)(1+|{\rm Ric}_\o|+|\na T|) \bigg\}.
\eea
Estimate
\be
2 |\na T| |T|^2 \leq {1 \over 8} |\na T|^2 + 8 |T|^4,
\ee
and
\be
C   \| \Omega \|^2 (|T|^4+|T|^3+ |T|^2+|T|) |\na T| \leq {1 \over 8} |\na T|^2 + 2 C^2   \| \Omega \|^4 (4)^2.
\ee
Using $0<\delta,\epsilon<1$,
\bea \label{p_t-|T|^2-deltaepsilon}
\p_t \, |T|^2 &\leq&   {1 \over 2 \| \Omega \|} \bigg\{ \Delta_F \,  | T |^2 - {1 \over 4} |\na T |^2 - {1 \over 8}  | {\rm Ric}_\o |^2 + C \epsilon \delta + C \delta  \bigg\}.
\eea
Combining (\ref{p_t-|Ric|^2-deltaepsilon}) and (\ref{p_t-|T|^2-deltaepsilon}), we obtain the desired estimate.
\begin{theorem} \label{ricci_estimate}
Start the flow with a constant function $u_0 = \log M$. There exists $M_0 \gg 1$ such that for every $M \geq M_0$, if
\be
\| \Omega \|^2 \leq  {C_2^2 \over M^2}, \ \  |T|^2 \leq {C_3 \over M},
\ee
along the flow, then
\be
 |\alpha' {\rm Ric}_\o| \leq {C_5 \over M^{1/2}},
\ee
where $C_5$ only depends on $(X,\hat{\omega})$, $\rho$, $\tilde{\mu}$ and $\alpha'$. Here, $C_2$ and $C_3$ are the constants given in Theorems \ref{C0-thm} and \ref{torsion_estimate} respectively.
\end{theorem}
{\it Proof:} Denote
\be
\epsilon = {1 \over M^{1/2}}, \ \ \delta = {C_3 \over M}.
\ee
Let $C_4$ denote the constant $C$ on the right-hand side of (\ref{lemma-Ric-estimate}), which only depends on $(X,\hat{\omega})$, $\rho$, $\tilde{\mu}$ and $\alpha'$. For $M_0$ large enough, we can simultaneously satisfy the hypothesis of Lemma \ref{Ric-delta-epsilon-lemma}, and the inequalities $2 \epsilon < {1 \over 2}$ and $(5+C_4 \delta^2)\epsilon \leq 1$. We will study the evolution equation of
\be
|\alpha' {\rm Ric}_\o |^2+ \Lambda  |T|^2,
\ee
where $\Lambda$ is a constant given by 
\be \label{Ric-Lambda}
\Lambda = \max \{ \, 4 |\alpha'|^{-1}, \, 8 |\alpha'|^2 (C_4+1) \, \}.
\ee
With this choice of $\Lambda$ and $M_0$, we have
\be
\bigg( {1 \over 2} - 2 \epsilon \bigg) \geq 0, \ \ \bigg( {\Lambda \over 4} - (5+C_4 \delta^2) \epsilon |\alpha'|^{-1} \bigg) \geq 0.
\ee
At $t=0$, $u_0= \log M$ and it follows that 
\be
\alpha'^2 | {\rm Ric}_\o |^2+ \Lambda |T|^2 =0.
\ee
Suppose that along the flow, we reach
\be
\alpha'^2 | {\rm Ric}_\o |^2+ \Lambda |T|^2 = (2 \Lambda C_3 + 1)\epsilon^2,
\ee
at some point $p \in X$ at a first time $t_0>0$. By Lemma \ref{Ric-delta-epsilon-lemma},
\bea 
\p_t ( | \alpha' {\rm Ric}_\o |^2+ \Lambda  |T|^2) &\leq&  {1 \over 2 \| \Omega \|} \bigg\{ - {\Lambda \over 8} | {\rm Ric}_\o |^2 + C_4(1+\Lambda) \epsilon \delta + C_4 \epsilon^2  + C_4 \Lambda \delta \bigg\}.
\eea
At $(p,t_0)$, we have
\be
|\alpha' {\rm Ric}_\o |^2 = (2 \Lambda C_3 + 1) \epsilon^2 - \Lambda |T|^2 \geq (2 \Lambda C_3 +1)\epsilon^2  - \Lambda \delta.
\ee
Thus
\bea 
\p_t ( |\alpha' {\rm Ric}_\o |^2+ \Lambda |T|^2) \leq  {1 \over 2 \| \Omega \|} \bigg\{ - {\Lambda \over 8 |\alpha'|^2 } \epsilon^2 + C_4 \epsilon^2 - { \Lambda^2 \over 8 |\alpha'|^2} (2 C_3 \epsilon^2 - \delta )  + C_4 \Lambda \delta + C_4 (1+\Lambda)\epsilon \delta  \bigg\}.\nonumber
\eea
After substituting the definition of $\epsilon$ and $\delta$, we obtain
\bea 
\p_t ( |\alpha' {\rm Ric}_\o |^2+ \Lambda |T|^2) &\leq&  {1 \over 2 \| \Omega \|} \bigg\{ - \bigg(  {\Lambda \over 8 |\alpha'|^2 } - C_4 \bigg) {1 \over M} - \bigg({\Lambda \over 8 |\alpha'|^2} - C_4  \bigg) {C_3 \Lambda  \over M} \nonumber\\
&& + C_3 C_4 (1+\Lambda) {1 \over M^{1/2}} {1 \over M} \bigg\}.
\eea
By our choice of $\Lambda$ (\ref{Ric-Lambda}), for $M_0 \gg 1$ depending only on $(X,\hat{\o})$, $\alpha'$, $\tilde{\mu}$, $\rho$, for all $M \geq M_0$ we have at $(p,t_0)$
\be
\p_t ( |\alpha' {\rm Ric}_\o |^2+ \Lambda |T|^2) \leq 0.
\ee
Hence along the flow, there holds
\be
| \alpha' {\rm Ric}_\o |^2+ \Lambda |T|^2 \leq (2 \Lambda C_3 + 1)\epsilon^2.
\ee
It follows that 
\be
|\alpha' {\rm Ric}_\o | \leq (2 \Lambda C_3 + 1)^{1/2} \epsilon
\ee
is preserved along the flow.

\section{Higher order estimates}
\setcounter{equation}{0}
\subsection{The evolution of derivatives of torsion}
\subsubsection{Covariant derivative of torsion}
Since $\na_{\bar{k}} T_j = {1 \over 2} R_{\bar{k} j}$, we only need to look at $\na_k T_j$. We will compute
\be
\p_t \na_i T_j = \na_i \p_t T_j - \p_t \Gamma^\lambda{}_{ij} T_\lambda.
\ee
First, using the standard formula for the evolution of the Christoffel symbols and (\ref{FY_parabolic}), we compute
\bea \label{p_t-Gamma}
\p_t \Gamma^\lambda{}_{ij} &=& g^{\lambda \bar{\mu}} \na_i \p_t g_{\bar{\mu} j} \nonumber\\
&=& \na_i  \bigg\{ {1 \over 2 \| \Omega \|} \bigg(- {R \over 2}  - {\alpha' \over 2} \| \Omega \|^3 \tilde{\rho}^{p \bar{q}} R_{\bar{q} p}+  |T|^2 + {\alpha' \over 4} \sigma_2(i {\rm Ric}_\o) + \alpha' \| \Omega \|^2 \, \nu  \bigg) \bigg\} \, \delta^\lambda{}_j \nonumber\\
&=& {1 \over 2 \| \Omega \|} \bigg\{ - {1 \over 2} \na_i R - {\alpha' \over 2} \| \Omega \|^3 \tilde{\rho}^{p \bar{q}} \na_i R_{\bar{q} p}+ {\alpha' \over 4} \sigma_2^{p \bar{q}} \na_i R_{\bar{q} p} +  g^{p \bar{q}} \na_i T_p \bar{T}_{\bar{q}} \nonumber\\
&&  +  {1 \over 2} g^{p \bar{q}}  T_p R_{\bar{q} i}  + {R \over 2}T_i -  |T|^2 T_i - {\alpha' \over 4} \sigma_2(i {\rm Ric}_\o) T_i - E_i \bigg\} \, \delta^\lambda{}_j.
\eea
We recall that the definition of $E_i$ is given in (\ref{defn_E}). Using (\ref{p_t-T})
\bea 
\p_t \na_i T_j &=& {1 \over 2 \| \Omega \|}   \na_i \bigg\{ F^{p \bar{q}} \na_p \na_{\bar{q}} T_j -  \na_j |T|^2 - {1 \over 2} T_j R + {\alpha' \over 4} T_j \sigma_2(i {\rm Ric}_\o) 
+  {1 \over 2} F^{p \bar{q}} T^\lambda{}_{pj} R_{\bar{q} \lambda}\nonumber\\
&&
 +  T_j |T|^2 + E_j \bigg\}  + \na_i \bigg\{ {1 \over 2 \| \Omega \|} \bigg\} \bigg\{  F^{p \bar{q}} \na_p \na_{\bar{q}} T_j  - g^{p \bar{q}} \na_j T_p \bar{T}_{\bar{q}} - {1 \over 2} g^{p \bar{q}}  T_p R_{\bar{q} j} 
 \nonumber\\
&&- {1 \over 2} T_j R + {\alpha' \over 4} T_j \sigma_2(i {\rm Ric}_\o) +  {1 \over 2} F^{p \bar{q}} T^\lambda{}_{pj} R_{\bar{q} \lambda} +  T_j |T|^2 + E_j \bigg\} \nonumber\\
&& -{1 \over 2 \| \Omega \|} \bigg\{ - F^{p \bar{q}} \na_p \na_{\bar{q}} T_i +  g^{p \bar{q}} \na_i T_p \bar{T}_{\bar{q}} +  {1 \over 2} g^{p \bar{q}}  T_p R_{\bar{q} i}\nonumber\\
&&  + {R \over 2}T_i -  |T|^2 T_i - {\alpha' \over 4} \sigma_2(i {\rm Ric}_\o) T_i - E_i  \bigg\} T_j.
\eea
First, we may rewrite
\be \label{p_t-T-id1}
F^{p \bar{q}} \na_p \na_{\bar{q}} T_j = {1 \over 2} F^{p \bar{q}} \na_p  R_{\bar{q} j} .
\ee
Next,
\bea \label{p_t-T-id2}
\na_i \{ F^{p \bar{q}} \na_p \na_{\bar{q}} T_j \}
&=&   F^{p \bar{q}} \na_i \na_p \na_{\bar{q}} T_j + \na_i \bigg( \alpha' \| \Omega \|^3 \tilde{\rho}^{p \bar{q}} - {\alpha' \over 2} \sigma_2^{p \bar{q}} \bigg)\na_p \na_{\bar{q}} T_j \nonumber\\
&=&  F^{p \bar{q}} \na_p \na_i \na_{\bar{q}} T_j +  F^{p \bar{q}} T^\lambda{}_{pi} \na_\lambda \na_{\bar{q}} T_j + \alpha' \na_i ( \| \Omega \|^3 \tilde{\rho}^{p \bar{q}} ) \na_p \na_{\bar{q}} T_j \nonumber\\
&&- {\alpha' \over 2}  \sigma_2^{p \bar{q}, r \bar{s}} \na_i R_{\bar{s} r} \na_p \na_{\bar{q}} T_j \nonumber\\
&=&  F^{p \bar{q}} \na_p  \na_{\bar{q}}\na_i T_j -   F^{p \bar{q}} \na_p ( R_{\bar{q} i}{}^\lambda{}_j T_\lambda) +  F^{p \bar{q}} T^\lambda{}_{pi} \na_\lambda R_{\bar{q} j} \nonumber\\
&& + {\alpha' \over 2} \na_i ( \| \Omega \|^3 \tilde{\rho}^{p \bar{q}} ) \na_p R_{\bar{q}j} - {\alpha' \over 4}  \sigma_2^{p \bar{q}, r \bar{s}} \na_i R_{\bar{s} r} \na_p R_{\bar{q} j}.
\eea
We also compute
\bea \label{p_t-T-id3}
\na_i \na_j|T|^2 &=& g^{p \bar{q}} \na_i \na_j T_p \bar{T}_{\bar{q}} + g^{p \bar{q}} \na_j T_p \na_i \bar{T}_{\bar{q}} + g^{p \bar{q}} \na_i T_p \na_j \bar{T}_{\bar{q}} + g^{p \bar{q}}  T_p \na_i \na_j \bar{T}_{\bar{q}}  \nonumber\\
&=&   g^{p \bar{q}} \na_i \na_j T_p \bar{T}_{\bar{q}} + {1 \over 2} g^{p \bar{q}} \na_j T_p R_{\bar{q} i} + {1 \over 2} g^{p \bar{q}} \na_i T_p R_{\bar{q} j} + {1 \over 2} g^{p \bar{q}}  T_p \na_i R_{\bar{q} j}.
\eea
We introduce the notation ${\cal E}$, which denotes any combination of terms involving only $Rm$, $T$, $g$, $\| \Omega \|$, $\alpha'$, $\rho$ and $\mu$, as well as any derivatives of $\rho$ and $\mu$. Note that $F^{p \bar{q}}$ is an element of ${\cal E}$. The notation $\ast$ refers to a contraction using the evolving metric $g$. The notation $D {\cal E}$ denotes any term which is a covariant derivative of a term in ${\cal E}$. For example, the group $D {\cal E}$ contains terms involving $\na T$, $\bar{\na} \bar{T}$, and $\na {\rm Ric}_\o$. Substituting (\ref{p_t-T-id1}), (\ref{p_t-T-id2}), (\ref{p_t-T-id3}) gives
\bea \label{p_t-naT}
\p_t \na_i T_j &=& {1 \over 2 \| \Omega \|}  \bigg\{ \Delta_F \na_i T_j - {\alpha' \over 4}  \sigma_2^{p \bar{q}, r \bar{s}} \na_i R_{\bar{s} r} \na_p R_{\bar{q} j} + \na \na T \ast {\cal E} + D {\cal E} \ast {\cal E} + {\cal E}  \bigg\}.
\eea
Here we also used that $\na_i E_j = \na \na T \ast {\cal E} + D {\cal E} \ast {\cal E} + {\cal E}$ which can be verified from the definition of $E_j$ given in (\ref{defn_E})
\subsubsection{Norm of covariant derivative of torsion}
We will compute
\be
\p_t |\na T|^2 = \p_t \{ g^{i \bar{j}} g^{k \bar{\ell}} \na_i T_k \na_{\bar{j}} \bar{T}_{\bar{\ell}} \}.
\ee
As in (\ref{p_t-normT1}), we have
\bea \label{p_t-|naT|^2-1}
\p_t | \na T|^2 &=&  2 \Re \langle \p_t \na T, \na T \rangle \nonumber\\
&&+  2 {| \na T |^2 \over 2 \| \Omega \|} \bigg( {R \over 2} +{\alpha' \over 2} \| \Omega \|^3 \tilde{\rho}^{p \bar{q}} R_{\bar{q} p}  - {\alpha' \over 4} \sigma_2(i {\rm Ric}_\o) -  |T|^2 -  \| \Omega \|^2 \, \nu  \bigg).
\eea
Next,
\bea
\Delta_F |\na T|^2 &=&  F^{p \bar{q}} g^{i \bar{j}} g^{k \bar{\ell}} \na_p \na_{\bar{q}} \na_i T_k \na_{\bar{j}} \bar{T}_{\bar{\ell}} + g^{i \bar{j}} g^{k \bar{\ell}}  \na_i T_k \overline{F^{q \bar{p}} \na_{\bar{p}} \na_q \na_j T_\ell} \nonumber\\
&&+ F^{p \bar{q}} g^{i \bar{j}} g^{k \bar{\ell}} \na_p  \na_i T_k \na_{\bar{q}} \na_{\bar{j}} \bar{T}_{\bar{\ell}} + F^{p \bar{q}} g^{i \bar{j}} g^{k \bar{\ell}} \na_{\bar{q}}  \na_i T_k  \na_p\na_{\bar{j}} \bar{T}_{\bar{\ell}} \nonumber\\
&=& 2 \Re \langle \Delta_F \na T, \na T \rangle +g^{i \bar{j}} g^{k \bar{\ell}}  \na_i T_k F^{p \bar{q}} R_{\bar{q} p \bar{j}}{}^{\bar{\lambda}} \na_{\bar{\lambda}} \bar{T}_{\bar{\ell}} \nonumber\\
&&+ g^{i \bar{j}} g^{k \bar{\ell}}  \na_i T_k F^{p \bar{q}} R_{\bar{q} p \bar{\ell}}{}^{\bar{\lambda}} \na_{\bar{j}} T_{\bar{\lambda}}+ |\na \na T|^2_{Fgg} +F^{p \bar{q}} g^{i \bar{j}} g^{k \bar{\ell}} \na_{\bar{q}}  \na_i T_k  \na_p\na_{\bar{j}} \bar{T}_{\bar{\ell}}.\nonumber
\eea
The last term can be written as a norm of $\na {\rm Ric}_\o$ plus commutator terms. Explicitly,
\bea
F^{p \bar{q}} g^{i \bar{j}} g^{k \bar{\ell}} \na_{\bar{q}}  \na_i T_k  \na_p\na_{\bar{j}} \bar{T}_{\bar{\ell}} 
&=& F^{p \bar{q}} g^{i \bar{j}} g^{k \bar{\ell}}  \na_i  \na_{\bar{q}} T_k  \overline{ \na_{\bar{p}} \na_{j} T_\ell} + F^{p \bar{q}} g^{i \bar{j}} g^{k \bar{\ell}}  R_{\bar{q} i}{}^\lambda{}_k T_\lambda  \na_p\na_{\bar{j}} \bar{T}_{\bar{\ell}} \nonumber\\
&=& F^{p \bar{q}} g^{i \bar{j}} g^{k \bar{\ell}}  \na_i  \na_{\bar{q}} T_k  \overline{\na_{j} \na_{\bar{p}}  T_\ell}+ F^{p \bar{q}} g^{i \bar{j}} g^{k \bar{\ell}}  \na_i  \na_{\bar{q}} T_k R_{\bar{j} p \bar{\ell}}{}^{\bar{\lambda}} T_{\bar{\lambda}} 
\nonumber\\
&&+ F^{p \bar{q}} g^{i \bar{j}} g^{k \bar{\ell}}  R_{\bar{q} i}{}^\lambda{}_k T_\lambda  \na_{\bar{j}} \na_p \bar{T}_{\bar{\ell}}+ F^{p \bar{q}} g^{i \bar{j}} g^{k \bar{\ell}}  R_{\bar{q} i}{}^\lambda{}_k T_\lambda  R_{\bar{j} p \bar{\ell}}{}^{\bar{\lambda}} \bar{T}_{\bar{\lambda}} \nonumber\\
&=& {1 \over 4} F^{p \bar{q}} g^{i \bar{j}} g^{k \bar{\ell}}  \na_i R_{\bar{q} k}  \overline{\na_{j} R_{\bar{p} \ell}}+ {1 \over 2}F^{p \bar{q}} g^{i \bar{j}} g^{k \bar{\ell}}  \na_i R_{\bar{q} k} R_{\bar{j} p \bar{\ell}}{}^{\bar{\lambda}} T_{\bar{\lambda}} \nonumber\\
&&+ {1 \over 2} F^{p \bar{q}} g^{i \bar{j}} g^{k \bar{\ell}}  R_{\bar{q} i}{}^\lambda{}_k T_\lambda  \na_{\bar{j}} R_{\bar{\ell} p}+ F^{p \bar{q}} g^{i \bar{j}} g^{k \bar{\ell}}  R_{\bar{q} i}{}^\lambda{}_k T_\lambda  R_{\bar{j} p \bar{\ell}}{}^{\bar{\lambda}} \bar{T}_{\bar{\lambda}}. \nonumber\\
\eea
Hence
\bea \label{p_t-|naT|^2-2}
\Delta_F |\na T|^2 
&=& 2 \Re \langle \Delta_F \na T, \na T \rangle+ |\na \na T|^2_{Fgg} + {1 \over 4} |\na {\rm Ric}_\o|^2_{Fgg}  \nonumber\\
&&+ g^{i \bar{j}} g^{k \bar{\ell}}  \na_i T_k F^{p \bar{q}} R_{\bar{q} p \bar{j}}{}^{\bar{\lambda}} \na_{\bar{\lambda}} \bar{T}_{\bar{\ell}} + g^{i \bar{j}} g^{k \bar{\ell}}  \na_i T_k F^{p \bar{q}} R_{\bar{q} p \bar{\ell}}{}^{\bar{\lambda}} \na_{\bar{j}} T_{\bar{\lambda}}   \nonumber\\
&&+ {1 \over 2}F^{p \bar{q}} g^{i \bar{j}} g^{k \bar{\ell}}  \na_i R_{\bar{q} k} R_{\bar{j} p \bar{\ell}}{}^{\bar{\lambda}} T_{\bar{\lambda}} + {1 \over 2} F^{p \bar{q}} g^{i \bar{j}} g^{k \bar{\ell}}  R_{\bar{q} i}{}^\lambda{}_k T_\lambda  \na_{\bar{j}} R_{\bar{\ell} p} \nonumber\\
&& + F^{p \bar{q}} g^{i \bar{j}} g^{k \bar{\ell}}  R_{\bar{q} i}{}^\lambda{}_k T_\lambda  R_{\bar{j} p \bar{\ell}}{}^{\bar{\lambda}} \bar{T}_{\bar{\lambda}}.
\eea
Therefore, by (\ref{p_t-naT}), (\ref{p_t-|naT|^2-1}) and (\ref{p_t-|naT|^2-2}), 
\bea \label{p_t-|naT|^2-3}
\p_t |\na T|^2 &=& {1 \over 2 \| \Omega \|}  \bigg\{ \Delta_F |\na T|^2 -|\na \na T|^2_{Fgg} - {1 \over 4} |\na  {\rm Ric}_\o|^2_{Fgg}   \nonumber\\
&&- {\alpha' \over 2} \Re \{ g^{i \bar{j}} g^{k \bar{\ell}} \sigma_2^{p \bar{q}, r \bar{s}} \na_i R_{\bar{s} r} \na_p R_{\bar{q} k}\na_{\bar{j}} \bar{T}_{\bar{\ell}} \} + \na \na T \ast \na T \ast {\cal E} \nonumber\\
&& + D {\cal E} \ast D {\cal E} \ast {\cal E}+ D {\cal E} \ast {\cal E} +{\cal E} \bigg\}.
\eea

\subsection{The evolution of derivatives of curvature}
\subsubsection{Derivative of Ricci curvature}
We will compute
\be
\p_t \na_i R_{\bar{k} j} = \na_i \p_t R_{\bar{k} j} -\p_t \Gamma^\lambda{}_{ij} R_{\bar{k} \lambda}.
\ee
Using (\ref{p_t-Ric2}) and (\ref{p_t-Gamma}), we obtain
\bea 
\p_t \na_i R_{\bar{k} j} &=&  {1 \over 2 \| \Omega \|} \bigg\{ \na_i ( F^{p \bar{q}}  \na_p \na_{\bar{q}} R_{\bar{k} j} )  -{\alpha' \over 2} \na_i ( \sigma_2^{p \bar{q}, r \bar{s}} \na_{\bar{k}} R_{\bar{s} r} \na_j R_{\bar{q} p})\nonumber\\
&& + (2g^{p \bar{q}} +2\alpha' \| \Omega \|^3 \tilde{\rho}^{p \bar{q}}) \ast \na \na T \ast \na T   +DD {\cal E} \ast {\cal E}  \nonumber\\
&& + D {\cal E} \ast D {\cal E} \ast {\cal E}  +D {\cal E} \ast {\cal E} + {\cal E}\bigg\}. 
\eea
Here, we used that $\na \bar{\na} \bar{T} = \bar{\na} {\rm Ric}_\o + Rm \ast \bar{T}$. Compute
\bea
\na_i ( F^{p \bar{q}}  \na_p \na_{\bar{q}} R_{\bar{k} j} ) &=&  F^{p \bar{q}} \na_i \na_p \na_{\bar{q}} R_{\bar{k} j} +\alpha' \na_i (\|\Omega\|^3 \tilde{\rho}^{p \bar{q}})\na_p \na_{\bar{q}} R_{\bar{k} j} -{\alpha' \over 2} \na_i (\sigma_2^{p \bar{q}})  \na_p \na_{\bar{q}} R_{\bar{k} j} \nonumber\\
&=& F^{p \bar{q}}  \na_p \na_i \na_{\bar{q}} R_{\bar{k} j}+  F^{p \bar{q}}  T^\lambda{}_{pi} \na_\lambda \na_{\bar{q}} R_{\bar{k} j} \nonumber\\
&&+\alpha' \na_i (\|\Omega\|^3 \tilde{\rho}^{p \bar{q}})\na_p \na_{\bar{q}} R_{\bar{k} j}-{\alpha' \over 2} \sigma_2^{p \bar{q}, r \bar{s}} \na_i R_{\bar{s} r}  \na_p \na_{\bar{q}} R_{\bar{k} j} \nonumber\\
&=& F^{p \bar{q}}  \na_p  \na_{\bar{q}} \na_i R_{\bar{k} j} + F^{p \bar{q}}  \na_p (R_{\bar{q} i \bar{k}}{}^{\bar{\lambda}} R_{\bar{\lambda} j}  -R_{\bar{q} i}{}^\lambda{}_j R_{\bar{k} \lambda} )  \nonumber\\
&&+  F^{p \bar{q}}  T^\lambda{}_{pi} \na_\lambda \na_{\bar{q}} R_{\bar{k} j}+\alpha' \na_i (\|\Omega\|^3 \tilde{\rho}^{p \bar{q}})\na_p \na_{\bar{q}} R_{\bar{k} j} \nonumber\\
&&-{\alpha' \over 2} \sigma_2^{p \bar{q}, r \bar{s}} \na_i R_{\bar{s} r}  \na_p \na_{\bar{q}} R_{\bar{k} j}. 
\eea
Hence, using that $\na_i  \sigma_2^{p \bar{q}, r \bar{s}} =0$ (\ref{sigma_2-pqrs}), we obtain
\bea 
\p_t \na_i R_{\bar{k} j} &=&  {1 \over 2 \| \Omega \|} \bigg\{ \Delta_F \na_i R_{\bar{k} j}   -{\alpha' \over 2}   \sigma_2^{p \bar{q}, r \bar{s}} \na_i \na_{\bar{k}} R_{\bar{s} r} \na_j R_{\bar{q} p}   \nonumber\\
&&-{\alpha' \over 2}   \sigma_2^{p \bar{q}, r \bar{s}}  \na_{\bar{k}} R_{\bar{s} r} \na_i \na_j R_{\bar{q} p} -{\alpha' \over 2} \sigma_2^{p \bar{q}, r \bar{s}} \na_i R_{\bar{s} r} \na_p \na_{\bar{q}} R_{\bar{k} j} \nonumber\\
&&  + (2g^{p \bar{q}} +2\alpha' \| \Omega \|^3 \tilde{\rho}^{p \bar{q}}) \ast \na \na T \ast \na T  \nonumber\\
&& +DD {\cal E} \ast {\cal E} + D {\cal E} \ast D {\cal E} \ast {\cal E} +D {\cal E} \ast {\cal E} + {\cal E}\bigg\}.
\eea

\subsubsection{Norm of derivative of Ricci curvature}
We will compute
\be
\p_t |\na {\rm Ric}_\o |^2 = \p_t \{ g^{i \bar{a}} g^{b \bar{k}} g^{j \bar{c}} \na_i R_{\bar{k} j} \overline{\na_a R_{\bar{b} c}} \}.
\ee
As in (\ref{p_t-normT1}), we have
\bea 
\p_t |\na {\rm Ric}_\o|^2 &=&  2 \Re \langle \p_t \na {\rm Ric}_\o, \na {\rm Ric}_\o \rangle \nonumber\\
&&+  3 | \na {\rm Ric}_\o |^2 {1 \over 2 \| \Omega \|} \bigg( {R \over 2} + {\alpha' \over 2} \| \Omega \|^3 \tilde{\rho}^{p \bar{q}} R_{\bar{q} p}  - {\alpha' \over 4} \sigma_2(i {\rm Ric}_\o) -  |T|^2 -  \| \Omega \|^2 \, \nu  \bigg). \nonumber
\eea
Next, compute
\bea
 \Delta_F |\na {\rm Ric}_\o |^2 &=& F^{p \bar{q}}  g^{i \bar{j}} g^{k \bar{\ell}} g^{m \bar{n}}\na_p \na_{\bar{q}} \na_i  R_{\bar{n} k} \overline{\na_j R_{\bar{m} \ell}} +   g^{i \bar{j}} g^{k \bar{\ell}} g^{m \bar{n}} \na_i  R_{\bar{n} k} \overline{ F^{q \bar{p}} \na_{\bar{p}} \na_q \na_j R_{\bar{m} \ell}}  \nonumber\\
&&+ |\na \na {\rm Ric}_\o|^2_{Fggg} +|\overline{\na} \na {\rm Ric}_\o|^2_{Fggg} \nonumber\\
&=& 2 \Re \langle \Delta_F \na {\rm Ric}_\o, \na  {\rm Ric}_\o \rangle + |\na \na {\rm Ric}_\o|^2_{Fggg} +|\overline{\na} \na {\rm Ric}_\o|^2_{Fggg}\nonumber\\
&&+  F^{p \bar{q}} g^{i \bar{j}} g^{k \bar{\ell}} g^{m \bar{n}} \na_i  R_{\bar{n} k} R_{\bar{q} p \bar{j}}{}^{\bar{\lambda}} \na_{\bar{\lambda}} R_{\bar{\ell} m}+  F^{p \bar{q}} g^{i \bar{j}} g^{k \bar{\ell}} g^{m \bar{n}} \na_i  R_{\bar{n} k} R_{\bar{q} p \bar{\ell}}{}^{\bar{\lambda}} \na_{\bar{j}} R_{\bar{\lambda} m} \nonumber\\
&&-  F^{p \bar{q}} g^{i \bar{j}} g^{k \bar{\ell}} g^{m \bar{n}} \na_i  R_{\bar{n} k} R_{\bar{q} p}{}^\lambda{}_m \na_{\bar{j}} R_{\bar{\ell} \lambda}.
\eea
Commuting covariant derivatives
\be
|\overline{\na} \na {\rm Ric}_\o|^2_{Fggg} = |\na \overline{\na} {\rm Ric}_\o|^2_{Fggg} +\na \overline{\na} {\cal E} \ast {\cal E}+{\cal E}.
\ee
Hence
\bea \label{p_t-|naRic|^2-0}
 \p_t |\na {\rm Ric}_\o|^2
&=&  {1 \over 2 \| \Omega \|} \bigg\{ \Delta_F |\na {\rm Ric}_\o|^2 -|\na \na {\rm Ric}_\o|^2_{Fggg} -|\na \overline{\na} {\rm Ric}_\o|^2_{Fggg}\bigg\} \nonumber\\
&&+{1 \over 2 \| \Omega \|} 2 \Re \bigg\{    -{\alpha' \over 2} g^{i \bar{a}} g^{b \bar{k}} g^{j \bar{c}} \sigma_2^{p \bar{q}, r \bar{s}} \na_i \na_{\bar{k}} R_{\bar{s} r} \na_j R_{\bar{q} p}\overline{\na_a R_{\bar{b} c}} \nonumber\\
&&-{\alpha' \over 2}g^{i \bar{a}} g^{b \bar{k}} g^{j \bar{c}} \sigma_2^{p \bar{q}, r \bar{s}}  \na_{\bar{k}} R_{\bar{s} r} \na_i \na_j R_{\bar{q} p} \overline{\na_a R_{\bar{b} c}} \nonumber\\
&&-{\alpha' \over 2}g^{i \bar{a}} g^{b \bar{k}} g^{j \bar{c}} \sigma_2^{p \bar{q}, r \bar{s}} \na_i R_{\bar{s} r} \na_p \na_{\bar{q}} R_{\bar{k} j} \overline{\na_a R_{\bar{b} c}}\nonumber\\
&&+ (2g^{p \bar{q}} +2\alpha' \| \Omega \|^3 \tilde{\rho}^{p \bar{q}}) \ast \na \na T \ast \na T \ast \na {\rm Ric}_\o \bigg\} \nonumber\\
&&  +DD {\cal E} \ast D {\cal E}\ast {\cal E}+ DD {\cal E} \ast {\cal E} + D {\cal E} \ast D {\cal E}\ast D {\cal E} \ast {\cal E} \nonumber\\
&&+ D {\cal E}\ast D {\cal E} \ast {\cal E}+ D {\cal E} \ast {\cal E} . 
\eea
\begin{lemma}
Suppose $|\alpha' {\rm Ric}_\o| \leq {1 \over 4}$ and $-{1 \over 8} g^{p \bar{q}} < \alpha' \| \Omega \|^3 \tilde{\rho}^{p \bar{q}} < {1 \over 8} g^{p \bar{q}}$. Then
\bea \label{p_t-|naRic|^2}
\p_t |\na {\rm Ric}_\o|^2 &\leq&  {1 \over 2 \| \Omega \|} \bigg\{ \Delta_F |\na {\rm Ric}_\o|^2 -{1 \over 2} |\na \na {\rm Ric}_\o|^2 - {1 \over 2} |\na \overline{\na} {\rm Ric}_\o|^2 \bigg\} \nonumber\\
&&+{1 \over 2 \| \Omega \|}  \bigg\{ 9 \alpha'^2 |\na {\rm Ric}_\o|^4 + 5 |\na \na T|  |\na T| |\na {\rm Ric}_\o| \nonumber\\
&&  +DD {\cal E} \ast D {\cal E}\ast {\cal E} + (D {\cal E} + {\cal E})^3 \bigg\}. 
\eea
\end{lemma}
{\it Proof:} By assumption, we may use
\be
|\na \na {\rm Ric}_\o|^2_{Fggg} +|\na \overline{\na} {\rm Ric}_\o|^2_{Fggg} \geq {3 \over 4} (|\na \na {\rm Ric}_\o|^2 +|\na \overline{\na} {\rm Ric}_\o|^2).
\ee
In coordinates where the evolving metric $g$ is the identity, we have $\sigma_2^{p \bar{q}, r \bar{s}}= \pm 1$. Using $2ab \leq a^2 + b^2$, estimate (\ref{p_t-|naRic|^2}) follows from (\ref{p_t-|naRic|^2-0}).
\subsection{Higher order estimates}

\begin{theorem} \label{higher_order_est}
There exists $0<\delta_1, \delta_2$ with the following property. Suppose 
\be
-{1 \over 8} g^{p \bar{q}} < \alpha' \| \Omega \|^3 \tilde{\rho}^{p \bar{q}} < {1 \over 8} g^{p \bar{q}}, \ \ \| \Omega \| \leq 1,
\ee
\be \label{Ric-small}
|\alpha' {\rm Ric}_\o| \leq \delta_1,
\ee
and 
\be 
|T|^2 \leq \delta_2,
\ee
along the flow. Then
\be
|\na {\rm Ric}_\o| \leq C, \ \ |\na T| \leq C,
\ee
where $C$ depends only on $\delta_1$, $\delta_2$, $\alpha'$, $\rho$, $\tilde{\mu}$, and $(X,\hat{\omega})$.
\end{theorem}
\par \noindent {\it Proof:} Let us assume that $\delta_1 < {1 \over 4}$. This will allow us to use the estimate
\be \label{F-2-g}
{3 \over 4} g_{\bar{k} j} \leq F^{j \bar{k}} \leq 2 g_{\bar{k} j}.
\ee
This follows from the definition of $F^{j \bar{k}}$, see (\ref{F^ii-in-coords}). From (\ref{p_t-Ric-full}), with assumptions (\ref{Ric-small}) and (\ref{F-2-g}) we may estimate
\bea \label{p_t-|Ric|^2-again}
\p_t |{\rm Ric}_\o|^2 &\leq& {1 \over 2 \| \Omega \|} \bigg\{ \Delta_F | {\rm Ric}_\o |^2  - {1 \over 2} |\na {\rm Ric}_\o |^2  \bigg\} \nonumber\\
&& + {1 \over 2 \| \Omega \|}  \Re \bigg\{ D {\cal E} \ast {\cal E} + 5 \, \na T \ast \na T \ast {\rm Ric} + {\cal E} \bigg\}.
\eea
Here we used
\be
 -\alpha' \Re \{ g^{j \bar{\ell}} g^{m \bar{k}} \sigma_2^{p \bar{q}, r \bar{s}}R_{\bar{\ell} m} \na_{\bar{k}} R_{\bar{s} r} \na_j R_{\bar{q} p} \} \leq \delta_1 | \na {\rm Ric}_\o|^2,
\ee
to absorb this term into the $-|\na {\rm Ric}_\o |^2$ term. We will compute the evolution of
\be
G = (|\alpha' {\rm Ric}_\o|^2+ \tau_1) |\na {\rm Ric}_\o|^2+ (|T|^2+ \tau_2) |\na T|^2,
\ee
where $\tau_1$ and $\tau_2$ are constants to be determined. First, we compute
\be
\p_t \{ (|\alpha' {\rm Ric}_\o|^2+ \tau_1 ) |\na {\rm Ric}_\o|^2 \} = \alpha'^2 \p_t|{\rm Ric}_\o|^2 |\na {\rm Ric}_\o|^2 + (|\alpha' {\rm Ric}_\o|^2+ \tau_1)\p_t|\na {\rm Ric}_\o|^2.
\ee
By (\ref{p_t-|naRic|^2}) and (\ref{p_t-|Ric|^2-again})
\bea  
&&\p_t \{ (|\alpha' {\rm Ric}_\o|^2+ \tau_1) |\na {\rm Ric}_\o|^2 \} \nonumber\\
 &\leq&  {1 \over 2 \| \Omega \|} \bigg\{ \Delta_F |\alpha' {\rm Ric}_\o |^2 |\na {\rm Ric}_\o|^2 - {\alpha'^2 \over 2} |\na {\rm Ric}_\o |^4  \bigg\} \nonumber\\
&& + {1 \over 2 \| \Omega \|}  \Re \bigg\{ D {\cal E} \ast {\cal E} + 5 \, \na T \ast \na T \ast {\rm Ric} + {\cal E} \bigg\}\alpha'^2 |\na {\rm Ric}_\o|^2 \nonumber\\
&& + {(|\alpha' {\rm Ric}_\o|^2+ \tau_1) \over 2 \| \Omega \|} \bigg\{ \Delta_F |\na {\rm Ric}_\o|^2 -{1 \over 2} |\na \na {\rm Ric}_\o|^2 - {1 \over 2} |\na \overline{\na} {\rm Ric}_\o|^2 \bigg\} \nonumber\\
&&+{(|\alpha'{\rm Ric}_\o|^2+ \tau_1) \over 2 \| \Omega \|}  \bigg\{ 9 \alpha'^2 |\na {\rm Ric}_\o|^4 + 5 |\na \na T|  |\na T| |\na {\rm Ric}_\o| \nonumber\\
&&  +\na \na {\cal E} \ast D {\cal E}\ast {\cal E}+ \na \overline{\na} {\cal E} \ast D {\cal E} \ast {\cal E} + (D {\cal E} + {\cal E})^3 \bigg\} .
\eea
Hence
\bea & \ &
\p_t \{ (|\alpha' {\rm Ric}_\o|^2+ \tau_1) |\na {\rm Ric}_\o|^2 \} \nonumber\\
 &\leq&  {1 \over 2 \| \Omega \|} \bigg\{ \Delta_F \{ (|\alpha' {\rm Ric}_\o|^2+ \tau_1) |\na {\rm Ric}_\o|^2 \} - \bigg( {1 \over 2}   - 9 |\alpha' {\rm Ric}_\o|^2 - 9 \tau_1 \bigg) \alpha'^2 |\na {\rm Ric}_\o|^4 \nonumber\\
&& -{1 \over 2} |\na \na {\rm Ric}_\o|^2(|\alpha' {\rm Ric}_\o|^2+ \tau_1) - {1 \over 2} |\na \overline{\na} {\rm Ric}_\o|^2(|\alpha' {\rm Ric}_\o|^2+ \tau_1)  \nonumber\\
&&- 2 \Re \, \{ F^{i \bar{j}} \na_i|\alpha' {\rm Ric}_\o|^2 \na_{\bar{j}} |\na {\rm Ric}_\o|^2 \} + 6 (\delta_1^2+\tau) |\na \na T|  |\na T| |\na {\rm Ric}_\o|  \bigg\} \nonumber\\
&& + {\alpha'^2 |\na {\rm Ric}_\o|^2 \over 2 \| \Omega \|}  \Re \bigg\{ 5 \, \na T \ast \na T \ast {\rm Ric} + D {\cal E} \ast {\cal E} +  {\cal E} \bigg\} \nonumber\\
&&+{(|\alpha' {\rm Ric}_\o|^2+ \tau_1) \over 2 \| \Omega \|}  \bigg\{  \na \na {\cal E} \ast D {\cal E}\ast {\cal E}+ \na \overline{\na} {\cal E} \ast D {\cal E} \ast {\cal E} + (D {\cal E} + {\cal E})^3 \bigg\} .
\eea
We estimate
\bea
&&- 2 \Re \, \{ F^{i \bar{j}} \na_i|\alpha' {\rm Ric}_\o|^2 \na_{\bar{j}} |\na {\rm Ric}_\o|^2 \} \nonumber\\
&\leq& 8 |\alpha'| \delta_1  |\na {\rm Ric}_\o|^2 (|\na \na {\rm Ric}_\o| + |\na \overline{\na} {\rm Ric}_\o| + {\cal E})\nonumber\\
&\leq& {\alpha'^2  \over 2^4} |\na {\rm Ric}_\o|^4 + 2^8 \delta_1^2( |\na \na {\rm Ric}_\o|^2+ |\na \overline{\na} {\rm Ric}_\o|^2) + C |\na {\rm Ric}_\o|^2,
\eea
\bea
&& 6 (\delta_1^2 +\tau_1) |\na \na T|  |\na T| |\na {\rm Ric}_\o| 
\nonumber\\
&\leq&  {1 \over 2} (\delta_1^2+\tau_1)  |\na \na T|^2 + 2^1 3^2 (\delta_1^2+\tau_1)  |\na T|^2 |\na {\rm Ric}_\o|^2 \nonumber\\
&\leq& {1 \over 2} (\delta_1^2+\tau_1) |\na \na T|^2  + {\alpha'^2 \over 2^4} |\na {\rm Ric}_\o|^4+ 2^4 3^4 \alpha'^{-2} (\delta_1^2+\tau_1)^2 |\na T|^4,
\eea
\bea
& \ & {\alpha'^2 |\na {\rm Ric}_\o|^2 \over 2 \| \Omega \|}  \Re \bigg\{ 5 \, \na T \ast \na T \ast {\rm Ric} + \na {\cal E} \ast {\cal E} +  {\cal E} \bigg\} \nonumber\\
&\leq& {1 \over 2 \| \Omega \|} \bigg\{ {\alpha'^2 \over 2^4} |\na {\rm Ric}_\o|^4 + 2^2 5^2 \delta_1^2  |\na T|^4+ C|\na {\rm Ric}_\o|^3 + C|\na T|^3 + C \bigg\}.
\eea
\bea & \ &
{(|\alpha' {\rm Ric}_\o|^2+ \tau_1) \over 2 \| \Omega \|}  \bigg\{  \na \na {\cal E} \ast D {\cal E}\ast {\cal E}+ \na \overline{\na} {\cal E} \ast D {\cal E} \ast {\cal E} + (D {\cal E} + {\cal E})^3 \bigg\} \nonumber\\
&\leq& {1 \over 2 \| \Omega \|} \bigg\{ {1 \over 4} |\na \na {\rm Ric}_\o|^2(|\alpha' {\rm Ric}_\o|^2+ \tau_1) + {1 \over 4} |\na \overline{\na} {\rm Ric}_\o|^2(|\alpha' {\rm Ric}_\o|^2+ \tau_1) \nonumber\\
&& + {1 \over 2} (\delta_1^2+\tau_1) |\na \na T|^2  + C|\na {\rm Ric}_\o|^3 + C|\na T|^3 + C \bigg\}.
\eea
Therefore 
\bea & \ & \label{G-first-piece}
\p_t \{ (|\alpha' {\rm Ric}_\o|^2+ \tau_1) |\na {\rm Ric}_\o|^2 \} \\
 &\leq&  {1 \over 2 \| \Omega \|} \bigg\{ \Delta_F \{ (|\alpha' {\rm Ric}_\o|^2+ \tau_1) |\na {\rm Ric}_\o|^2 \} - \bigg( {1 \over 4}   - 9 \delta_1^2 - 9 \tau_1 \bigg) \alpha'^2 |\na {\rm Ric}_\o|^4 \nonumber\\
&& - (|\na \na {\rm Ric}_\o|^2 + |\na \overline{\na} {\rm Ric}_\o|^2)({\tau_1 \over 4} - 2^8 \delta_1^2)  + (\delta_1^2 +\tau_1)  |\na \na T|^2  \nonumber\\
&& + \bigg( 2^4 3^4 \alpha'^{-2} (\delta_1^2+\tau_1)^2 +2^2 5^2 \delta_1^2\bigg)  |\na T|^4 + C_{\alpha',\tau,\delta} |\na {\rm Ric}_\o|^3 + C_{\alpha',\tau,\delta} |\na T|^3+  C_{\alpha',\tau,\delta}  \bigg\}.\nonumber
\eea
Next, we compute
\be
\p_t \{ (|T|^2+ \tau_2) |\na T|^2 \} = \p_t|T|^2 |\na T|^2 + (|T|^2+ \tau_2 ) \p_t|\na T|^2.
\ee
By (\ref{p_t-|T|^2}), we have
\be
\p_t |T|^2 \leq {1 \over 2 \| \Omega \|} \bigg\{ \Delta_F | T |^2 - |\na T |^2_{Fg} + C |\na T| + C \bigg\}.
\ee
By (\ref{p_t-|naT|^2-3}), we have
\bea
\p_t |\na T|^2 &\leq& {1 \over 2 \| \Omega \|}  \bigg\{ \Delta_F |\na T|^2 -|\na \na T|^2_{Fgg} +|\alpha'| |\na T| |\na {\rm Ric}_\o|^2 \nonumber\\
&& +C |\na \na T| |\na T| +C |\na T|^2 + C |\na {\rm Ric}_\o|^2 +C  \bigg\}.
\eea
By our assumption $| \alpha' {\rm Ric}_\o| \leq {1 \over 4}$, we have $|\na \na T|^2_{Fgg} \geq {1 \over 2} |\na \na T|^2$ and $|\na T|^2_{Fg} \geq {1 \over 2} |\na T|^2$. Therefore
\bea & \ &
\p_t \{ (|T|^2+ \tau_2 ) |\na T|^2 \} \nonumber\\
&\leq& {1 \over 2 \| \Omega \|} \bigg\{ \Delta_F \{ (|T|^2+ \tau_2 )  |\na T|^2 \} - 2  \Re \{ F^{i \bar{j}} \na_i |T|^2 \na_{\bar{j}} |\na T|^2 \} \nonumber\\
&& - {1 \over 4} |\na T |^4 -(|T|^2+ \tau_2 ) {1 \over 4} |\na \na T|^2 +C |\na {\rm Ric}_\o|^3 + C |\na T|^3 +C  \bigg\}.
\eea
Here we used Young's inequality $|\na T| |\na {\rm Ric}_\o|^2 \leq {1 \over 3} |\na T|^3 + {2 \over 3} |\na {\rm Ric}_\o|^3$. In the following, we will use that $\overline{\na} T$ can be expressed as Ricci curvature. We estimate
\bea
&&- 2 \Re \{ F^{i \bar{j}} \na_i |T|^2 \na_{\bar{j}} |\na T|^2 \} \nonumber\\
&\leq& 4 |T| |\na T| (|\na T| + |\overline{\na} T|)(|\na \na T| + |\overline{\na} \na T|) \nonumber\\
&\leq&  4 |T| |\na T|^2 |\na \na T| + 4 |T| |\na T|^2 |\na {\rm Ric}_\o| + 4 |T| |\na T| |{\rm Ric}_\o| |\na \na T|  \nonumber\\
&&+ 4 |T| |\na T| |{\rm Ric}_\o||\na {\rm Ric}_\o| + 4 |T| |\na T| (|\na T| + |\overline{\na} T|) |R \ast T|.
\eea
We may estimate the first term in the following way
\bea
4 |T| |\na T|^2 |\na \na T| \leq 4 |\na T|^2  (\delta_2)^{1/2} |\na \na T| \leq {1 \over 2^3}  |\na T|^4 + 2^5 \delta_2 |\na \na T|^2.
\eea
The other terms may be estimated using Young's inequality, and we can derive
\bea
- 2 \Re \{ F^{i \bar{j}} \na_i |T|^2 \na_{\bar{j}} |\na T|^2 \} \leq {1 \over 2^3}  |\na T|^4 + 2^6 \delta_2 |\na \na T|^2+ C |\na T|^3 + C |\na {\rm Ric}_\o|^3 +C.\nonumber
\eea
Hence
\bea \label{G-second-piece}
\p_t \{ (|T|^2+ \tau_2) |\na T|^2 \} 
&\leq& {1 \over 2 \| \Omega \|} \bigg\{ \Delta_F \{ (|T|^2+ \tau_2) |\na T|^2 \}  - {1 \over 8} |\na T |^4  \nonumber\\
&&- ({\tau_2 \over 4} -2^6 \delta_2)  |\na \na T|^2 +C |\na {\rm Ric}_\o|^3 + C |\na T|^3 +C  \bigg\}.
\eea
Combining (\ref{G-first-piece}) and (\ref{G-second-piece}) gives
\bea 
\p_t G &\leq&  {1 \over 2 \| \Omega \|} \bigg\{ \Delta_F G - \bigg( {1 \over 4}   - 9 \delta_1^2 - 9 \tau_1 \bigg) \alpha'^2|\na {\rm Ric}_\o|^4 \nonumber\\
&& - \bigg( {\tau_1 \over 4} - 2^8 \delta_1^2 \bigg)  (|\na \na {\rm Ric}_\o|^2 + |\na \overline{\na} {\rm Ric}_\o|^2)  - \bigg({\tau_2 \over 4} -2^6 \delta_2 - \delta_1^2 -\tau_1 \bigg) |\na \na T|^2 \nonumber\\
&&- \bigg( {1 \over 8}- 2^4 3^4 \alpha'^{-2} (\delta_1^2+\tau_1)^2-2^2 5^2 \delta_1^2\bigg) |\na T|^4 \nonumber\\
&&+ C_{\alpha',\tau,\delta} |\na {\rm Ric}_\o|^3 + C_{\alpha',\tau,\delta} |\na T|^3+  C_{\alpha',\tau,\delta}  \bigg\}.
\eea
We may choose $\tau_1= \min \{2^{-7}, 2^{-5} 3^{-2} |\alpha'| \}$ and $\tau_2=1$. Then for any $\delta_1, \ \delta_2 >0$ such that
\be
\delta_1, \delta_2 \leq 2^{-6} \tau_1 \ll \tau_2=1,
\ee
we have the estimate
\be
\p_t G \leq  {1 \over 2 \| \Omega \|} \bigg\{ \Delta_F G -  {1 \over 8} \alpha'^2 |\na {\rm Ric}_\o|^4 - {1 \over 16}  |\na T|^4+ C_{\alpha',\tau,\delta}  \bigg\}.
\ee
Now, suppose $G$ attains its maximum at a point $(z,t)$ where $t>0$. From the above estimate, at this point we have
\be
 {1 \over 8} \alpha'^2 |\na {\rm Ric}_\o|^4 +  {1 \over 16} |\na T|^4 \leq C_{\alpha',\tau,\delta}.
\ee
It follows that $G$ is uniformly bounded along the flow, and hence 
\be
|\na {\rm Ric}_\o| \leq C, \ \ |\na T| \leq C,
\ee
along the flow.
\begin{corollary} \label{higher_order_estimates}
There exists $0< \delta_1, \delta_2$ with the following property. Suppose 
\be
-{1 \over 8} g^{p \bar{q}} < \alpha' \| \Omega \|^3 \tilde{\rho}^{p \bar{q}} < {1 \over 8} g^{p \bar{q}},
\ee
\be
|\alpha'{\rm Ric}_\o| \leq \delta_1,
\ee
and 
\be 
|T|^2 \leq \delta_2,
\ee
along the flow. If there exists $\delta_0>0$ such that $0< \delta_0 \leq \| \Omega \| \leq 1$ along the flow, then
\be
|D^k {\rm Ric}_\o| \leq C, \ \ |D^k T| \leq C,
\ee
where $C$ depends only on $\delta_0$, $\delta_1$, $\delta_2$, $\alpha'$, $\rho$, $\tilde{\mu}$, and $(X,\hat{\omega})$.
\end{corollary}
{\it Proof:} Since $\| \Omega \| = e^{-u}$, we are assuming that $|u|$ stays bounded, and that the metrics $\hat{g}$ and $g= e^u \hat{g}$ are equivalent. We are also assuming that $ e^{-u}|Du|_{\hat{g}}^2 \ll 1$ and $e^{-u} |\alpha' u_{\bar{k} j}|_{\hat{g}} \ll 1$. By Theorem \ref{higher_order_est}, there exists $\delta_1$ and $\delta_2$ such that $|\na \na u|$ and $|\na \bar{\na} \na u|$ stay bounded along the flow. We will estimate partial derivatives in coordinate charts. Since
\be
\p_i \bar{\p}_j \p_k u = \na_i \bar{\na}_j \na_k u+ \Gamma^\lambda{}_{ik} u_{\bar{j} \lambda}, \ \ \p_i \p_j u = \na_i \na_j u + \Gamma^\lambda{}_{ij} u_\lambda,
\ee
and the Christoffel symbol
\be
\Gamma^\lambda{}_{ik}  = e^{-u} \hat{g}^{\lambda \bar{\gamma}} \p_i (e^u \hat{g}_{\bar{\gamma} k}) = u_i \delta^\lambda{}_k + \hat{\Gamma}^\lambda{}_{ik}
\ee
stays bounded, we have that 
\be
|u|, \ |\p u|, \ |\p \p u|, \ |\p \bar{\p} u|, \ |\p \bar{\p} \p u| \leq C.
\ee
The scalar equation is
\be \label{scalar_higherorder}
\p_t u = \Delta_{\hat{\o}} u +\alpha' e^{-2u} \tilde{\rho}^{p \bar{q}} u_{\bar{q} p} + \alpha' e^{-u} \hat{\sigma_2}(i \ddb u)  + |Du|^2_{\hat{\o}} + e^{-u} \nu.
\ee
where $\nu(x,u,Du)$. Differentiating once gives
\be \label{linear_parabolic}
\p_t Du = \hat{F}^{p \bar{q}} D u_{\bar{q} p} + \alpha' D(e^{-2u} \tilde{\rho}^{p \bar{q}}) u_{\bar{q} p} + D|Du|^2_{\hat{g}} - \alpha' e^{-u} \hat{\sigma_2}(i \ddb u) Du + D (e^{-u} \nu),
\ee
where
\be \label{F^pq-hat-defn}
\hat{F}^{p \bar{q}} = \hat{g}^{p \bar{q}} + \alpha' e^{-2u} \tilde{\rho}^{p \bar{q}} + \alpha' e^{-u} \hat{\sigma_2}^{p \bar{q}}.
\ee
We note that $\hat{F}^{j \bar{k}}$ only differs from $F^{j \bar{k}}$ (\ref{F^pq-defn}) by a factor of $e^u$. From our assumptions on $|\alpha' {\rm Ric}_\o| = e^{-u} | \alpha' \p \bar{\p} u|_{\hat{g}}$ and $\| \Omega \| = e^{-u}$, we have uniform ellipticity of $\hat{F}^{j \bar{k}}$. Differentiating twice yields
\be
\p_t u_{\bar{k} j} = \hat{F}^{p \bar{q}} \p_p \p_{\bar{q}} u_{\bar{k} j} + \Psi(x,u,\p u, \p \p u, \p \bar{\p} u, \p \bar{\p} \p u),
\ee
where $\Psi$ is uniformly bounded along the flow. By the Krylov-Safonov theorem \cite{KrSa}, we have that $u_{\bar{k} j}$ is bounded in the $C^{\alpha/2,\alpha}$ norm. The function $u$ and the spacial gradient $Du$ are also bounded in the $C^{\alpha/2,\alpha}$ norm since the right-hand sides of (\ref{scalar_higherorder}) and (\ref{linear_parabolic}) are bounded. We may now apply parabolic Schauder theory (for example, in \cite{Kr}) to the linearized equation (\ref{linear_parabolic}). Standard theory and a bootstrap argument give higher order estimates of $u$, and hence we obtain estimates on derivatives of the curvature and torsion of $g=e^u \hat{g}$.

\section{Long time existence}
\setcounter{equation}{0}
\begin{proposition} \label{all_estimates}
Let $C_1$, $C_2$, $C_5$ be the named constants as before, which only depend on $(X,\hat{g})$, $\tilde{\mu}$, $\rho$, and $\alpha'$. There exists $M_0 \gg 1$ such that for all $M \geq M_0$, the following statement holds. If the flow exists on $[0,t_0)$, and initially starts with $u_0=\log M$, then along the flow
\be
{1 \over C_1 M} \leq e^{-u} \leq {C_2 \over M}, \ \ |T|^2 \leq {C_3 \over M}, \ \ |\alpha' {\rm Ric}_\o| \leq {C_5 \over M^{1/2}},
\ee
and
\be
|D^k u|^2_{\hat{g}} \leq \tilde{C}_k, \ \ {1 \over 2} \hat{g}^{j \bar{k}} \leq \hat{F}^{j \bar{k}} \leq 2 \hat{g}^{j \bar{k}},
\ee
where $\tilde{C}_k$ only depends on $(X,\hat{g})$, $\mu$, $\rho$, $\alpha'$, $M$.
\end{proposition}
{\it Proof:} Let $\delta_1$ and $\delta_2$ be the constants from Corollary \ref{higher_order_estimates}, and choose a smaller $\delta_1$ if necessary to ensure $\delta_1 < 10^{-6}$. Recall that from Theorem \ref{C0-thm},
\be
{1 \over C_1 M} \leq \| \Omega \| = e^{-u} \leq {C_2 \over M}
\ee
along the flow for $M$ large enough. Consider the set
\be
I = \{ t \in [0,t_0) \ {\rm such} \ {\rm that} \ |\alpha' {\rm Ric}_\o| \leq \delta_1, \   |T|^2 \leq \delta_2 \ {\rm holds} \ {\rm on} \ [0,t] \}.
\ee
Since at $t=0$ we have $|\alpha' {\rm Ric}_\o| = |T|^2=0$, we know that $I$ is non-empty. By definition, $I$ is relatively closed. We now show that $I$ is open. Suppose $\hat{t} \in I$. By definition of $I$, the hypothesis of Theorem \ref{torsion_estimate} is satisfied, hence $|T|^2 \leq {C_3 \over M} < \delta_2$ at $\hat{t}$ as long as $M$ is large enough. It follows that the hypothesis of Theorem \ref{ricci_estimate} is satisfied as long as $M$ is large enough, hence $|\alpha' {\rm Ric}_\o| \leq {C_5 \over M^{1/2}} <\delta_1$ at $\hat{t}$. We can conclude the existence of $\e>0$ such that $[\hat{t}+\e) \subset I$, and hence $I$ is open.
\par It follows that $I=[0,t_0)$. We know that $-C \hat{g}^{p \bar{q}} \leq \tilde{\rho}^{p \bar{q}} \leq C \hat{g}^{p \bar{q}}$ since $\tilde{\rho}$ can be bounded using the background metric. For $M$ large enough, we can conclude
\be
-{1 \over 8} e^{-u} \hat{g}^{p \bar{q}} < \alpha' e^{-3u} \tilde{\rho}^{p \bar{q}} < {1 \over 8} e^{-u} \hat{g}^{p \bar{q}},
\ee
and we can apply Corollary \ref{higher_order_estimates} to obtain higher order estimates of $u$. Uniform ellipticity follows from the definition of $\hat{F}^{j \bar{k}}$ (\ref{F^pq-hat-defn}) and the estimates on $|\alpha' {\rm Ric}_\o| = e^{-u} | \alpha' \p \bar{\p} u|_{\hat{g}}$ and $\| \Omega \|$. Q.E.D.

\begin{theorem} \label{long_time}
There exists $M_0 \gg 1$ such that for all $M \geq M_0$, if the flow initially starts with $u_0 = \log M$, then the flow exists on $[0,\infty)$.
\end{theorem}
{\it Proof:} By short-time existence \cite{PPZ}, we know the flow exists for some maximal time interval $[0,T)$. If $T<\infty$, we may apply the previous proposition to extend the flow to $[0,T+\epsilon)$, which is a contradiction. Q.E.D.

\section{Convergence of the flow}
\setcounter{equation}{0}

We may apply Theorem \ref{long_time} to construct solutions to the Fu-Yau equation.
\begin{theorem}
There exists $M_0 \gg 1$ such that for all $M \geq M_0$, if the flow initially starts with $u_0 = \log M$, then the flow exists on $[0,\infty)$ and converges smoothly to a function $u_\infty$, where $u_\infty$ solves
\be \label{Fu-Yau-eqn}
0 = i \ddb (e^{u_\infty} \hat{\o} - \alpha' e^{-u_\infty} \rho) + {\alpha' \over 2} i \ddb u_\infty \wedge i \ddb u_\infty + \mu , \ \ \ \int_X e^{u_\infty} = M .
\ee
\end{theorem}
{\it Proof:} Since we will work with the scalar equation, all norms in this section will be with respect to the background metric $\hat{\o}$. Let $v = \p_t e^u$. Recall that
\be
\int_X v = 0,
\ee
along the flow. Differentiating equation (\ref{flow_forms}) with respect to time gives
\be
2 \p_t v {\hat{\o}^2 \over 2!} = i \ddb (v \hat{\o} + \alpha' e^{-2u} v \rho) + \alpha' i \ddb u \wedge i \ddb (e^{-u}v).
\ee
Consider the functional
\be
J(t) =  \int_X v^2 \,  {\hat{\o}^2 \over 2!}.
\ee
Compute
\bea
{d J \over dt} &=& \int_X v \, i \ddb (v \hat{\o} + \alpha' e^{-2u} v \rho) + \alpha' \int_X v \,i \ddb u \wedge i \ddb (e^{-u}v)  \\
&=& - \int_X i \p v \wedge \bar{\p} v \wedge \hat{\o}- \alpha' \int_X i \p v \wedge \bar{\p} (e^{-2u} v \rho) - \alpha' \int_X i \ddb u \wedge i \p v \wedge i \bar{\p} (e^{-u} v) \nonumber\\
&=& - \int_X |\na v|^2 - \alpha' \int_X e^{-2u} \, i \p v \wedge \bar{\p}  v \wedge \rho +2 \alpha' \int_X e^{-2u}v \, i \p v \wedge \bar{\p} u \wedge  \rho   \nonumber\\
&&  - \alpha' \int_X e^{-2u}v \, i \p v \wedge \bar{\p}  \rho - \alpha' \int_X e^{-u} \, i \ddb u \wedge i \p v \wedge i \bar{\p} v  + \alpha' \int_X e^{-u}v \,i \ddb u \wedge i \p v \wedge i \bar{\p} u. \nonumber
\eea
We may estimate
\bea
{d J \over dt} &\leq&  - \int_X |\na v|^2 + \alpha' \|\rho\| \int_X  e^{-2u} |\na v|^2 +2 \alpha' \| \rho \| \|\na u \| \int_X e^{-2u} |v| \, |\na v| \\
&&  + \alpha' \| \p \rho \| \int_X e^{-2u} |v| \, |\na v| + \| \alpha' e^{-u} i \ddb u \| \int_X |\na v|^2 \nonumber\\
&&  +  \|\na u\| \, \|\alpha' e^{-u} i \ddb u \| \int_X |v| \, |\na v|.  \nonumber
\eea
By Proposition \ref{all_estimates}, we know that on $[0,\infty)$ we have the estimates
\be
e^{-u} \leq {C_2 \over M} \ll 1, \ \ |\na u|^2_{\hat{g}} \leq C_3 C_1, \ \ | \alpha' e^{-u} u_{\bar{k} j} |_{\hat{g}} \leq {C_5 \over M^{1/2}}.
\ee
Hence for any $\e >0$, we can choose $M$ large enough such that
\bea
{d J \over dt} \leq - {1 \over 2} \int_X |\na v|^2 +  \e \int_X |v| \, |\na v| \leq  - \bigg({1 \over 2} -  {\e \over 2} \bigg) \int_X |\na v|^2 + {\e \over 2} \int_X |v|^2. 
\eea
Since $\int_X v = 0$, we may use the Poincar\'e inequality to obtain, for $\e>0$ small enough,
\bea\label{Jinequality}
{d J \over dt} \leq - \eta \int_X v^2 = - \eta J,
\eea
with $\eta>0$. This implies that
\bea
J(t) \leq C e^{-\eta t}, 
\eea 
that is, 
\bea \label{v-L2-decay}
\int_X v^2 \leq Ce^{-\eta t}.
\eea
From this estimate, we see that for any sequence $v(t_j)$ converging to $v_{\infty}$, we have $v_{\infty} =0$.

\

We can now show convergence of the flow. Following the argument given in Proposition 2.2 in \cite{Cao}, we have
\bea
\int_X | e^u(x, s') - e^u(x, s)| &\leq & \int_X \int^{s'}_{s} | \p_t e^u(x, t)| = \int^{s'}_{s}\int_X  |v(x, t)| \nonumber\\
&\leq & 
 \int^{s'}_{s} \left(\int_X v^2 \right)^{1\over 2} dt
\leq \int^{+\infty}_{s} \left(\int_X v^2 \right)^{1\over 2} dt
\nonumber\\
&\leq & 
C \int^{+\infty}_{s} e^{- {\eta \over 2} t} dt
\eea
Recall that we normalized the background metric such that $\int_X {\hat{\o}^2 \over 2} =1$. This estimate shows that, as $t\rightarrow +\infty$, $e^u(x, t)$ are Cauchy in $L^1$ norm. Thus $e^u(x, t)$ converges in the $L^1$ norm to some function $e^{u_\infty}(x)$ as $t\rightarrow \infty$.
\medskip
\par By our uniform estimates, $e^{u_\infty}$ is bounded in $C^\infty$, and a standard argument shows that $e^u$ converges in $C^\infty$. Indeed, if there exist a sequence of times such that $\| e^{-u(x,t_j)}- e^{-u_\infty(x)} \|_{C^k} \geq \epsilon$, then by our estimates a subsequence converges in $C^k$ to $e^{-u'_\infty}$. Then $\| e^{-u'_\infty(x)}- e^{-u_\infty(x)} \|_{L^1}=0$ but $\| e^{-u'_\infty(x)}- e^{-u_\infty(x)} \|_{C^k} \geq \epsilon$, a contradiction.
\medskip
\par It follows from (\ref{v-L2-decay}) that $e^{u_\infty}$ satisfies the Fu-Yau equation (\ref{Fu-Yau-eqn}).

\

\bigskip
Department of Mathematics, Columbia University, New York, NY 10027, USA

\smallskip

phong@math.columbia.edu

\bigskip
Department of Mathematics, Columbia University, New York, NY 10027, USA

\smallskip
picard@math.columbia.edu

\bigskip
Department of Mathematics, University of California, Irvine, CA 92697, USA

\smallskip
xiangwen@math.uci.edu

\end{document}